\documentclass{amsart}

\usepackage{amsmath,amssymb,amsthm,mathtools}
\usepackage{xcolor}
\usepackage{graphicx}
\usepackage{booktabs}
\usepackage{float}
\usepackage{placeins}
\usepackage{mathrsfs}
\usepackage[colorlinks=true,linkcolor=blue,citecolor=blue,urlcolor=blue]{hyperref}

\theoremstyle{plain}

\newtheorem{theorem}{Theorem}

\newtheorem{lemma}{Lemma}[section]
\newtheorem{proposition}[lemma]{Proposition}
\newtheorem{corollary}[lemma]{Corollary}
\newtheorem{theoremcorollary}[lemma]{Corollary}

\theoremstyle{definition}
\newtheorem{assumption}[lemma]{Assumption}
\newtheorem{setting}[lemma]{Setting}

\newtheorem{example}[lemma]{Example}

\theoremstyle{remark}
\newtheorem{remark}[lemma]{Remark}

\numberwithin{equation}{section}

\newcommand{\R}{\mathbb R}
\newcommand{\N}{\mathbb N}
\newcommand{\QT}{\mathcal O_T}
\newcommand{\OT}{\overline{\mathcal O_T}}

\newcommand{\Hc}{\mathcal H}
\newcommand{\Ec}{\mathcal E}
\newcommand{\Xc}{\mathcal X}
\newcommand{\ud}{\,\mathrm d}
\newcommand{\U}{\Upsilon}
\newcommand{\Un}{\Upsilon_n}
\newcommand{\J}{\mathcal J}

\newcommand{\Jny}{\mathcal J_y^n}
\newcommand{\tnorm}[1]{\lvert\!\lvert\!\lvert #1\rvert\!\rvert\!\rvert}
\newcommand{\Kone}{\hyperref[cond:K1]{(K1)}}
\newcommand{\Ktwo}{\hyperref[cond:K2]{(K2)}}
\newcommand{\Kthree}{\hyperref[cond:K3]{(K3)}}
\newcommand{\Kfour}{\hyperref[cond:K4]{(K4)}}
\newcommand{\Kfive}{\hyperref[cond:K5]{(K5)}}

\newcommand{\Konefour}{\hyperref[cond:K1]{(K1)}--\hyperref[cond:K4]{(K4)}}

\begin{document}

\title[Scheme-Induced Minimum Action Methods for SPDEs]{Scheme-Induced Minimum Action Methods for SPDEs: A Variational Framework for Convergence Analysis}

\author{Jialin Hong}
\address{State Key Laboratory of Mathematical Sciences, Academy of Mathematics and Systems Science, Chinese Academy of
	Sciences; School of Mathematical Sciences, University of Chinese
Academy of Sciences, Beijing 100049, China}
\email{hjl@lsec.cc.ac.cn}
\thanks{This work is supported by the National Key R\&D Program of China under Grant No.\ 2020YFA0713701,  the National Natural Science Foundation of China (Nos.\ 12031020, 12461160278, 12471386, 12671472, 12201228), and the Fundamental Research Funds for the Central Universities  (No.\ 2026BRSXB006).}

\author{Diancong Jin}
\address{School of Mathematics and Statistics, Huazhong University of Science and
Technology, Wuhan 430074, China; Hubei Key Laboratory of Engineering Modeling
and Scientific Computing, Huazhong University of Science and Technology, Wuhan
430074, China}
\email{jindc@hust.edu.cn (Corresponding author)}

\author{Derui Sheng}
\address{Department of Applied Mathematics, The Hong Kong Polytechnic University,
Hung Hom, Kowloon 999077, Hong Kong}
\email{sdr@lsec.cc.ac.cn}

%60F10：Large deviations；
%49J45：Methods involving semicontinuity and convergence; relaxation
%60H15：Stochastic partial differential equations；
%65M06：Finite difference methods for time-dependent PDEs；
%65M12：Stability and convergence of numerical methods for time-dependent PDEs。

\begin{abstract}
The minimum action method (MAM) is an important numerical
tool for computing the most probable transition paths and the leading exponential rates of rare-event probabilities for stochastic systems with small noise. To the best of our knowledge, no rigorous convergence analysis is available for fully discrete MAMs for stochastic partial differential equations (SPDEs). This paper establishes the convergence of the
constrained minima of fully discrete scheme-induced action functionals for terminal observations of scalar semilinear SPDEs. A main difficulty is that a consistent discrete approximation of an admissible control generally fails to satisfy the discrete terminal constraint exactly. To overcome this difficulty, we construct an asymptotically negligible correction that restores the constraint, and incorporate this argument into an abstract variational framework. We then propose a  criterion that reduces the assumptions of this framework to verifiable conditions involving the continuous and discrete Green functions. Applying this criterion to fully discrete schemes for a stochastic wave equation and a fourth-order parabolic equation, we prove convergence of the constrained minima along arbitrary spatial and temporal refinement sequences, without a mesh-ratio condition.
\end{abstract}

\maketitle

\section{Introduction}

\subsection{Research background}
Small random perturbations can induce qualitative changes in the long-time
behavior of a dynamical system.  While a deterministic trajectory may remain near a stable state, small noise can drive the system out of its basin
of attraction and toward another stable state.  Such noise-induced transitions
arise in a variety of applications, including nucleation, chemical reactions,
and transitions between distinct physical regimes \cite{FW,Kampen81}.  In the
small-noise regime, these transition events are rare but could have a decisive influence on
the long-time dynamics.

In this paper, we study small-noise-induced transitions for the following scalar semilinear
stochastic partial differential equation (SPDE):
\begin{align}\label{SPDE}
\mathcal L u_\varepsilon(t,x)
=b(u_\varepsilon(t,x))
+\sqrt{\varepsilon}\,\sigma(u_\varepsilon(t,x))\dot W(t,x),~(t,x)\in \mathcal O_T=(0,T)\times\mathcal O, 
\end{align}
for a given spatial domain $\mathcal O$, where $W$ is a Brownian sheet
on $\mathcal O_T$ and the noise intensity $\varepsilon$ is sufficiently small. Here $\mathcal L$ denotes a linear spatiotemporal differential operator such as  $\partial_t-\Delta$,
$\partial_{tt}-\Delta$, and $\partial_t+\Delta^2$. 
For a fixed observation point $x_0\in\mathcal O$ and a Borel
measurable set $A\subset\mathbb R$, we are interested in the transition event that  $u_\varepsilon$ reaches $A$ at $(T,x_0)$ starting from a given initial state.

The  Freidlin--Wentzell theory of large deviations provides an effective framework for
quantifying the probability $\mathbb P\bigl(u_\varepsilon(T,x_0)\in A\bigr)$ (cf.\ \cite{Dembo,FW}). When the large
deviation upper and lower bounds agree for $A$, this theory yields
\[
\lim_{\varepsilon\downarrow0}\varepsilon
\log\mathbb P\bigl(u_\varepsilon(T,x_0)\in A\bigr)
=-\inf_{y\in A}I(y).
\]
The variational representation of $I$ is formulated through a controlled
deterministic equation, also called the skeleton equation.  Formally, replacing $\sqrt{\varepsilon}\dot W$ by a control
$h\in L^2(\mathcal O_T)$ defines the solution  $\Upsilon(h)$ of the skeleton equation, whose
Green-function representation is
\begin{equation*}%\label{eq:skeleton}
\Upsilon(h)(t,x)
=g(t,x)+\int_0^t\int_{\mathcal O}G_{t-s}(x,z)
[b(\Upsilon(h)(s,z))
+\sigma(\Upsilon(h)(s,z))h(s,z)] \ud z\ud s. 
\end{equation*}
Here $g$ is the solution of the homogeneous linear problem and $G$ is the corresponding
Green function.  The terminal rate function then admits the variational
representation
\begin{equation}\label{eq:I-intro}
I(y)=\inf_{\{h\in L^2(\mathcal O_T):
		\Upsilon(h)(T,x_0)=y\}}
\frac12\|h\|_{L^2(\mathcal O_T)}^2.
\end{equation}
The precise assumptions and a detailed derivation of the formulation of $I$ are given in
Section~\ref{sec:spde-notation}.
Whenever the infimum is attained, we call any
\[
h_y^*\in
\operatorname*{arg\,min}_{\{h\in L^2(\mathcal O_T):
		\Upsilon(h)(T,x_0)=y\}}
\frac12\|h\|_{L^2(\mathcal O_T)}^2
\]
an optimal control and set $u_y^*:=\Upsilon(h_y^*)$.  The path $u_y^*$ is called a
minimum action path associated with the terminal value $y$.  Thus, at the
large-deviation scale, $I(y)$ quantifies the difficulty of $u_{\varepsilon}$ reaching $y$ at $(T,x_0)$ starting from an initial state,
whereas a minimizer $u_y^*$ describes a most likely realization of the corresponding
rare event.  This variational characterization forms the basis of minimum
action methods (MAMs), which approximate rare transitions through deterministic
constrained optimization; see, e.g., \cite{ERenVandenEijnden2004,LDPofonepoint}.

Numerous numerical methods have been proposed for minimum action problems; below, we mention a few key approaches without attempting an exhaustive review.
E, Ren, and Vanden-Eijnden introduced the original MAM
\cite{ERenVandenEijnden2004}, followed by several variants designed to address
different numerical difficulties in computing transition paths.  The geometric
MAM eliminates the dependence on the time parametrization of the path
\cite{HeymannVandenEijnden2008}, while the adaptive MAM redistributes temporal
grid points to better resolve portions of the path with rapid motion
\cite{ZhouRenE2008}.  The MAM with optimal linear time scaling determines the transition time optimally from the transition path through a linear rescaling of time \cite{Wan2015} and was subsequently extended to dynamical systems with constant time delays by introducing an auxiliary path \cite{WanZhai2021}.
\subsection{Motivation and objective}
Despite substantial progress in algorithm design, rigorous convergence results
for MAMs remain limited. We are only aware of \cite{WanYuZhai2018,LDPofonepoint} concerning stochastic ordinary differential
equations. In \cite{WanYuZhai2018}, a conforming finite element approximation
of the Freidlin--Wentzell action functional with additive noise was analyzed by
means of $\Gamma$-convergence. Our earlier work \cite{LDPofonepoint} proved
convergence of minimizers for finite-difference discretizations of the action
functional and derived convergence rates for the minimum values.  
However, to the best of our knowledge, no rigorous convergence analysis is available for fully discrete MAMs for stochastic partial differential equations (SPDEs).

Most existing MAMs construct a numerical action problem by directly discretizing the continuous action functional; see, e.g., \cite{ERenVandenEijnden2004,HeymannVandenEijnden2008,Wan2015,ZhouRenE2008}. In this work, we study the convergence of fully discrete \emph{scheme-induced minimum action methods}.   A closely related construction was introduced by Wan and Yu under the name dynamic-solver-consistent MAM \cite{WanYu2017}. They first discretized the SPDE in space and then defined the action functional as the Freidlin--Wentzell rate function of the resulting finite-dimensional stochastic system. Their construction is therefore scheme-induced at the spatially semidiscrete level.  The purpose of the present paper is not merely to extend this construction from a semidiscrete to a fully discrete setting. Instead, we aim to develop a variational framework for proving the  convergence of fully discrete MAMs for \eqref{SPDE}.

Following this scheme-induced principle, we first discretize the SPDE \eqref{SPDE} and then use the rate function induced by a fully discrete numerical solution to approximate \(I\). More precisely, for each fixed grid \((n,m_n)\), let $u_\varepsilon^{n,m_n}$ denote the corresponding fully discrete approximation of $u^\varepsilon$, where $m_n$ is the associated number of time steps satisfying $m_n\to\infty$ as $n\to\infty$. We
consider numerical schemes for \eqref{SPDE} such that, for each fixed $n$, $\{u_\varepsilon^{n,m_n}(T,x_0)\}_{\varepsilon>0}$ satisfies a large
deviation principle (LDP) as $\varepsilon\to0$ with the rate function
\begin{equation}\label{eq:Inmn}
I^{n,m_n}(y)=
\inf_{\{h_n\in\mathcal X_n:
\Upsilon_n(h_n)(T,x_0)=y\}}
\frac12\|h_n\|_{L^2(\mathcal O_T)}^2,
\qquad y\in\mathbb R.
\end{equation}
Here, $\mathcal X_n\subset L^2(\mathcal O_T)$ denotes the corresponding finite-dimensional discrete control space, and $\Upsilon_n(h_n)$ denotes the solution of the discrete skeleton equation driven by the control $h_n$.
The detailed formulation of $I^{n,m_n}$ is given in
Section~\ref{sec:spde-notation}.

%
% In detail, for each fixed mesh, let $u_\varepsilon^{n,m_n}$ denote the
%resulting fully discrete solution of $u_\varepsilon$.  Here \(n\) indexes the spatial mesh, \(m_n\) is the number of time steps associated with that mesh satisfying $\lim_{n\to \infty}m_n=\infty$, and \(\tau_n=T/m_n\) is the time-step size.  We
%consider numerical schemes of \eqref{SPDE} for which $u_\varepsilon^{n,m_n}(T,x_0)$ satisfies a large
%deviation principle (LDP) with rate function
%\[
%I^{n,m_n}(y)=
%\inf_{\{h_n\in\mathcal X_n:
%\Upsilon_n(h_n)(T,x_0)=y\}}
%\frac12\|h_n\|_{L^2(\mathcal O_T)}^2,
%\qquad y\in\mathbb R.
%\]
%Here, $\Upsilon_n(h_n)$ denote the
%corresponding solution of the discrete skeleton equation driven by a discrete control $h_n$.
%Here $\mathcal X_n\subset L^2(\mathcal O_T)$ is the corresponding finite-dimensional discrete control
%space.
%The detailed formulation of $I^{n,m_n}$ is given in
%Section~\ref{sec:spde-notation}.

%This construction has two important features. First, \(I^{n,m_n}\) is not obtained by discretizing the continuous constrained variational problem \eqref{eq:I-intro}; rather, it is the large deviation rate function governing the rare-event probabilities of the fully discrete numerical solution \(u^{n,m_n}_\varepsilon(T,x_0)\). 
This construction has two important features. First, the discrete variational problem \eqref{eq:Inmn} is induced by the numerical scheme of \eqref{SPDE} itself rather than through a direct discretization of the continuous constrained problem \eqref{eq:I-intro}. The rate function $I^{n,m_n}$ itself governs the rare-event probabilities of the fully discrete numerical solution \(u^{n,m_n}_\varepsilon\).
Second, the discrete control spaces \(\mathcal X_n\) are naturally embedded in the control space \(L^2(\mathcal O_T)\), so the rate functions $I$ and $I^{n,m_n}$ share the same control-based variational structure, which facilitates their comparison. It is worth emphasizing that the convergence of the numerical solution \(u^{n,m_n}_\varepsilon(T,x_0)\) alone does not imply the convergence of the constrained minima \(I^{n,m_n}\). 
%The difficulty is that the control spaces $\mathcal{X}_n$ vary with the mesh, while the terminal equality constraint must be satisfied exactly at every discrete level. 
In this work,
our main objective is to prove the convergence
\[
I^{n,m_n}(y)\rightarrow I(y),
\qquad y\in\mathbb R,
\]
along an arbitrary diagonal refinement sequence $\{(n,m_n)\}_{n\ge1}$ with
$m_n\to\infty$. 

\subsection{Main results and organization}
Our convergence analysis of $I^{n,m_n}$ is based on the theory of $\Gamma$-convergence and involves two main difficulties. First, the controls $h_n$
belong to mesh-dependent spaces $\mathcal X_n$, and the boundedness in $L^2(\mathcal O_T)$-norm generally yields,
after taking a subsequence, only $h_n\rightharpoonup h$ in $L^2(\mathcal O_T)$.
To match the terminal constraint $\Upsilon(h)(T,x_0)=y$, one must nevertheless prove
\begin{align} \label{eq:weaktostrong}
h_n\rightharpoonup h
\quad\Longrightarrow\quad
\Upsilon_n(h_n)\to\Upsilon(h)
\quad\text{in }C(\overline{\mathcal O_T}).
\end{align}
Second, if $h$ is an admissible continuous control, its projection $Q_nh$ onto $\mathcal{X}_n$
usually satisfies only
\[
\Upsilon_n(Q_nh)(T,x_0)\rightarrow \Upsilon(h)(T,x_0)=y,
\]
rather than the exact discrete constraint $\Upsilon_n(Q_nh)(T,x_0)=y$.
Here $Q_n:L^2(\mathcal O_T)\to\mathcal X_n$ is the $L^2$-orthogonal projection
obtained by averaging the control over each space-time grid cell.

In Section \ref{sec:kernel-criterion}, we  resolve these difficulties through an abstract variational argument, which is  formulated under four abstract assumptions concerning the qualitative properties of the skeleton maps $\Upsilon_n$ and $\Upsilon$ (see Assumptions \ref{ass:abstract-compact-stability}--\ref{ass:abstract-terminal-nondegeneracy}). 
These assumptions include the compactness and weak stability of the discrete skeleton maps $\Upsilon_n$, approximation of controls, directional equicontinuity of $\Upsilon_n$, and  local nondegeneracy of the terminal map $h\mapsto \Upsilon(h)(T,x_0)$.
 In particular, 
the compactness and weak stability of $\Upsilon_n$
allow us to obtain the property \eqref{eq:weaktostrong}. 
To deal with the issue \(\Upsilon_n(Q_nh)(T,x_0)\neq \Upsilon(h)(T,x_0)\), we construct an admissible recovery sequence by adding an asymptotically negligible correction $r_n\in\mathcal X_n$ to the projected control \(Q_nh\) such that (see Lemma~\ref{lem:abstract-direction-correction})
%
%More precisely, under Assumptions \ref{ass:abstract-compact-stability}--\ref{ass:abstract-terminal-nondegeneracy}, Lemma~\ref{lem:abstract-direction-correction} yields corrections \(r_n\in\mathcal X_n\) such that
\[
\Upsilon_n(Q_nh+r_n)(T,x_0)=\Upsilon(h)(T,x_0),
\qquad
\|r_n\|_{L^2(\mathcal O_T)}\to0.
\]
The variational convergence framework in Theorem~\ref{thm:abstract-framework} is independent of any particular SPDE or discretization and can therefore be applied once the corresponding abstract assumptions have been verified. Moreover, Corollary~\ref{cor:subsequential-minimizers} yields subsequential convergence of asymptotically minimizing controls $h_n$ and their associated controlled paths $u_n=\Upsilon_n(h_n)$ for the discrete problem \eqref{eq:Inmn}.
% implies that, for every target $y$
%with $I(y)<+\infty$, any sequence of asymptotically minimizing controls admits
%a subsequence converging strongly in $L^2(\mathcal O_T)$, while the associated
%controlled paths converge in $C(\overline{\mathcal O_T})$ to a minimum action
%path of the continuous constrained problem.

In Section~\ref{S:KB}, we apply the abstract convergence framework of Theorem~\ref{thm:abstract-framework} to scheme-induced MAMs for the SPDE \eqref{SPDE}. Motivated by the mild formulations of the skeleton maps $\Upsilon_n$ and $\Upsilon$, we derive a more directly verifiable criterion that reduces the compactness and weak stability of $\Upsilon_n$ to conditions \Kone{}--\Kfour{} on the continuous and discrete Green functions; see Proposition~\ref{prop:kernel-compactness-stability}.
A more delicate step is to verify the local nondegeneracy of the terminal map
$
h\mapsto \Upsilon(h)(T,x_0).
$
To this end, we first establish resolvent estimates for the forward linearized equation and the backward Volterra equation in Lemmas~\ref{lem:linear-resolvent} and \ref{lem:backward-volterra-wellposed}, respectively. We then prove the $C^1$-regularity of the skeleton map $\Upsilon$ in Lemma~\ref{lem:skeleton-C1} and derive an adjoint representation of $D\Upsilon(\cdot)(T,x_0)[q]$ in Lemma~\ref{lem:adjoint-representation}. Combining this representation with the nondegeneracy condition \Kfive{} for the continuous Green function, Proposition~\ref{prop:terminal-nondegeneracy} establishes the local nondegeneracy of the terminal map. 
 Building on Propositions~\ref{prop:kernel-compactness-stability} and \ref{prop:terminal-nondegeneracy}, Theorem~\ref{thm:closed-diagonal} provides a criterion that reduces the abstract assumptions of Theorem~\ref{thm:abstract-framework} to conditions \Kone{}--\Kfive{} on the continuous and discrete Green functions.
 %
%
%In Section~\ref{S:KB}, we apply the abstract convergence framework of Theorem~\ref{thm:abstract-framework} to scheme-induced MAMs for the SPDE \eqref{SPDE}. Motivated by the mild formulations of the skeleton maps $\Upsilon_n$ and $\Upsilon$, we derive a more directly verifiable criterion that reduces the compactness and weak stability of $\Upsilon_n$ to conditions on the continuous and discrete Green functions, as shown in Proposition \ref{prop:kernel-compactness-stability}.
%Another delicate step is to verify
%the local nondegeneracy of the terminal map $h\mapsto \Upsilon(h)(T,x_0)$. In Proposition~\ref{prop:terminal-nondegeneracy}.  We first establish resolvent estimates for the
%forward linearized equation and the backward Volterra equation in
%Lemmas~\ref{lem:linear-resolvent} and
%\ref{lem:backward-volterra-wellposed}.  We then prove the $C^1$-regularity of
%the skeleton map in Lemma~\ref{lem:skeleton-C1} and obtain the adjoint
%representation of $D\Upsilon(h)(T,x_0)[q]$
%in Lemma~\ref{lem:adjoint-representation}.  Combining this representation with the nondegeneracy condition \Kfive{} for the continuous Green function, we prove terminal nondegeneracy and establish the key ingredient needed for
%Theorem~\ref{thm:closed-diagonal}.

In Section~\ref{sec:applications}, we apply the criterion to two fully discrete schemes. For the one-dimensional stochastic wave equation, we use a spatial finite difference method (FDM) combined with an exponential Euler method (EEM) in time, while for the fourth-order parabolic equation in dimensions $d\in\{1,2,3\}$, we use a spatial FDM combined with the implicit Euler method. By verifying conditions \Kone{}--\Kfive{}, Corollaries~\ref{cor:wave-diagonal-convergence} and \ref{cor:fourth-diagonal-convergence} establish the convergence of the corresponding scheme-induced MAMs for every sequence $m_n\to\infty$, without imposing any prescribed ratio between the spatial and temporal mesh sizes. The proofs of these corollaries rely on the verification of the Green-function conditions and the fixed-grid LDPs, which are provided in Appendices~\ref{app:kernel-verifications} and \ref{app:fixed-grid-LDPs}, respectively.

In Section~\ref{sec:numerical-experiments}, we present the numerical experiments and
 illustrate the convergence of the scheme-induced MAMs through
an exactly solvable linear example and a nonlinear example for the stochastic
wave equation. Finally, concluding remarks
are given in Section~\ref{sec:concluding-remarks}.

\subsection{Contributions} The main contributions of this paper can be summarized as follows.
First, Theorem~\ref{thm:abstract-framework} establishes a variational
convergence framework for scheme-induced MAMs on
changing control spaces,  with the terminal-correction argument of Lemma~\ref{lem:abstract-direction-correction} providing the key tool for enforcing the discrete terminal constraint exactly.
Second, Theorem~\ref{thm:closed-diagonal} reduces the abstract framework in Theorem \ref{thm:abstract-framework} 
to  verifiable Green-function conditions.  Finally,
Corollaries~\ref{cor:wave-diagonal-convergence} and
\ref{cor:fourth-diagonal-convergence} establish fully discrete diagonal
convergence for the stochastic wave and fourth-order parabolic equations without imposing
a space-time mesh-ratio condition.

A preliminary version of this work, posted as arXiv:2209.08341v1, treated a spatial finite difference semidiscretization of the one-dimensional stochastic wave equation. The present manuscript supersedes that version, which will not be submitted separately, and replaces its scheme-specific inverse construction by the general terminal-correction argument developed here.

%
%\emph{Organization.} The remainder of the paper is organized as follows.
%Section~\ref{sec:kernel-criterion} introduces the scalar SPDE setting and
%develops the abstract variational framework. Section~\ref{S:KB} establishes
%the Green-function criterion, and Section~\ref{sec:applications} applies the
%criterion to the two fully discrete schemes. Numerical experiments are
%presented in Section~\ref{sec:numerical-experiments}, and concluding remarks
%are given in Section~\ref{sec:concluding-remarks}.
%Appendices~\ref{app:kernel-verifications} and
%\ref{app:fixed-grid-LDPs} contain, respectively, the verification of the
%Green-function conditions and the fixed-grid LDPs.

\section{\texorpdfstring{{Scalar SPDE setting and abstract variational framework}}{Scalar SPDE setting and abstract variational framework}}
\label{sec:kernel-criterion}

In this part, we first present the continuous and discrete skeleton equations together with their Green-function representations in Section \ref{sec:spde-notation}. Then we  develop the abstract variational framework for the convergence analysis of scheme-induced MAMs in Section \ref{sec:abstract-framework}.

\subsection{\texorpdfstring{{A motivating scalar SPDE setting and
scheme-induced rate functions}}{A motivating scalar SPDE setting and
scheme-induced rate functions}}
\label{sec:spde-notation}

We now introduce the continuous and fully discrete skeleton equations and fix the notation used below. At this stage, we assume that these equations are well posed and that the corresponding continuous and fixed-grid LDPs hold; Section \ref{sec:applications} verifies these properties for the schemes considered here.

{ Let $\mathcal O\subset\mathbb R^d$ be a bounded domain satisfying
$|\partial\mathcal O|=0$, and set
$$
\QT:=(0,T)\times\mathcal O,
\qquad
\OT:=[0,T]\times\overline{\mathcal O}.
$$
All $L^2$-spaces below are taken over $\QT$, while continuity, supremum norms,
and compactness statements are understood on $\OT$.  Consider the
SPDE on $\QT$ with small noise:}
\begin{equation}\label{eq:abstract-spde}
\mathcal{L} u_\varepsilon(t,x)=b(u_\varepsilon(t,x))+\sqrt{\varepsilon}\sigma(u_\varepsilon(t,x))\dot{W}(t,x),
\end{equation}
with suitable initial and boundary conditions independent of $\varepsilon$. $W$ is a Brownian sheet on $\QT$.
We denote by $g(t,x)$ the solution of the homogeneous linear problem
$\mathcal L g=0$ with the prescribed initial values.
{Throughout the paper, unless stated otherwise, $b$ and $\sigma$ are
assumed to be globally Lipschitz on $\R$.  In particular, both $b$ and $\sigma$
satisfy the linear growth condition.}
Let $G_t(x,z)$ denote the Green function of the linear problem associated
with $\mathcal L$ and these boundary conditions: if $v$ solves
$\mathcal L v=\varphi$ with zero initial values, then
$$
{        v(t,x)=\int_0^t\int_{\mathcal O}G_{t-s}(x,z)\varphi(s,z) \ud z \ud s.}
$$
Formally, the skeleton equation is obtained by replacing the small-noise
forcing $\sqrt{\varepsilon}\dot W$ with a deterministic control
$h\in L^2(\mathcal O_T)$.  More precisely, whenever the corresponding
Freidlin--Wentzell LDP holds on $C(\OT)$, its good rate function is
$$
        J(f)
        =
        \inf_{\{h\in L^2(\mathcal O_T):\ \U(h)=f\}}
        \frac12\|h\|_{L^2(\mathcal O_T)}^2,
        \qquad f\in C(\OT),
$$
where $\Upsilon:L^2(\mathcal O_T)\to C(\OT)$ maps $h$ to the solution
$\Upsilon(h)=u$ of the skeleton equation
\begin{equation}\label{eq:continuous-skeleton}
\begin{aligned}
u(t,x)
={}&g(t,x)
 {+\int_0^t\int_{\mathcal O}
G_{t-s}(x,z)[b(u(s,z))+\sigma(u(s,z))h(s,z)] \ud z \ud s} .
\end{aligned}
\end{equation}
For each fixed $x_0\in\mathcal O$, the terminal evaluation map
$f\mapsto f(T,x_0)$ is continuous on $C(\OT)$.  Hence the contraction
principle shows that $\{u_\varepsilon(T,x_0)\}_{\varepsilon>0}$ satisfies an
LDP on $\R$ with rate function
$$
        I(y)
        =
        \inf_{\{h\in L^2(\mathcal O_T):\ \U(h)(T,x_0)=y\}}
        \frac12\|h\|_{L^2(\mathcal O_T)}^2,
        \qquad y\in\R;
$$
see \cite{Dupuis08} and the references therein.
{The following examples illustrate how the homogeneous term $g$
and the Green kernel $G$ are realized for several concrete linear operators.}
{
\begin{example}[Stochastic heat equation] \label{Example1}
Let $\mathcal L=\partial_t-\Delta$, with $u(0)=u_0$ and homogeneous
Dirichlet boundary conditions. Then $g(t)=e^{t\Delta}u_0$. 
Moreover, $G_t$ is the integral kernel of the heat semigroup, in the sense that
$$
\bigl(e^{t\Delta}\varphi\bigr)(x)
=\int_{\mathcal O}G_t(x,z)\varphi(z)\ud z,
\qquad t>0,\quad \varphi\in L^2(\mathcal O).
$$
\end{example}
\begin{example}[Stochastic wave equation]\label{Example2}
Let $\mathcal L=\partial_{tt}-\Delta$, with $u(0)=u_0$,
$\partial_tu(0)=v_0$, and homogeneous Dirichlet boundary conditions. Then
$$
g(t)=\cos\bigl(t\sqrt{-\Delta}\bigr)u_0+
(-\Delta)^{-1/2}\sin\bigl(t\sqrt{-\Delta}\bigr)v_0.
$$
Moreover, $G_t$ is the integral kernel of the sine propagator, in the sense that
$$
\Bigl((-\Delta)^{-1/2}\sin\bigl(t\sqrt{-\Delta}\bigr)\varphi\Bigr)(x)
=\int_{\mathcal O}G_t(x,z)\varphi(z)\ud z,
\qquad t>0,\quad \varphi\in L^2(\mathcal O).
$$
\end{example}
\begin{example}[Fourth-order parabolic equation] \label{Example3}
Let $\mathcal L=\partial_t+\Delta^2$, with $u(0)=u_0$ and boundary
conditions \(u\big|_{\partial\mathcal O}=0,~
\Delta u\big|_{\partial\mathcal O}=0\). Then
$
g(t)=e^{-t\Delta^2}u_0
$. 
Moreover, $G_t$ is the integral kernel of the fourth-order semigroup, namely
$$
\bigl(e^{-t\Delta^2}\varphi\bigr)(x)
=\int_{\mathcal O}G_t(x,z)\varphi(z)\ud z,
\qquad t>0,\quad \varphi\in L^2(\mathcal O).
$$
\end{example}}

We now introduce the time-grid and space-grid projection maps used by the fully discrete
schemes for \eqref{eq:abstract-spde}.  Let $\{m_n\}_{n\ge1}\subset\N^+$ be an arbitrary sequence of
time-step numbers with $m_n\to+\infty$, and write $\tau_n:=T/m_n$.  If
$t_j=j\tau_n$, we denote by $\lfloor s\rfloor_n$ the left endpoint of the
time interval containing $s$, that is,
$
        \lfloor s\rfloor_n:=t_j,~ s\in[t_j,t_{j+1}),
$
{and set $\lfloor T\rfloor_n:=t_{m_n-1}$.}
{ For each $n\ge1$, let $\mathscr T_n$ be the spatial mesh at
level $n$ and let $\kappa_n:\mathcal O\to\overline{\mathcal O}$ be a map that
is constant on every cell of $\mathscr T_n$.  We impose the mesh-consistency
condition
$$
        \delta_n:=\sup_{z\in\mathcal O}|z-\kappa_n(z)|\to0
        \qquad\text{as }n\to\infty.
$$
For the one-dimensional FDM considered later,
$\mathscr T_n$ is the uniform partition of $\mathcal O=(0,1)$ and $\kappa_n$ maps each
cell to its left endpoint, namely
$\kappa_n(z)=\lfloor nz\rfloor/n$.
More generally, in higher dimensions, write $\mathscr T_n=\{K\}$ and choose,
for each cell $K\in\mathscr T_n$, a representative point
$\xi_K\in\overline K$, such as its center, a vertex, or a quadrature point.
Then define
$   \kappa_n(z):=\xi_K,~z\in K.$}
Let $u_\varepsilon^{n,m_n}$ denote the fully discrete numerical
solution of \eqref{eq:abstract-spde} on the grid with spatial mesh
$\mathscr T_n$ and time step $\tau_n$.  

In what follows, we restrict attention
to fully discrete schemes for which, at each fixed grid, the terminal family
$\{u_\varepsilon^{n,m_n}(T,x_0)\}_{\varepsilon>0}$ satisfies an LDP on $\R$
and its good rate function is obtained from the controlled version of the same
scheme.  To state this control representation, we introduce the scheme-induced
skeleton map.  We write $\mathcal X_n$ for the finite-dimensional
space of controls that are piecewise constant on the space-time grid; in
particular, $\mathcal X_n$ is a linear subspace of $L^2(\mathcal O_T)$.  For
$h_n\in\mathcal X_n$, the discrete skeleton map is defined by  $ \Un(h_n):=u_n,$
where $u_n\in C(\OT)$ is the chosen continuous
interpolation in time and space of the discrete skeleton.  {In the formula below, set
$F_1=b$, $F_2=\sigma$, $h_n^1\equiv1$, and $h_n^2=h_n$.}  Then $u_n$ satisfies
\begin{equation}\label{eq:fully-discrete-mild}
\begin{aligned}
u_n(t,x)
={}&g_n(t,x)
{+\sum_{k=1}^2\int_0^t\int_{\mathcal O}
G_n(t,\lfloor s\rfloor_n;x,z)
F_k\!\left(u_n(\lfloor s\rfloor_n,\kappa_n(z))\right)
h_n^k(s,z) \ud z \ud s .}
\end{aligned}
\end{equation}
Here $G_n(t,r;x,z)$ denotes the Green function induced by
the numerical scheme.  With this notation, the scheme-induced rate function of
$u_\varepsilon^{n,m_n}(T,x_0)$ is
$$
        I^{n,m_n}(y)
        =
        \inf_{\{h_n\in\mathcal X_n:\, \Un(h_n)(T,x_0)=y\}}
        \frac12\|h_n\|_{L^2(\mathcal O_T)}^2,
        \qquad y\in\R.
$$
This property is verified for the two fully discrete schemes in
Section~\ref{subsec:concrete-fully-discrete-schemes}: the finite difference method combined with the exponential Euler method (FDM--EEM) for the stochastic wave equation and the FDM--implicit Euler scheme for the fourth-order stochastic parabolic equation.  Their scheme-specific skeleton maps and
fixed-grid LDPs, including the identification of their rate functions with
the above control representation, are given in
Appendix~\ref{app:fixed-grid-LDPs}.
\begin{remark}
For any integrator whose discrete propagator is invariant under time
translations, $G_n(t,r;x,z)$ depends only on the time difference $t-r$.  This
includes standard exponential time integrators applied with a fixed step size
to autonomous problems.
In this case, it may be written as the one-parameter Green function
$\mathsf G_n(t-r;x,z)$. We retain the two-parameter notation for schemes whose continuous-time interpolation is not invariant under arbitrary time translations, including general \(\theta\)-schemes with nonuniform time steps or time-dependent coefficients.
\end{remark}

\subsection{\texorpdfstring{{Abstract variational framework and
minimum action convergence}}{Abstract variational framework and minimum action
convergence}}
\label{sec:abstract-framework}

The preceding subsection provides the continuous and discrete skeleton maps
arising from the scalar SPDE.  Only the control spaces, skeleton maps, and the
terminal observation enter the variational argument.  We now isolate these
objects so that the convergence proof does not depend on the form of a
particular Green function.

Let $\Hc$ be a Hilbert space and, for each $n$, let $\Xc_n$ be a linear
subspace of $\Hc$. Let $\Ec$ be a metric space, and let
$$
        \U:\Hc\to\Ec, \U_n:\Xc_n\to\Ec, \Phi:\Ec\to\R.
$$
We assume throughout the abstract framework that $\Phi$ is
continuous.

Before introducing the abstract action functionals, we record how these
objects are realized in the scalar SPDE framework of
Section~\ref{sec:spde-notation}.

\begin{setting}[Scalar SPDE realization of the abstract framework]
\label{set:scalar-spde-realization}
In the scalar SPDE framework of Section~\ref{sec:spde-notation}, take
$$
\Hc=L^2(\QT),\qquad
\Ec=C(\OT),\qquad
\Phi(f)=f(T,x_0),
$$
and let $\Xc_n$ be the finite-dimensional space of controls that are
piecewise constant on the space-time grid.  The maps $\U$ and $\Un$ are the
skeleton maps defined by \eqref{eq:continuous-skeleton} and
\eqref{eq:fully-discrete-mild}, respectively, and the map
$\Un:\Xc_n\to C(\OT)$ is continuous.
\end{setting}

We next define the abstract action functionals associated with these maps.
%For a fixed $y\in\R$, consider the following minization problems
%$$
%        I(y)
%        :=
%        \inf_{\{h\in\Hc:\,\Phi(\U(h))=y\}}
%        \frac12\|h\|_{\Hc}^2 ,\qquad   \widehat I^n(y)
%        :=
%        \inf_{\{h_n\in\Xc_n:\,\Phi(\U_n(h_n))=y\}}
%        \frac12\|h_n\|_{\Hc}^2 .
%$$
For $y\in\R$, define functionals on $\Ec$ by
\begin{gather}
\J_y(f):=
        \inf_{\{h\in\Hc:\ \U(h)=f,\ \Phi(f)=y\}}
        \frac12\|h\|_{\Hc}^2, \label{Jy}\\
\Jny(f):=
        \inf_{\{h_n\in\Xc_n:\ \U_n(h_n)=f,\ \Phi(f)=y\}}
        \frac12\|h_n\|_{\Hc}^2. \label{Jny}
\end{gather}
Then we consider the following two minimization problems
\begin{align} \label{ratefunction}
	  I(y)=\inf_{f\in\Ec}\J_y(f),
	\qquad
	\widehat I^n(y)=\inf_{f\in\Ec}\Jny(f).
\end{align}
Under Setting~\ref{set:scalar-spde-realization}, these minimization problems
are precisely the continuous and fully discrete one-point rate functions
introduced in Section~\ref{sec:spde-notation}; in particular,
$\widehat I^n(y)=I^{n,m_n}(y)$.

To study the convergence of $\widehat I^n$, we will use the theory of $\Gamma$-convergence.  We
first recall the variational principle used below; see, e.g.,
\cite{DalMaso93}.  Let $(\mathcal M,d)$ be a metric space and let
$F_n,F:\mathcal M\to[0,+\infty]$.  We say that $F_n$ $\Gamma$-converges to
$F$ on $\mathcal M$ if the following two conditions hold:
$$
\text{if }x_n\to x\text{ in }\mathcal M,\quad\text{then}~
        F(x)\le \liminf_{n\to\infty}F_n(x_n)\qquad \text{($\Gamma$-liminf~inequality)},
$$
and, for every $x\in\mathcal M$, there exists a recovery sequence
$x_n\to x$ such that
$$
        F(x)\ge \limsup_{n\to\infty}F_n(x_n)  \qquad \text{($\Gamma$-limsup~inequality)}.
$$
The family $\{F_n\}_{n\ge1}$ is equi-coercive if, for every $a\ge0$, the
set
$
        \cup_{n\ge1}\{x\in\mathcal M:\ F_n(x)\le a\}
$
is relatively compact in $\mathcal M$.  We shall use the standard consequence
of equi-coercive $\Gamma$-convergence; see, e.g.,
\cite[Theorem~7.8]{DalMaso93}.
\begin{proposition}\label{prop:gamma}
If $F_n$ $\Gamma$-converges to $F$ and $\{F_n\}$ is equi-coercive, then
$$
        \inf_{x\in\mathcal M}F(x)
        =
        \lim_{n\to\infty}\inf_{x\in\mathcal M}F_n(x).
$$
\end{proposition}

We now state the assumptions that yield equi-coercivity, the lower-bound
inequality, and an admissible recovery sequence.
%The stochastic wave equation, its concrete Green function, and the precise normalization of $\mathcal X_n$ are specified later in Section~\ref{sec:wave-example}, where the kernel requirements are verified.

\begin{assumption}[Compactness and weak stability of $\Upsilon_n$]
\label{ass:abstract-compact-stability}
The following two conditions hold.
\begin{enumerate}
\item[(i)] For every $R>0$, the set
$\{\U_n(h):\,h\in\Xc_n,\ \|h\|_{\Hc}\le R,\ n\ge1\}$
is relatively compact in $\Ec$.

\item[(ii)] Whenever $n_j\to+\infty$ and
$h_j\in\Xc_{n_j}$ satisfy $h_j\rightharpoonup h$ in $\Hc$, one has
$$
        \U_{n_j}(h_j)\to\U(h)
        \qquad\text{in }\Ec\text{ as }j\to\infty.
$$
Hereafter $\rightharpoonup$ denotes weak convergence in $\Hc$.
\end{enumerate}
\end{assumption}

\begin{assumption}[Control approximation]
\label{ass:abstract-control-approx}
There exist operators $Q_n:\Hc\to\Xc_n$, $n\in\mathbb{N}^+$, such that
$
        Q_nh\to h$ as $n\to\infty$ in $\Hc$
for every $h\in\Hc$.
\end{assumption}

\begin{assumption}[Directional equicontinuity of $\Upsilon_n$]
\label{ass:abstract-directional-consistency}
Let $h_n^0\in\Xc_n$, $h\in\Hc$, and
$h_n^0\to h$ in $\Hc$.  For every direction $q\in\Hc$, with
$q_n:=Q_nq$, and every $a_0>0$, the family of maps
$
a\mapsto \Phi(\U_n(h_n^0+a q_n))
$
is equicontinuous on $[-a_0,a_0]$.
\end{assumption}

\begin{assumption}[Local terminal nondegeneracy]
\label{ass:abstract-terminal-nondegeneracy}
The terminal observation map
$$
\mathcal T:\Hc\to\R,
\qquad
\mathcal T(\widetilde h):=\Phi(\U(\widetilde h)),
$$
is $C^1$ on $\Hc$.  Moreover, for every $h\in\Hc$, there exists
$q_h\in\Hc$ such that the directional derivative
$
D\mathcal T(h)[q_h]\ne0.
$
\end{assumption}

\begin{lemma}[Abstract terminal correction]
\label{lem:abstract-direction-correction}
Assume that Assumption~\ref{ass:abstract-compact-stability}(ii) and
Assumptions~\ref{ass:abstract-control-approx},
\ref{ass:abstract-directional-consistency}, and
\ref{ass:abstract-terminal-nondegeneracy} hold.  Let $h_n^0\in\Xc_n$ satisfy
$h_n^0\to h$ in $\Hc$, and let $c_n\to0$.  {Then there
exist $n_0\in\N^+$ and $r_n\in\Xc_n$, $n\ge n_0$, such that}
$
\Phi(\U_n(h_n^0+r_n))
=
\Phi(\U_n(h_n^0))+c_n$ and
$
\lim_{n\to\infty}\|r_n\|_{\Hc}=0.
$
\end{lemma}

\begin{proof}

Let $q_h$ be the direction given by
Assumption~\ref{ass:abstract-terminal-nondegeneracy}.  Set
$q:=q_h$, $q_n:=Q_nq$, and
$$
\Phi_n(a):=\Phi(\U_n(h_n^0+a q_n)),
\qquad
\Phi_\ast(a):=\mathcal T(h+a q)=\Phi(\U(h+a q)).
$$
For every fixed $a\in\R$, Assumption~\ref{ass:abstract-control-approx}
gives
$
h_n^0+a q_n\to h+a q
$
in $\Hc$, and hence also weakly in $\Hc$.  It follows from
Assumption~\ref{ass:abstract-compact-stability}(ii) and the continuity of
$\Phi$ that $\Phi_n(a)\to\Phi_\ast(a)$.  On every compact interval,
Assumption~\ref{ass:abstract-directional-consistency} gives the
equicontinuity of $\{\Phi_n\}_{n\ge1}$.  Therefore, the pointwise convergence
is uniform on that interval.  Thus
$\Phi_n\rightarrow\Phi_\ast$ locally uniformly on $\R$ as $n\to\infty$.  By
Assumption~\ref{ass:abstract-terminal-nondegeneracy},
$\Phi_\ast$ is $C^1$ near $0$, and
$
\Phi_\ast'(0)=D\mathcal T(h)[q]\ne0.
$
Without loss of generality, replacing $q$ by $-q$ if necessary, assume
that $\Phi_\ast'(0)>0$.  Choose $\delta>0$ sufficiently small so that
$\Phi_\ast$ is strictly increasing on
$[-\delta,\delta]$, which implies 
$
\Phi_\ast(-\delta)<\Phi_\ast(0)<\Phi_\ast(\delta).
$
For this fixed $\delta$, choose $\gamma>0$ such that
\begin{equation}\label{eq:Phiast}
\Phi_\ast(-\delta)+3\gamma<\Phi_\ast(0)<\Phi_\ast(\delta)-3\gamma.
\end{equation}
Since $c_n\rightarrow0$ and
$\Phi_n\rightarrow\Phi_\ast$ uniformly on the fixed interval
$[-\delta,\delta]$, for all sufficiently large $n$,
$$
|c_n|<\gamma,
\qquad
\sup_{|a|\le\delta}|\Phi_n(a)-\Phi_\ast(a)|<\gamma.
$$
Together with \eqref{eq:Phiast}, it follows that
$
\Phi_n(-\delta)<\Phi_n(0)+c_n<\Phi_n(\delta).
$
By continuity of $\Phi_n$, there exists $a_n\in(-\delta,\delta)$
such that
$
\Phi_n(a_n)=\Phi_n(0)+c_n.
$ Taking $r_n:=a_nq_n$, we obtain $
\Phi\bigl(\U_n(h_n^0+r_n)\bigr)=
\Phi\bigl(\U_n(h_n^0)\bigr)+c_n.
$  

By
Assumption~\ref{ass:abstract-control-approx}, $Q_nq\to q$ in $\Hc$, which implies that $\{q_n\}$ is bounded in $\Hc$. 
We claim that $\lim_{n\to \infty}a_n=0$.  Otherwise, along a subsequence,
$a_n\rightarrow a_\ast\in[-\delta,\delta]$ with $a_\ast\ne0$.  Since
$\Phi_\ast$ is continuous, $\Phi_n\rightarrow\Phi_\ast$ uniformly on
$[-\delta,\delta]$, and $c_n\to 0$, we have 
%$$
%\begin{aligned}
%|\Phi_n(a_n)-\Phi_\ast(a_\ast)|
%&\le
%\sup_{|a|\le\delta}|\Phi_n(a)-\Phi_\ast(a)|
%+|\Phi_\ast(a_n)-\Phi_\ast(a_\ast)|\rightarrow0
%\qquad\text{as }n\to\infty.
%\end{aligned}
%$$
%Therefore
$$
\Phi_\ast(a_\ast)
=\lim_{n\to\infty}\Phi_n(a_n)
=\lim_{n\to\infty}\bigl(\Phi_n(0)+c_n\bigr)
=\Phi_\ast(0),
$$
contradicting the strict monotonicity of $\Phi_\ast$.  This proves $\lim_{n\to\infty}a_n=0$, which along with the boundedness of $\{q_n\}$ in $\mathcal{H}$ proves 
$
\lim_{n\to\infty}\|r_n\|_{\Hc}=0.
$ \end{proof}

\begin{theorem}[Abstract variational convergence]
\label{thm:abstract-framework}
Under
Assumptions~\ref{ass:abstract-compact-stability}--\ref{ass:abstract-terminal-nondegeneracy}, for every $y\in\R$,
$$
        \widehat I^n(y)\rightarrow I(y)
        \qquad\text{as }n\to\infty.
$$
\end{theorem}

\begin{proof}
By Proposition \ref{prop:gamma} and the definitions of $I$ and \(\widehat I^n\) (see \eqref{ratefunction}), it suffices to prove that $\{\Jny\}_{n\ge1}$ is equi-coercive and $\Gamma$-convergent.

For $a\ge0$, let
$
\mathcal K_n^a:=\{f\in\Ec:\ \Jny(f)\le a\}.
$
We first prove that
$
\mathcal K^a:=
\cup_{n\ge1}\mathcal K_n^a
$
is relatively compact in $\Ec$.  Let $\{f_j\}_{j\ge1}$ be any
sequence in $\mathcal K^a$.  Then, for each $j$, there
exists an index $n_j$ such that $f_j\in\mathcal K_{n_j}^a$.  It follows from \eqref{Jny} that there exists
$h_j\in\Xc_{n_j}$ such that
$$
\U_{n_j}(h_j)=f_j,
\qquad
\frac12\|h_j\|_{\Hc}^2\le a+1.
$$
Hence $\{h_j\}$ is bounded in $\Hc$, and
Assumption~\ref{ass:abstract-compact-stability}(i) yields a convergent subsequence
of $\{f_j\}$ in $\Ec$, which proves that
$\mathcal K^a$ is relatively compact. Hence, the family $\{\Jny\}_{n\ge1}$ is equi-coercive.

	To prove the $\Gamma$-liminf inequality, let $f_n\rightarrow f$ in
$\Ec$ as $n\to\infty$, and set
$
L:=\liminf_{n\to\infty}\Jny(f_n).
$ It suffices to consider the case $L<+\infty$.
Choose a subsequence $\{n_k\}_{k\ge1}$ such that
$
\J_{y}^{n_k}(f_{n_k})\rightarrow L$
as $k\to\infty$. Then for every $k$, select controls $h_{n_k}\in\Xc_{n_k}$ satisfying
\begin{equation}\label{eq:Unkhnk}
\U_{n_k}(h_{n_k})=f_{n_k},
\qquad \Phi(f_{n_k})=y,\qquad
\frac12\|h_{n_k}\|_{\Hc}^2
\le \J_{y}^{n_k}(f_{n_k})+\frac1k.
\end{equation}
In particular, this implies that $\{h_{n_k}\}_{k\ge1}$ is bounded in
$\Hc$.  Since $\Hc$ is a Hilbert space, after passing to a further
subsequence, still denoted by $\{h_{n_k}\}_{k\ge1}$, there exists
$h\in\Hc$ such that $h_{n_k}\rightharpoonup h$ in $\Hc$ as
$k\to\infty$.  Using $f_n\to f$, \eqref{eq:Unkhnk}, and Assumption~\ref{ass:abstract-compact-stability}(ii), we obtain
$
        f=\lim_{k\to\infty}f_{n_k}
        =
        \lim_{k\to\infty}\U_{n_k}(h_{n_k})
        =
        \U(h).
$
Since $\Phi(f_{n_k})=y$ and $\Phi$ is continuous, one has $\Phi(f)=y$,
and thus
\begin{equation*}
\J_y(f)
\le \frac12\|h\|_{\Hc}^2
\le \liminf_{k\to\infty}\frac12\|h_{n_k}\|_{\Hc}^2
\le L=\liminf_{n\to\infty}\Jny(f_n).
\tag{$\Gamma$-liminf inequality}
\end{equation*}

We now prove the $\Gamma$-limsup inequality. Let $f\in\Ec$. If
$\J_y(f)=+\infty$, the constant sequence $f_n:=f$ already satisfies the
$\Gamma$-limsup requirement. We therefore assume 
$\J_y(f)<+\infty$. For $\eta>0$, it follows from \eqref{Jy} that there exists
$h^\eta\in\Hc$ such that
$$
\U(h^\eta)=f,\qquad \Phi(f)=y,
\qquad
\frac12\|h^\eta\|_{\Hc}^2\le \J_y(f)+\eta.
$$
By setting $h_n^0:=Q_nh^\eta\in\mathcal{X}_n$, 
Assumption~\ref{ass:abstract-control-approx} leads to
$h_n^0\rightarrow h^\eta$ in $\Hc$ as $n\to\infty$.  Therefore, by Assumption~\ref{ass:abstract-compact-stability}(ii), we have
$
\U_n(h_n^0)\rightarrow \U(h^\eta)=f$
in $\Ec$ as $n\to\infty$.
Consequently, by the continuity of $\Phi$ again,
$
c_n:=\Phi(f)-\Phi(\U_n(h_n^0))=y-\Phi(\U_n(h_n^0))\rightarrow0$
as $n\to\infty$.
Applying Lemma~\ref{lem:abstract-direction-correction} gives {an
index $n_0$ and a sequence $r_n\in \Xc_n$, $n\ge n_0$, converging
to $0$ in $\Hc$} such that
$
\Phi\bigl(\U_n(h_n^0+r_n)\bigr)=y$ for every $n\ge n_0$.
{For the finitely many indices $n<n_0$, set
$\bar h_n:=h_n^0$ and $f_n^\eta:=\U_n(h_n^0)$. For $n\ge n_0$, set}
$\bar h_n:=h_n^0+r_n$ and $f_n^\eta:=\U_n(\bar h_n)$ so that
$
\bar h_n\rightarrow h^\eta$ in $\Hc$ as $n\to\infty$. 
Consequently, 
\begin{equation*}
\limsup_{n\to\infty}\Jny(f_n^\eta)
\le \lim_{n\to\infty}\frac12\|\bar h_n\|_{\Hc}^2
=\frac12\|h^\eta\|_{\Hc}^2
\le \J_y(f)+\eta.
\tag{$\Gamma$-limsup inequality}
\end{equation*}
Moreover, by Assumption~\ref{ass:abstract-compact-stability}(ii),
$f_n^\eta\rightarrow \Upsilon(h^\eta)=f$ in $\Ec$ as $n\to\infty$.
{For each $j\ge1$, apply the preceding construction with
$\eta=j^{-1}$ and denote the resulting recovery sequence $\{f_n^\eta\}$ by
$\{f_n^j\}_{n\ge1}$.  Choose integers $1\le N_1<N_2<\cdots$ such that,
for every $j\ge1$ and $n\ge N_j$,
$$
d_{\Ec}(f_n^j,f)\le j^{-1},
\qquad
\J_y^n(f_n^j)\le \J_y(f)+2j^{-1}.
$$
For $n\ge N_1$, set
$j(n):=\max\{j\ge1:N_j\le n\}$.  Then $j(n)\to\infty$ and
$f_n^{j(n)}\to f$ in $\Ec$, while
$$
\limsup_{n\to\infty}\J_y^n\bigl(f_n^{j(n)}\bigr)\le \J_y(f).
$$
This proves the $\Gamma$-limsup inequality with the recovery sequence $\{f_n^{j(n)}\}$.}
%Hence $\Jny$ $\Gamma$-converges to $\J_y$.  
\end{proof}

\begin{theoremcorollary}[Subsequential convergence of minimizing controls and paths]
\label{cor:subsequential-minimizers}
Assume the hypotheses of Theorem~\ref{thm:abstract-framework}, and fix
$y\in\R$ such that $I(y)<+\infty$.  Let $\epsilon_n\to0$ with
$\epsilon_n\ge0$, and let $h_n\in\Xc_n$ be asymptotically minimizing
controls satisfying
\[
        \Phi(\U_n(h_n))=y,
        \qquad
        \frac12\|h_n\|_{\Hc}^2
        \le \widehat I^n(y)+\epsilon_n.
\]
Define the associated controlled paths by $u_n:=\U_n(h_n)$.
Then there exist a subsequence $\{n_k\}_{k\ge1}$ and a control $h^*\in\Hc$ such that, with $u^*:=\U(h^*)$,
\[
        h_{n_k}\to h^* \quad\text{in }\Hc,
        \qquad
        u_{n_k}\to u^* \quad\text{in }\Ec,
\]
and
\[
        \Phi(u^*)=y,
        \qquad
        \frac12\|h^*\|_{\Hc}^2=I(y).
\]
Consequently, $h^*$ is a minimizing control and 
 $u^*$ is the associated minimum action path for the continuous constrained
problem $\inf_{f\in\Ec}\J_y(f)$.  If the minimizing control for the continuous problem is unique, then the
full sequences $h_n\to h^*$ in $\Hc$ and $u_n\to u^*$ in $\Ec$ without
extracting a subsequence.
\end{theoremcorollary}

\begin{proof}
By the definition of $\J_y^n$,
\[
0\le \J_y^n(u_n)-\widehat I^n(y)
\le \frac12\|h_n\|_{\Hc}^2-\widehat I^n(y)
\le \epsilon_n\rightarrow0.
\]
Thus $\{u_n\}$ is an asymptotically minimizing sequence of controlled paths.
Theorem~\ref{thm:abstract-framework} and the asymptotic minimality imply that
$\{h_n\}_{n\ge1}$ is bounded in $\Hc$.  Hence, after passing to a
subsequence,
$h_{n_k}\rightharpoonup h^*$ in $\Hc$.  By
Assumption~\ref{ass:abstract-compact-stability}(ii),
        $u_{n_k}\to u^*~\text{in }\Ec$. 
The continuity of $\Phi$ gives $\Phi(u^*)=y$.  Therefore,
\[
\begin{aligned}
 I(y)
 \le \frac12\|h^*\|_{\Hc}^2
 \le \liminf_{k\to\infty}\frac12\|h_{n_k}\|_{\Hc}^2
 &\le \limsup_{k\to\infty}\frac12\|h_{n_k}\|_{\Hc}^2\\
 \le \lim_{k\to\infty}
       \bigl(\widehat I^{n_k}(y)+\epsilon_{n_k}\bigr)
 = I(y).
\end{aligned}
\]
Consequently, $h^*$ is a minimizing control for the continuous constrained problem and
$\|h_{n_k}\|_{\Hc}\to\|h^*\|_{\Hc}$.  Weak convergence together with
convergence of the norms in the Hilbert space $\Hc$ yields
$h_{n_k}\to h^*$ strongly in $\Hc$.  Since $u^*=\U(h^*)$ and $h^*$ attains
$I(y)$,
\[
        I(y)\le \J_y(u^*)
        \le \frac12\|h^*\|_{\Hc}^2=I(y).
\]
Thus $u^*$ is a minimum action path for the continuous constrained problem.
If the minimizing control for the continuous problem is unique, the same argument
applies to every subsequence and proves convergence of both full sequences.
\end{proof}

%Theorem~\ref{thm:abstract-framework} separates the lower-bound mechanism from
%the terminal correction used in the recovery sequence.  Compactness and weak
%stability control the former; directional consistency and terminal
%nondegeneracy control the latter.  We next derive both sets of properties from
%estimates on the continuous and discrete Green functions.

\section{A Green-function criterion for diagonal convergence}\label{S:KB}
In this section, we introduce five kernel conditions that provide a verifiable criterion for applying
Theorem~\ref{thm:abstract-framework}.  
%Conditions \Kone--\Kfour{} provide consistency, uniform
%control estimates, and compactness of skeleton trajectories.  Condition
%\Kfive{} enters only in the argument for terminal nondegeneracy.  

\subsection{Green functions}
{Unless otherwise specified, all suprema over spatial
variables, time variables and the discretization index are taken over
$x,z\in\overline{\mathcal O}$, $t,s\in[0,T]$ and $n\in\mathbb N^+$,
respectively. Throughout the paper, $C$ denotes a generic positive constant
{that may vary from one occurrence to the next and may depend on
$g$, $T$, $b$, and $\sigma$,} and we do not normally display such
dependencies. When necessary, we write $C_R$ or $C(R,\ldots)$ to emphasize the dependence of $C$ upon the parameter $R$.}
{For later use in the compactness estimate, we introduce two discrete
kernel increments and their limiting counterparts.  If $t\ge t'$, set
{\begin{align*}
\mathcal K_{n,1}(t,t',x,x')
&:=
{\int_0^{t'}\int_{\mathcal O}
\bigl|G_n(t,\lfloor s\rfloor_n;x,z)
-G_n(t',\lfloor s\rfloor_n;x',z)\bigr|^2 \ud z\ud s,}\\
\mathcal K_{1}(t,t',x,x')
&:=
{\int_0^{t'}\int_{\mathcal O}
\bigl|G_{t-s}(x,z)-G_{t'-s}(x',z)\bigr|^2 \ud z\ud s,
}\\
\mathcal K_{n,2}(t,t',x)
&:=
{\int_{t'}^{t}\int_{\mathcal O}
\bigl|G_n(t,\lfloor s\rfloor_n;x,z)\bigr|^2 \ud z\ud s,}\\
\mathcal K_{2}(t,t',x)
&:={
\int_{t'}^{t}\int_{\mathcal O}
\bigl|G_{t-s}(x,z)\bigr|^2 \ud z\ud s .}
\end{align*}}
If $t<t'$, these quantities are defined by exchanging $(t,x)$ and
$(t',x')$ on the right-hand side above.}

\begin{description}
\item[\phantomsection\label{cond:K1}(K1)] \emph{[Initial-term consistency]}
$
        g_n\to g$ as $n\to \infty$ in $C(\OT)$. 

\item[\phantomsection\label{cond:K2}(K2)] \emph{[Green-function consistency]}
\begin{equation}\label{eq:green-consistency}
\varepsilon_n:=\sup_{(t,x)\in\OT}
{\int_0^t\int_{\mathcal O}
\left|
G_n(t,\lfloor s\rfloor_n;x,z)-G_{t-s}(x,z)
\right|^2 \ud z \ud s}
\rightarrow0
~\text{as }n\to\infty .
\end{equation}

\item[\phantomsection\label{cond:K3}(K3)] \emph{[Integrable kernel bound]}
There exists a nonincreasing function
$\varrho\in L^1(0,T;[0,+\infty))$, independent of $n$, such that
$$
        \sup_{n\in\N^+}
{        \sup_{\substack{0\le r<t\le T, \, x\in\overline{\mathcal O}}}
        \int_{\mathcal O}|G_n(t,r;x,z)|^2 \ud z\le \varrho(t-r).}
$$

\item[\phantomsection\label{cond:K4}(K4)] \emph{[Continuity modulus]}
{
$  \lim_{\delta\downarrow0}
        \sup_{\substack{|t-t'|+|x-x'|\le\delta}}
        \mathcal K_{1}(t,t',x,x')=0.$
}

\item[\phantomsection\label{cond:K5}(K5)] \emph{[Nondegeneracy]}
{For the fixed observation point $x_0\in\mathcal O$,}
$$
{\int_0^T\int_{\mathcal O}
\left|G_{T-s}(x_0,z)\right|^2
\ud z\ud s>0.}
$$
\end{description}

\begin{lemma}
\label{lem:transfer-limiting-modulus}
{
{Under the conditions \Ktwo--\Kfour, we have}
$$
\lim_{\delta\downarrow0}\limsup_{n\to\infty}
\sup_{\substack{|t-t'|+|x-x'|\le\delta}}
\mathcal K_{n,1}(t,t',x,x')=0.
$$
\begin{equation}\label{eq:discrete-time-strip-estimate}
\lim_{\delta\downarrow 0}\sup_{n}
{\sup_{\substack{0\le t-t'\le\delta,\, x\in\overline{\mathcal O}}}
\mathcal K_{n,2}(t,t',x)=0,\quad
\lim_{\delta\downarrow 0}
\sup_{\substack{0\le t-t'\le\delta,\, x\in\overline{\mathcal O}}}
        \mathcal K_2(t,t',x)=0.}
\end{equation}
}
\end{lemma}

\begin{proof}
{
For every $t,t'\in[0,T]$ and $x,x'\in\overline{\mathcal O}$, insert
$G_{t-s}(x,z)$ and $G_{t'-s}(x',z)$ into the definition of
$\mathcal K_{n,1}$ and use $|a+b+c|^2\le3(|a|^2+|b|^2+|c|^2)$.  This gives
$$
\mathcal K_{n,1}(t,t',x,x')
\le 6\varepsilon_n+3\mathcal K_1(t,t',x,x').
$$
Taking the supremum over $|t-t'|+|x-x'|\le\delta$, then taking the upper
limit as $n\to\infty$, and finally letting $\delta\downarrow0$ proves
the assertion.
}

{
It remains to prove \eqref{eq:discrete-time-strip-estimate}. Fix
$0\le t-t'\le\delta$ and $x\in\overline{\mathcal O}$. By \Kthree{},
$\lfloor s\rfloor_n\le s$, and the monotonicity of $\varrho$,
\begin{align*}
\mathcal K_{n,2}(t,t',x)
&\le\int_{t'}^t\varrho\bigl(t-\lfloor s\rfloor_n\bigr)\ud s
\le\int_{t'}^t\varrho(t-s)\ud s
\le\int_0^\delta\varrho(r)\ud r.
\end{align*}
Moreover, by $|a|^2\le2|a-b|^2+2|b|^2$ and \Ktwo{},
\begin{align*}
\mathcal K_2(t,t',x)
&\le 2\int_{t'}^t\int_{\mathcal O}
\left|G_{t-s}(x,z)-G_n(t,\lfloor s\rfloor_n;x,z)\right|^2\ud z\ud s
    +2\mathcal K_{n,2}(t,t',x)\\
&\le2\varepsilon_n+2\int_0^\delta\varrho(r)\ud r.
\end{align*}
Consequently, we arrive at 
\begin{align}\label{eq:limiting-time-strip-rho}
{\sup_{n\ge1}\sup_{\substack{0\le t-t'\le\delta,\,x\in\overline{\mathcal O}}}
\mathcal K_{n,2}(t,t',x)}
&\le\int_0^\delta\varrho(r)\ud r,\\
\sup_{\substack{0\le t-t'\le\delta,\,x\in\overline{\mathcal O}}}
\mathcal K_2(t,t',x)
&\le2\int_0^\delta\varrho(r)\ud r.
\label{eq:limiting-time-strip-rho-continuous}
\end{align}
}
The integrability of $\varrho$ now gives \eqref{eq:discrete-time-strip-estimate}.
\end{proof}
For the limiting
kernel $G$, taking $t'=0$ and  $\delta=T$ 
in {\eqref{eq:limiting-time-strip-rho-continuous}}, we have
\begin{equation}\label{eq:limiting-kernel-global-l2-bound}
{
{        \sup_{t,x}\int_0^t\int_{\mathcal O}
        |G_{t-s}(x,z)|^2 \ud z \ud s
        \le 2\int_0^T\varrho(r)\ud r<+\infty .}}
\end{equation}
{We also need the pointwise-in-time counterpart of
{\Kthree{} for the Green function $G$.  Fix $(t,x)\in(0,T]\times\overline{\mathcal O}$
and define
$
 d_n(s):=\int_{\mathcal O}
 \left|G_n(t,\lfloor s\rfloor_n;x,z)-G_{t-s}(x,z)\right|^2\ud z,
~ s\in(0,t).
$
By \Ktwo{}, $d_n\to0$ in $L^1(0,t)$.  Hence there is a subsequence
$\{{n_j}\}_{j\ge1}$ for which $d_{n_j}(s)\to0$ for almost every
$s\in(0,t)$.  Equivalently, for almost every such $s$,
$
 G_{n_j}(t,\lfloor s\rfloor_{n_j};x,\cdot)
 \to G_{t-s}(x,\cdot)
~\text{in }L^2(\mathcal O).
$
In particular, we have
{
\begin{equation}\label{eq:limiting-kernel-pointwise-rho}
 \begin{aligned}
 \int_{\mathcal O}|G_{t-s}(x,z)|^2\ud z
 &=\lim_{j\to\infty}\int_{\mathcal O}
 |G_{n_j}(t,\lfloor s\rfloor_{n_j};x,z)|^2\ud z\\
 &\le \limsup_{j\to\infty}\varrho(t-\lfloor s\rfloor_{n_j})
 \le \varrho(t-s)
 \quad\text{for a.e. }s\in(0,t),
 \end{aligned}
\end{equation}}
where we have also used
$\lfloor s\rfloor_{n_j}\le s$, the monotonicity of $\varrho$, and \Kthree{}}.}

{
\begin{lemma}%[The continuous Green operator]
\label{lem:green-convolution}
Let the conditions \Ktwo--{\Kfour} hold. Define
\begin{align}\label{eq:G}
(\mathcal G\phi)(t,x):=
\int_0^t\int_{\mathcal O}G_{t-s}(x,z)\phi(s,z)\ud z\ud s,
\qquad \phi\in L^2(\mathcal O_T).
\end{align}
Then $\mathcal G$ is a compact linear operator from $L^2(\mathcal O_T)$ into
$C(\OT)$.
\end{lemma}

\begin{proof}
The linearity of $\mathcal G$ is immediate from its definition. For
$\phi\in L^2(\mathcal O_T)$, the Cauchy--Schwarz inequality and
\eqref{eq:limiting-kernel-global-l2-bound} give
\begin{equation}\label{eq:mathphi}
\sup_{(t,x)}|(\mathcal G\phi)(t,x)|^2
\le
2\|\phi\|_{L^2(\mathcal O_T)}^2\int_0^T\varrho(r)\ud r.
\end{equation}
Moreover, if $T\ge t\ge t'\ge0$, then
$$
\begin{aligned}
|(\mathcal G\phi)(t,x)-(\mathcal G\phi)(t',x')|
\le
\sqrt{2}\|\phi\|_{L^2(\mathcal O_T)}
\left(\mathcal K_1(t,t',x,x')+\mathcal K_2(t,t',x)\right)^{1/2}.
\end{aligned}
$$
The preceding increment estimate, together with \Kfour{} and
\eqref{eq:discrete-time-strip-estimate}, shows that
$\mathcal G\phi\in C(\OT)$ for every $\phi\in L^2(\mathcal O_T)$ and that the
image under $\mathcal G$ of the unit ball of $L^2(\mathcal O_T)$ is
equicontinuous. Equation~\eqref{eq:mathphi} then shows that this image
is uniformly bounded in $C(\OT)$. Since $\OT$ is compact, the
Arzel\`a--Ascoli theorem implies that the image is relatively compact
in $C(\OT)$. Therefore, $\mathcal G:L^2(\mathcal O_T)\to C(\OT)$ is compact.
\end{proof}
}

\subsection{Compactness and weak stability of the skeleton maps}
{
\begin{proposition}[Volterra--Gronwall inequality]
\label{prop:volterra-gronwall}
Let $A,B\ge0$, let $\varrho\in L^1(0,T)$ be nonnegative, and let
$f:[0,T]\to[0,\infty)$ be bounded and measurable.  If
$$
f(t)\le A+B\int_0^t\varrho(t-s)f(s)\ud s,
\qquad 0\le t\le T,
$$
then
$
\sup_{0\le t\le T}f(t)\le C(B,T,\varrho)A.
$
\end{proposition}}

The finite constant $C(B,T,\varrho)$ may be expressed in terms of
the Volterra resolvent generated by $B\varrho$; see
\cite[Lemma~2.2]{Zhang2010Volterra} for the proof of Proposition \ref{prop:volterra-gronwall}.
We next present a uniform a priori bound for the discrete skeleton maps $\Upsilon_n$, $n\ge1$.
{
\begin{lemma}
\label{lem:discrete-apriori-bound}
Assume that $b$ and $\sigma$ are Lipschitz continuous.  If conditions \Kone{} and \Kthree{} hold, then for every $R>0$, there exists $C_R>0$ such that
$$
\sup_{n\ge1}\sup_{\{h_n\in\mathcal X_n:
\|h_n\|_{L^2(\mathcal O_T)}\le R\}}
\|\Upsilon_n(h_n)\|_{C(\OT)}\le C_R.
$$
\end{lemma}

\begin{proof}
Fix $R>0$, let $h_n\in\mathcal X_n$ satisfy
$\|h_n\|_{L^2(\mathcal O_T)}\le R$. Set $u_n:=\Upsilon_n(h_n)$ and 
$$
M_n(t):=\sup_{0\le r\le t,\ x\in\overline{\mathcal O}}|u_n(r,x)|^2.
$$
It follows from \Kone{}, the linear growth of $b$ and $\sigma$, the
Cauchy--Schwarz inequality, \Kthree{}, 
$\lfloor s\rfloor_n\le s$, and \eqref{eq:fully-discrete-mild} that for
$0\le r\le t$ and $x\in\overline{\mathcal O}$,
$$
|u_n(r,x)|^2
\le C_R+C_R\int_0^r\varrho(r-s)M_n(s)\ud s.
$$
  Since $M_n$ is non-decreasing and
$\varrho\ge0$,
$$
\begin{aligned}
\int_0^r\varrho(r-s)M_n(s)\ud s
&=\int_0^r\varrho(v)M_n(r-v)\ud v\\
&\le\int_0^t\varrho(v)M_n(t-v)\ud v
=\int_0^t\varrho(t-s)M_n(s)\ud s.
\end{aligned}
$$
Taking the supremum over $r$ and $x$ therefore yields
$$
M_n(t)\le C_R+C_R\int_0^t\varrho(t-s)M_n(s)\ud s,\quad 0< t\le T.
$$
Proposition~\ref{prop:volterra-gronwall} yields
$\sup_{n\ge1}\|u_n\|_{C(\OT)}\le C_R$, as claimed.
\end{proof}
In the following lemma, we prove the local Lipschitz continuity of the discrete skeleton maps $\Upsilon_n$, $n\ge1$.
\begin{lemma}
\label{lem:lipschitz-strong-consistency}
Assume that $b$ and $\sigma$ are Lipschitz continuous.  Then under
 conditions \Kone{} and \Kthree{}, for every $R>0$, there exists $C_R>0$ such that for any $n\ge1$ and any $h_n,\tilde h_n\in\mathcal X_n$ satisfying
$\|h_n\|_{L^2(\mathcal O_T)}+\|\tilde h_n\|_{L^2(\mathcal O_T)}\le R$,
$$
\|\Upsilon_n(h_n)-\Upsilon_n(\tilde h_n)\|_{C(\OT)}
\le C_R\|h_n-\tilde h_n\|_{L^2(\mathcal O_T)}.
$$
%In particular, every fixed-grid map $\Upsilon_n:\mathcal X_n\to C(\OT)$
%is continuous.
\end{lemma}
The conclusion of Lemma~\ref{lem:lipschitz-strong-consistency} follows from Lemma~\ref{lem:discrete-apriori-bound}, \Kthree{}, and Proposition~\ref{prop:volterra-gronwall}. We omit the routine details.

%\begin{proof}
%Let $u_n:=\Upsilon_n(h_n)$ and
%$\tilde u_n:=\Upsilon_n(\tilde h_n)$, and define
%$$
%D_n(t):=\sup_{0\le r\le t,\ x\in\overline{\mathcal O}}
%|u_n(r,x)-\tilde u_n(r,x)|^2.
%$$
%By Lemma~\ref{lem:discrete-apriori-bound}, $\|u_n\|_{C(\OT)}+\|\tilde u_n\|_{C(\OT)}\le C_R$.  Subtracting
%their discrete mild equations, using the Lipschitz continuity of $b$ and
%$\sigma$, and applying the Cauchy--Schwarz inequality together with
%\Kthree{}, gives
%$$
%D_n(t)
%\le
%C_R\|h_n-\tilde h_n\|_{L^2(\mathcal O_T)}^2
%+C_R\int_0^t\varrho(t-s)D_n(s)\ud s.
%$$
%Proposition~\ref{prop:volterra-gronwall} proves the asserted estimate.
%\end{proof}
}

{
\begin{corollary}
\label{lem:continuous-skeleton-local-bounds}
Assume that $b$ and $\sigma$ are Lipschitz continuous and that conditions
\Ktwo{} and \Kthree{} hold.  Then, for every $R>0$, there exists $C_R>0$
such that:
\begin{enumerate}
\item[(i)]
$\sup\limits_{\{h\in L^2(\mathcal O_T): \|h\|_{L^2(\mathcal O_T)}\le R\}}
\|\U(h)\|_{C(\OT)}
\le C_R$.
\item[(ii)] For any $h_1,h_2\in L^2(\mathcal O_T)$ satisfying
$\|h_1\|_{L^2(\mathcal O_T)}+\|h_2\|_{L^2(\mathcal O_T)}\le R$,
\[
\|\U(h_1)-\U(h_2)\|_{C(\OT)}
\le C_R\|h_1-h_2\|_{L^2(\mathcal O_T)}.
\]
\end{enumerate}
\end{corollary}

The conclusion of Corollary~\ref{lem:continuous-skeleton-local-bounds} follows from \eqref{eq:limiting-kernel-global-l2-bound}, \eqref{eq:limiting-kernel-pointwise-rho}, and Proposition~\ref{prop:volterra-gronwall}; we omit the routine details.

%\begin{proof}
%For (i), apply the Cauchy--Schwarz inequality and the linear growth of
%$b$ and $\sigma$ to \eqref{eq:continuous-skeleton}, using the pointwise
%Green-function bound \eqref{eq:limiting-kernel-pointwise-rho}; the same
%a priori argument as in Lemma~\ref{lem:discrete-apriori-bound}, followed
%by Proposition~\ref{prop:volterra-gronwall}, gives the asserted uniform
%bound.  Having established (i), subtract the continuous skeleton
%equations for $\U(h_1)$ and $\U(h_2)$.  Part (i), the Cauchy--Schwarz
%inequality, \eqref{eq:limiting-kernel-pointwise-rho}, and the global bound
%\eqref{eq:limiting-kernel-global-l2-bound}, followed again by
%Proposition~\ref{prop:volterra-gronwall}, give (ii), exactly as in
%Lemma~\ref{lem:lipschitz-strong-consistency}.
%\end{proof}
}

\begin{lemma}
\label{lem:compactness-kernel}
Assume that $b$ and $\sigma$ are Lipschitz continuous. 
{Then under Setting~\ref{set:scalar-spde-realization}, the conditions \Konefour{} imply}
Assumption~\ref{ass:abstract-compact-stability}(i).
\end{lemma}

\begin{proof}
{For $R>0$, set
$
 E_R:=\{\Un(h):\ h\in\mathcal X_n,\ \|h\|_{L^2(\mathcal O_T)}\le R,\ n\ge1\}.
$
Let
$$
        u_j=\Upsilon_{n_j}(h_j),\qquad
        h_j\in\mathcal X_{n_j},\qquad
        \|h_j\|_{L^2(\mathcal O_T)}\le R,\quad j\in\mathbb{N}^+
$$
be an arbitrary sequence from the set $E_R$.  It suffices to
extract a subsequence of $\{u_j\}$ converging in $C(\OT)$.
{We distinguish two cases.  If $\{n_j\}_{j\ge1}$ is bounded,
then, after passing to a subsequence, we may assume that
$n_j\equiv n_*$.}  Since the grid-control space
$\mathcal X_{n_*}$ is finite dimensional, the ball
$
        \{h\in\mathcal X_{n_*}:\|h\|_{L^2(\mathcal O_T)}\le R\}
$
is compact.  By the fixed-grid continuity supplied by
{Lemma~\ref{lem:lipschitz-strong-consistency}}, the image of this
ball under $\Upsilon_{n_*}$ is compact.  Hence $\{u_j\}$ has a convergent
subsequence.

{If $\{n_j\}_{j\ge1}$ is unbounded, it has a subsequence
tending to $+\infty$.  Passing to and relabeling this subsequence, we may
assume that $n_j\to+\infty$; for notational simplicity, write}}
$$
        u_n=\Un(h_n),\qquad
        h_n\in\mathcal X_n,\qquad
        \|h_n\|_{L^2(\mathcal O_T)}\le R,\quad n\in \mathbb{N}^+ .
$$
{By Lemma~\ref{lem:discrete-apriori-bound}, the sequence
$\{u_n\}$ is uniformly bounded in $C(\OT)$.  Since $F_1=b$ and $F_2=\sigma$ have
linear growth, $F_1(u_n)$ and $F_2(u_n)$ are uniformly bounded in
$C(\OT)$.}  Set
$
        w_n(t,x):=u_n(t,x)-g_n(t,x).
$
We next prove that $\{w_n\}$ is equicontinuous.  Let
$|t-t'|+|x-x'|\le\delta$, and write
$
        w_n(t,x)-w_n(t',x')
        =
        {\sum_{k=1}^2\Delta_{n,k}(t,t',x,x')},
$
where
\begin{align*}
&\Delta_{n,k}(t,t',x,x')
:={}
\int_0^t\int_{\mathcal O}
G_n(t,\lfloor s\rfloor_n;x,z)
F_k\!\left(u_n(\lfloor s\rfloor_n,\kappa_n(z))\right)h_n^k(s,z)\ud z\ud s\\
&-
\int_0^{t'}\int_{\mathcal O}
G_n(t',\lfloor s\rfloor_n;x',z)
F_k\!\left(u_n(\lfloor s\rfloor_n,\kappa_n(z))\right)h_n^k(s,z)\ud z\ud s,
\quad k=1,2.
\end{align*}
Since $F_k(u_n)$ is uniformly bounded and
$\|h_n^k\|_{L^2(\mathcal O_T)}\le C(R,T)$, we obtain from the Cauchy--Schwarz inequality that for $k=1,2$, 
$$
\begin{aligned}
|\Delta_{n,k}(t,t',x,x')|
&\le C_R\left(\mathcal K_{n,1}(t,t',x,x')
        +\mathcal K_{n,2}(t,t',x)\right)^{1/2}\\
\end{aligned}
$$
%Therefore,
%$$
%        |w_n(t,x)-w_n(t',x')|
%        \le C_R
%        \left(\mathcal K_{n,1}(t,t',x,x')
%        +\mathcal K_{n,2}(t,t',x)\right)^{1/2}.
%$$
{For every $\epsilon>0$, Lemma~\ref{lem:transfer-limiting-modulus}
first provides $\delta>0$ and then $N_\epsilon\in\mathbb N^+$ such that the
preceding increment is smaller than $\epsilon$ whenever $n\ge N_\epsilon$
and $|t-t'|+|x-x'|\le\delta$.  Each of the finitely many functions
$w_1,\ldots,w_{N_\epsilon-1}$ is uniformly continuous on the compact set
$\OT$; after decreasing $\delta$ if necessary, the same increment bound
holds for these indices as well.  Thus $\{w_n\}$ is equicontinuous.}
Moreover,
$\{w_n\}$ is uniformly bounded because $\{u_n\}$ is uniformly bounded by
{Lemma~\ref{lem:discrete-apriori-bound}}, while
$\{g_n\}$ is uniformly bounded by \Kone.  Hence the Arzel\`a--Ascoli theorem yields a
subsequence of $\{w_n\}$ converging in $C(\OT)$.  Utilizing \Kone{}  again, 
$u_n=g_n+w_n$ has a convergent subsequence in $C(\OT)$.  This completes the proof.
\end{proof}

\begin{lemma}\label{lem:weak-to-strong}
Assume that $b$ and $\sigma$ are Lipschitz continuous. {Then under Setting~\ref{set:scalar-spde-realization},  the conditions \Konefour{} and
Assumption~\ref{ass:abstract-compact-stability}(i) imply Assumption~\ref{ass:abstract-compact-stability}(ii).}
\end{lemma}

\begin{proof}
Let $n_j\to+\infty$ and $h_j\in\mathcal X_{n_j}$ satisfy
$h_j\rightharpoonup h$ in $L^2(\mathcal O_T)$.  Since \Konefour{} are inherited by
subsequences, we relabel
\[
(\mathcal X_{n_j},\Upsilon_{n_j},G_{n_j},\tau_{n_j},\kappa_{n_j})
\quad\text{by}\quad
(\mathcal X_n,\Upsilon_n,G_n,\tau_n,\kappa_n),
\]
and set
$
        u_n:=\Un(h_n).
$
By \eqref{eq:fully-discrete-mild}, $u_n$ satisfies
\begin{equation}\label{eq:un}
{
\begin{aligned}
u_n(t,x)
={}&g_n(t,x)+\sum_{k=1}^2\int_0^t\int_{\mathcal O}
G_n(t,\lfloor s\rfloor_n;x,z)\\
&\quad\times
F_k\!\left(u_n(\lfloor s\rfloor_n,\kappa_n(z))\right)h_n^k(s,z)
\ud z \ud s,
\end{aligned}}
\end{equation}
where $h_n^1\equiv 1$, $h_n^2=h_n$, $F_1=b$, and $F_2=\sigma$.
Since weakly convergent sequences are bounded, $\{h_n\}$ is bounded in
$L^2(\mathcal O_T)$. By Assumption~\ref{ass:abstract-compact-stability}(i), every
subsequence of $\{u_n\}$ has a further subsequence, still denoted by
$u_n$, such that
$
        u_n\to v$ in $C(\OT;\R)$.
{We claim that $v$ satisfies the continuous skeleton equation with control
$h$, namely, that}
\begin{equation}\label{eq:v}
{
\begin{aligned}
v(t,x)
={}&g(t,x)+\sum_{k=1}^2\int_0^t\int_{\mathcal O}
G_{t-s}(x,z)F_k(v(s,z))h^k(s,z)\ud z\ud s,
\end{aligned}}
\end{equation}
{where $h^1\equiv1$ and
$h^2=h$.}
{Let $\mathscr R_h(v)$ denote the right-hand side of
\eqref{eq:v}.  Without assuming \eqref{eq:v}, subtracting
$\mathscr R_h(v)$ from the discrete equation \eqref{eq:un} gives the
identity}
\begin{gather}\label{eq:nu-v}
{u_n(t,x)-\mathscr R_h(v)(t,x)}=g_n(t,x)-g(t,x)+\sum_{k=1}^2\bigl(A_{n,k}(t,x)+B_{n,k}(t,x)\bigr)+E_n(t,x),
\end{gather}
where
\begin{gather*}
A_{n,k}(t,x)
:=\int_0^t\int_{\mathcal O}
\bigl[
G_n(t,\lfloor s\rfloor_n;x,z)-G_{t-s}(x,z)
\bigr]
F_k\!\left(u_n(\lfloor s\rfloor_n,\kappa_n(z))\right)h_n^k(s,z)
\ud z \ud s,\\
B_{n,k}(t,x)
:=\int_0^t\int_{\mathcal O}G_{t-s}(x,z)
\bigl[F_k\!\left(u_n(\lfloor s\rfloor_n,\kappa_n(z))\right)
-F_k(v(s,z))\bigr]h_n^k(s,z)\ud z \ud s,\\
E_n(t,x)
:=
\int_0^t\int_{\mathcal O}
G_{t-s}(x,z)F_2(v(s,z))(h_n-h)(s,z) \ud z \ud s .
\end{gather*}
For $k=1,2$, the sequence $\{h_n^k\}_{n\ge1}$ is bounded in
$L^2(\mathcal O_T)$. Since $u_n\to v$ in $C(\OT)$,
$\sup_{s\in[0,T]}|s-\lfloor s\rfloor_n|\le\tau_n$, and
$\sup_{z\in\mathcal O}|z-\kappa_n(z)|\le\delta_n$, the uniform continuity of
$v$ yields
\begin{equation}\label{eq:etatn}
   \sup_{\substack{t,x}}
        \bigl|u_n(\lfloor t\rfloor_n,\kappa_n(x))-v(t,x)\bigr|
        \to0 \qquad\text{as } n\to\infty .
\end{equation}
By the Cauchy--Schwarz inequality, the linear growth of $F_k$, the boundedness of $\{u_n\}$,
and \Ktwo{},
$$
        \sup_{t,x}|A_{n,k}(t,x)|\to0
        \qquad\text{as } n\to\infty,
        \quad k=1,2.
$$
Similarly, the Lipschitz continuity of $F_k$,
\eqref{eq:limiting-kernel-global-l2-bound}, and \eqref{eq:etatn} yield
$$
        \sup_{t,x}|B_{n,k}(t,x)|\to0
        \qquad\text{as } n\to\infty,
        \quad k=1,2.
$$
Together with the convergence of $g_n$ to $g$ in \Kone{}, it remains to estimate $E_n$.
{Since $F_2(v)\in C(\OT)$,  the weak convergence
$h_n\rightharpoonup h$ in $L^2(\mathcal O_T)$ implies
$F_2(v)(h_n-h)\rightharpoonup0
~\text{in }L^2(\mathcal O_T).$
By definition,
$
E_n=\mathcal G\bigl(F_2(v)(h_n-h)\bigr).
$
Since $\mathcal G:L^2(\mathcal O_T)\to C(\OT)$ is compact by
Lemma~\ref{lem:green-convolution}, we have
$\|E_n\|_{C(\OT)}\to 0.$}
{The preceding estimates show that
$\|u_n-\mathscr R_h(v)\|_{C(\OT)}\to0$.  Since $u_n\to v$ in
$C(\OT)$, it follows that $v=\mathscr R_h(v)$; equivalently, $v$
satisfies \eqref{eq:v} with control $h$.}

By uniqueness of the continuous skeleton equation \eqref{eq:v}, which follows from Proposition \ref{prop:volterra-gronwall} and the pointwise kernel bound \eqref{eq:limiting-kernel-pointwise-rho}, we have \(v=\Upsilon(h)\). Since every subsequence of \(\{u_n\}\) has a further subsequence converging to \(\Upsilon(h)\), the whole sequence converges to \(\Upsilon(h)\).
\end{proof}
%By Lemma~\ref{lem:compactness-kernel}, \Kone, \Kthree{} and \Kfour{} give Assumption~\ref{ass:abstract-compact-stability}(i).  Together with \Ktwo,
%Lemma~\ref{lem:weak-to-strong} gives
%Assumption~\ref{ass:abstract-compact-stability}(ii).
Combining Lemmas~\ref{lem:compactness-kernel} and
\ref{lem:weak-to-strong} leads to the following proposition.
\begin{proposition}%[Verification of compactness and weak stability]
\label{prop:kernel-compactness-stability}
Assume that $b$ and $\sigma$ are Lipschitz continuous. {Then under Setting~\ref{set:scalar-spde-realization},  the conditions \Konefour{} imply Assumption~\ref{ass:abstract-compact-stability}.}
\end{proposition}

\subsection{Directional equicontinuity}

%Let $u_i:=\U(h_i)$ and $D(t):=\sup_{0\le r\le t,\ x\in\overline{\mathcal O}}
%|u_1(r,x)-u_2(r,x)|^2$.  Subtracting the two continuous mild equations,
%using the Lipschitz continuity of $b$ and $\sigma$, and applying the
%Cauchy--Schwarz inequality together with
%\eqref{eq:limiting-kernel-pointwise-rho} and
%\eqref{eq:limiting-kernel-global-l2-bound}, gives
%$$
%D(t)\le C_R\|h_1-h_2\|_{L^2(\mathcal O_T)}^2
%+C_R\int_0^t\varrho(t-s)D(s)\ud s,
%\qquad 0\le t\le T.
%$$
%Proposition~\ref{prop:volterra-gronwall} yields the stated estimate.
%\end{proof}}

\begin{proposition}%[Directional consistency from the kernel requirements]
\label{lem:directional-consistency-verification}
Assume that $b$ and $\sigma$ are Lipschitz continuous. 
Under Setting~\ref{set:scalar-spde-realization}, the conditions \Kone{} and
\Kthree{}, together with Assumption~\ref{ass:abstract-control-approx}, imply
Assumption~\ref{ass:abstract-directional-consistency}.
\end{proposition}

\begin{proof}
Let $h_n^0\in \mathcal{X}_n$, $n=1,2,\ldots$, satisfy
$h_n^0\to h$ in $L^2(\mathcal O_T)$.  Let $q_n=Q_nq$ and fix $a_0>0$.
Since both $\{h_n^0\}_{n\ge1}$ and $\{q_n\}_{n\ge1}$ are bounded in
$L^2(\mathcal O_T)$, set
\[
R_0:=2\sup_{n\ge1}\|h_n^0\|_{L^2(\mathcal O_T)}
      +2a_0\sup_{n\ge1}\|q_n\|_{L^2(\mathcal O_T)}<+\infty.
\]
Here we have used Assumption~\ref{ass:abstract-control-approx} to obtain
$q_n\to q$.  Define
$
\Phi_n(a):=\Un(h_n^0+a q_n)(T,x_0).
$
For $a,a'\in[-a_0,a_0]$, Lemma~\ref{lem:lipschitz-strong-consistency}
gives
\[
\sup_{n\ge1}|\Phi_n(a)-\Phi_n(a')|
\le C_{R_0}\sup_{n\ge1}\|q_n\|_{L^2(\mathcal O_T)}|a-a'|.
\]
Thus $\{\Phi_n\}_{n\ge1}$ is equi-Lipschitz, and hence equicontinuous,
on $[-a_0,a_0]$.  This proves
Assumption~\ref{ass:abstract-directional-consistency}.
\end{proof}

\subsection{Local terminal
nondegeneracy}
\label{sec:terminal-nondegeneracy}
In this subsection, we first establish the
differentiability of the scalar terminal observation
$
\mathcal T(h):=\U(h)(T,x_0)
$
 and present the adjoint representation of the directional derivative $D\mathcal{T}(h)[q]$.  The
additional Green-function condition \Kfive{} then guarantees the local terminal
nondegeneracy in Assumption \ref{ass:abstract-terminal-nondegeneracy}.

If $b$ and $\sigma$ are continuously differentiable with bounded derivatives, for any 
$h\in L^2(\mathcal O_T)$, we define
\begin{equation}\label{eq:ah}
        a_h(s,z):=b'(\U(h)(s,z))+\sigma'(\U(h)(s,z))h(s,z),\quad (s,z)\in\OT.
\end{equation}
{We also introduce the forward Volterra operator $\mathcal V_h$
and the backward Volterra operator ${\mathcal V_h^\dagger}$ used in the
adjoint-state representation:}
\begin{align}
\label{eq:forward-volterra-operator}
(\mathcal V_h\phi)(t,x)
&:=
\int_0^t\int_{\mathcal O}
G_{t-s}(x,z)a_h(s,z)\phi(s,z) \ud z \ud s,
\qquad \phi\in C(\OT),\\
\label{eq:backward-volterra-operator}
({\mathcal V_h^\dagger}p)(s,z)
&:=
\int_s^T\int_{\mathcal O}
G_{r-s}(\xi,z)a_h(r,\xi)p(r,\xi)\ud\xi\ud r,
\qquad p\in L^2(\mathcal O_T).
\end{align}

We first present in Lemmas \ref{lem:linear-resolvent} and \ref{lem:backward-volterra-wellposed} resolvent estimates for $\mathcal V_h$ and $\mathcal V_h^\dagger$, respectively.
{
\begin{lemma}%[A resolvent estimate for $\mathcal V_h$]
\label{lem:linear-resolvent}
Assume that $b$ and $\sigma$ are continuously differentiable with bounded derivatives.
Let the conditions \Ktwo--{\Kfour} hold and
 $h\in L^2(\mathcal O_T)$. {Then
$\mathcal V_h:C(\OT)\to C(\OT)$ is a compact linear operator.} Moreover,
for every $f\in C(\OT)$, the equation
\begin{equation}\label{eq:vLn}
v=\mathcal V_hv+f
\end{equation}
admits a unique solution.  Also, for every $R>0$, there exists $C_R>0$
such that for any  $h\in L^2(\mathcal O_T)$ satisfying $\|h\|_{L^2(\mathcal O_T)}\le R$,
\begin{equation}\label{eq:linear-resolvent-estimate}
\|(I-\mathcal V_h)^{-1}f\|_{C(\OT)}
\le C_R\|f\|_{C(\OT)}.
\end{equation}
\end{lemma}

\begin{proof}
Fix $R>0$ and suppose that
$
        \|h\|_{L^2(\mathcal O_T)}\le R.
$ 
Then by the boundedness of $b'$ and $\sigma'$, we have
\begin{equation}\label{eq:boundah}
\|a_h\|_{L^2(\mathcal O_T)}
\le (T|\mathcal O|)^{1/2}\|b'\|_\infty+\|\sigma'\|_\infty\|h\|_{L^2(\mathcal O_T)}
\le \widetilde C_R,
\end{equation}
where $\|\cdot\|_\infty$ denotes the supremum norm.
{The multiplication operator
$$
M_{a_h}:C(\OT)\to L^2(\mathcal O_T),\qquad M_{a_h}v:=a_hv,
$$
is bounded by \eqref{eq:boundah}. Since
$\mathcal V_h=\mathcal G\circ M_{a_h}$ and $\mathcal G$ is compact by
Lemma~\ref{lem:green-convolution}, $\mathcal V_h$ is a compact linear
operator from $C(\OT)$ to $C(\OT)$.}

We next prove the existence and uniqueness of the solution to \eqref{eq:vLn} by a contraction argument.  For $\lambda>0$, equip $C(\OT)$ with the  norm
$$
        \|v\|_\lambda
        :=
        \sup_{(t,x)\in\OT}e^{-\lambda t}|v(t,x)|,\quad v\in C(\OT).
$$
This norm is equivalent to the usual norm $ \|\cdot\|_{C(\OT)}$.
Set
$
        A_1:=b', A_2:=\sigma',
         h^1\equiv1, h^2:=h
$ throughout this proof. Then by \eqref{eq:forward-volterra-operator}, for $v\in C(\OT)$, we can reformulate
$$
        \mathcal V_hv(t,x)
        =\sum_{k=1}^2\int_0^t\int_{\mathcal O}G_{t-s}(x,z)
        A_k(u(s,z))h^k(s,z)v(s,z)\ud z\ud s,
$$
where $u:=\U(h)$.
For $k=1,2$, the Cauchy--Schwarz inequality and
\eqref{eq:limiting-kernel-pointwise-rho} yield
\begin{align*}
&\;e^{-\lambda t}\biggl|\int_0^t\int_{\mathcal O}G_{t-s}(x,z)
A_k(u(s,z))h^k(s,z)v(s,z)\ud z\ud s\biggr|\\
\le&\; C\|h^k\|_{L^2(\mathcal O_T)}\|v\|_\lambda
\left(\int_0^T e^{-2\lambda r}\varrho(r)\ud r\right)^{1/2}.
\end{align*}
Consequently, for any $v\in C(\OT)$,
\begin{equation}\label{eq:Lhv}
        \|\mathcal V_hv\|_\lambda
        \le \kappa_R(\lambda)\|v\|_\lambda,\quad \kappa_R(\lambda):=
C_R\left(\int_0^T e^{-2\lambda r}\varrho(r)\ud r\right)^{1/2}.
\end{equation}
Since $\varrho\in L^1(0,T)$, the dominated convergence theorem implies
$
        \kappa_R(\lambda)\to0$ as $\lambda\to+\infty$.
We may therefore choose $\lambda_R$ sufficiently large that
$\kappa_R(\lambda_R)<1$.
Thus, for each $f\in C(\OT)$, the map
$
        v\mapsto \mathcal V_hv+f
$
is a contraction on $C(\OT)$ endowed with $\|\cdot\|_{\lambda_R}$.
The Banach fixed-point theorem therefore yields a unique solution to the equation \eqref{eq:vLn}.
%Thus $I-\mathcal V_h$ is bijective.
  By the equivalence of the  norm $\|\cdot\|_{\lambda}$ and the
usual uniform norm $\|\cdot\|_{C(\OT)}$,
the estimate \eqref{eq:Lhv} with $\lambda=\lambda_R$ proves \eqref{eq:linear-resolvent-estimate} with \(C_R=\frac{e^{\lambda_RT}}
{1-\kappa_R(\lambda_R)}\), for any $h\in L^2(\mathcal O_T)$ with $\|h\|_{ L^2(\mathcal O_T)}\le R$.
%\begin{equation}\label{eq:linear-resolvent-estimate}
%        \|{(I-\mathcal V_h)^{-1}} f\|_{C(\OT)}
%      \le
%C_R\|f\|_{C(\OT)}.
%\end{equation}
\end{proof}

{
\begin{lemma}%[A resolvent estimate for the backward Volterra operator]
\label{lem:backward-volterra-wellposed}
{
Assume that $b$ and $\sigma$ are continuously differentiable with bounded derivatives.
Let the conditions \Ktwo--{\Kthree} hold and
 $h\in L^2(\mathcal O_T)$. Then for every $f\in L^2(\mathcal O_T)$, the equation
\begin{equation}\label{eq:general-backward-volterra}
 p=f+{\mathcal V_h^\dagger}p
\end{equation}
has a unique solution $p\in L^2(\mathcal O_T)$.  
Moreover, for every $R>0$, there exists $C_R>0$ such that for any  $h\in L^2(\mathcal O_T)$ satisfying $\|h\|_{L^2(\mathcal O_T)}\le R$,
\begin{equation}\label{eq:ILhdag}
\|(I-\mathcal V_h^\dagger)^{-1}f\|_{L^2(\mathcal O_T)}
\le C_R\|f\|_{L^2(\mathcal O_T)}.
\end{equation}
}
\end{lemma}

\begin{proof}
{
For $\lambda>0$, equip $L^2(\mathcal O_T)$ with the  norm
$$
\tnorm{p}_\lambda
:=
\left(
\int_0^T\int_{\mathcal O} e^{2\lambda s}|p(s,z)|^2\ud z\ud s
\right)^{1/2}.
$$
This norm is equivalent to the usual $L^2(\mathcal O_T)$ norm.  For almost every
$s\in[0,T]$, Minkowski's inequality, the definition of ${\mathcal V_h^\dagger}p$ in \eqref{eq:backward-volterra-operator},
\eqref{eq:limiting-kernel-pointwise-rho}, and the Cauchy--Schwarz inequality
yield
\begin{align}
e^{\lambda s}
\|({\mathcal V_h^\dagger}p)(s,\cdot)\|_{L^2(\mathcal O)}
&\le
\int_s^T e^{-\lambda(r-s)}
\left\|
\int_{\mathcal O}G_{r-s}(\xi,\cdot)a_h(r,\xi)p(r,\xi)\ud\xi
\right\|_{L^2(\mathcal O)}
e^{\lambda r}\ud r
\notag\\
&\quad\le
\int_s^T e^{-\lambda(r-s)}\varrho(r-s)^{1/2}\|a_h(r,\cdot)\|_{L^2(\mathcal O)}
e^{\lambda r}\|p(r,\cdot)\|_{L^2(\mathcal O)}\ud r.
\label{eq:backward-pointwise-weighted}
\end{align}
Define
\[
\begin{aligned}
\check k_\lambda(r)&:=\mathbf 1_{[-T,0]}(r)
        e^{\lambda r}\varrho(-r)^{1/2},\\
F_{h,p}(r)&:=e^{\lambda r}
        \|a_h(r,\cdot)\|_{L^2(\mathcal O)}
        \|p(r,\cdot)\|_{L^2(\mathcal O)}.
\end{aligned}
\]
Then inequality~\eqref{eq:backward-pointwise-weighted} becomes
$$
e^{\lambda s}\|({\mathcal V_h^\dagger}p)(s,\cdot)\|_{L^2(\mathcal O)}
\le (\check k_\lambda*F_{h,p})(s),
$$
{where $F_{h,p}$ is extended by zero outside $[0,T]$.} Thus, applying Young's
convolution inequality, the Cauchy--Schwarz inequality, and \eqref{eq:boundah} yields
\begin{align}\label{eq:backward-weighted-contraction}
\tnorm{{\mathcal V_h^\dagger}p}_\lambda
&\le \|\check k_\lambda\|_{L^2(\mathbb R)}\|F_{h,p}\|_{L^1(0,T)}\\\notag
&\le
\left(\int_0^T e^{-2\lambda r}\varrho(r)\ud r\right)^{1/2}
\left(\int_0^T\int_{\mathcal O}e^{2\lambda r}|p(r,z)|^2\ud z\ud r\right)^{1/2}\|a_h\|_{L^2(\mathcal O_T)}  \\\notag
&=
\widetilde{C}_R
\left(\int_0^T e^{-2\lambda r}\varrho(r)\ud r\right)^{1/2}
\tnorm{p}_\lambda.
\end{align}
Since $\varrho\in L^1(0,T)$, the dominated convergence theorem gives
$
\int_0^T e^{-2\lambda r}\varrho(r)\ud r\to0$ as $\lambda\to+\infty$.
Hence, for a fixed $R>0$, one can choose $\lambda'_R>0$ such that
$
\kappa_R'
:=
\widetilde C_R\sqrt{\int_0^T e^{-2\lambda'_R r}\varrho(r)\ud r}
<1,
$ which together with \eqref{eq:backward-weighted-contraction} implies that 
\begin{equation}\label{eq:Lhdp}
\tnorm{{\mathcal V_h^\dagger}p}_{\lambda'_R}\le
\kappa_R'\tnorm{p}_{\lambda'_R}.
\end{equation}
Hence, for any $\|h\|_{L^2(\mathcal O_T)}\le R$, the map
$
p\mapsto f+{\mathcal V_h^\dagger}p
$
is a contraction on $L^2(\mathcal O_T)$ endowed with the norm $\tnorm{\cdot}_{\lambda'_R}$.
The Banach fixed-point theorem gives a unique solution of
\eqref{eq:general-backward-volterra}.  
%Thus $I-\mathcal V_h^\dagger$ is bijective.
Finally, the equivalence between the norms $\tnorm{\cdot}_{\lambda'_R}$ and $\|\cdot\|_{L^2(\mathcal O_T)}$, along with \eqref{eq:Lhdp}, yields \eqref{eq:ILhdag}.
}
\end{proof}}
Define the inhomogeneous
control operator $\mathcal B_h$ by
\begin{equation}\label{eq:inhomogeneous-control-operator}
(\mathcal B_hq)(t,x)
:=
\int_0^t\int_{\mathcal O}
G_{t-s}(x,z)\sigma(\U(h)(s,z))q(s,z) \ud z \ud s,
\qquad q\in L^2(\mathcal O_T).
\end{equation}
{Since $\sigma(\U(h))\in C(\OT)$, the multiplication
operator
$$
M_{\sigma(\U(h))}:L^2(\mathcal O_T)\to L^2(\mathcal O_T),
\qquad M_{\sigma(\U(h))}q:=\sigma(\U(h))q,
$$
is bounded. Hence
$\mathcal B_h=\mathcal G\circ M_{\sigma(\U(h))}$ is a compact linear
operator from $L^2(\mathcal O_T)$ to $C(\OT)$ by
Lemma~\ref{lem:green-convolution}.}
%
%
%\begin{corollary}[The linearized skeleton map]
%\label{lem:linearized-skeleton-estimate}
%Let $h,q\in L^2(\mathcal O_T)$. We denote by $\mathcal{V}_h q:={(I-\mathcal V_h)^{-1}}(\mathcal B_h q)$ the linearized skeleton map. Then for every $R>0$, there exists $C_R>0$ such that for any $\|h\|_{L^2(\mathcal O_T)}\le R$,
%\begin{equation}\label{eq:vq}
%        \|{{(I-\mathcal V_h)^{-1}}\circ \mathcal B_h}q\|_{C(\OT)}
%        \le C_R\|q\|_{L^2(\mathcal O_T)}
%\end{equation}
%\end{corollary}

%\begin{proof}
%The global $L^2$-bound \eqref{eq:limiting-kernel-global-l2-bound},
%the local bound in Corollary~\ref{lem:continuous-skeleton-local-bounds}(i),
%and the boundedness of $\sigma'$ imply
%$$
%        \|\mathcal B_hq\|_{C(\OT)}\le C_R\|q\|_{L^2(\mathcal O_T)}.
%$$
%The assertion \eqref{eq:vq} follows by applying Lemma~\ref{lem:linear-resolvent} with
%$f=\mathcal B_hq$.  
%\end{proof}

}
\begin{lemma}[$C^1$-regularity of the skeleton map]
\label{lem:skeleton-C1}
Assume that $b$ and $\sigma$ are continuously differentiable with bounded derivatives.
Let the conditions \Ktwo--{\Kfour} hold. Then
$
        \U:L^2(\mathcal O_T)\to C(\OT)
$
is continuously Fréchet differentiable, and for any
 $h\in L^2(\mathcal O_T)$,
$$
        D\U(h)={(I-\mathcal V_h)^{-1}}\mathcal B_h
        \quad\text{in }\mathcal L\bigl(L^2(\mathcal O_T),C(\OT)\bigr).
$$
Equivalently, for $h,q\in L^2(\mathcal O_T)$,
$D\U(h)[q]={(I-\mathcal V_h)^{-1}}(\mathcal B_hq)$. 
{Moreover, for every $h\in L^2(\mathcal O_T)$,
$D\U(h):L^2(\mathcal O_T)\to C(\OT)$ is compact.}
\end{lemma}

{
\begin{proof}

Fix $h\in L^2(\mathcal O_T)$ and let $R:=\|h\|_{L^2(\mathcal O_T)}$.  We first prove that $\U$ is
Fr\'echet differentiable at $h$.  {Define the candidate derivative
$A_h:={(I-\mathcal V_h)^{-1}}\mathcal B_h.$
The operator $\mathcal B_h:L^2(\mathcal O_T)\to C(\OT)$ is compact, while
$(I-\mathcal V_h)^{-1}:C(\OT)\to C(\OT)$ is bounded by
Lemma~\ref{lem:linear-resolvent}.  Hence
$A_h:L^2(\mathcal O_T)\to C(\OT)$ is a compact linear operator and, in particular,
$A_h\in\mathcal L\bigl(L^2(\mathcal O_T),C(\OT)\bigr)$.}

Set
$$
     u:=\U(h),\qquad   \widetilde{u}_q:=\U(h+q),
        \qquad R_q:=\widetilde{u}_q-u-A_hq,
$$
where $q\in \mathbb{B}_1:=\{q\in L^2(\mathcal O_T):\|q\|_{L^2(\mathcal O_T)}\le1\}$.  
Denote $\mathcal E_q:=(I-\mathcal V_h)R_q$. Then the definition yields $\mathcal E_q=(I-\mathcal V_h)(\widetilde{u}_q-u)-\mathcal B_hq$. 
By \eqref{eq:continuous-skeleton},   \eqref{eq:forward-volterra-operator}, and \eqref{eq:inhomogeneous-control-operator}, we obtain 
$$
\begin{aligned}
\mathcal E_q(t,x)
={}&\int_0^t\int_{\mathcal O}G_{t-s}(x,z)
\bigl[b(\widetilde{u}_q)-b(u)-b'(u)(\widetilde{u}_q-u)\bigr](s,z)\ud z\ud s\\
&+\int_0^t\int_{\mathcal O}G_{t-s}(x,z)
\bigl[\sigma(\widetilde{u}_q)-\sigma(u)-\sigma'(u)(\widetilde{u}_q-u)\bigr](s,z)
h(s,z)\ud z\ud s\\
&+\int_0^t\int_{\mathcal O}G_{t-s}(x,z)
\bigl[\sigma(\widetilde{u}_q)-\sigma(u)\bigr](s,z)q(s,z)\ud z\ud s.
\end{aligned}
$$
Further,
applying Lemma \ref{lem:linear-resolvent} leads to
\begin{equation}\label{eq:Rq}
        \|R_q\|_{C(\OT)}=\|(I-\mathcal V_h)^{-1}\mathcal{E}_q\|_{C(\OT)}\le C_R\|\mathcal E_q\|_{C(\OT)}.
\end{equation}
{By Corollary~\ref{lem:continuous-skeleton-local-bounds}(i),
there exists $M_R>0$ such that
$\|u\|_{C(\OT)}+\|\widetilde{u}_q\|_{C(\OT)}\le M_R$ for all
$q\in \mathbb{B}_1$.  In addition, Corollary~\ref{lem:continuous-skeleton-local-bounds}(ii)
gives}
\begin{equation}\label{eq:skeleton-local-lipschitz}
        \delta_q:=\|\widetilde{u}_q-u\|_{C(\OT)}
        \le C_R\|q\|_{L^2(\mathcal O_T)}.
\end{equation}
Since $b$ and $\sigma$ are continuously differentiable with bounded derivatives, we have 
\begin{equation}\label{omegahr}
\omega_R(r):=
\sup_{\substack{|a|,|c|\le M_R,\,|a-c|\le r}}
\bigl(|b'(a)-b'(c)|+|\sigma'(a)-\sigma'(c)|\bigr)\to 0,\quad \text{as }r\downarrow0,
\end{equation}
and by Taylor's formula,
$$
\max\bigl\{|b(\widetilde{u}_q)-b(u)-b'(u)(\widetilde{u}_q-u)|,|\sigma(\widetilde{u}_q)-\sigma(u)-\sigma'(u)(\widetilde{u}_q-u)|\bigr\}
\le \omega_R(\delta_q)|\widetilde{u}_q-u|.
$$
The global $L^2$-bound \eqref{eq:limiting-kernel-global-l2-bound} on $G$,
the boundedness of $\sigma'$, and
\eqref{eq:skeleton-local-lipschitz} yield
\begin{equation}\label{eq:remainder-source-estimate}
\begin{aligned}
\|\mathcal E_q\|_{C(\OT)}
&\le C_R\bigl(\omega_R(\delta_q)\delta_q
        +\delta_q\|q\|_{L^2(\mathcal O_T)}\bigr)\\
&\le C_R\bigl(\omega_R(C_R\|q\|_{L^2(\mathcal O_T)})
        +\|q\|_{L^2(\mathcal O_T)}\bigr)\|q\|_{L^2(\mathcal O_T)}.
\end{aligned}
\end{equation}
Consequently, we deduce from \eqref{eq:Rq}, \eqref{eq:remainder-source-estimate}, and \eqref{omegahr} that
{
$$
\frac{\|\U(h+q)-\U(h)-A_hq\|_{C(\OT)}}
     {\|q\|_{L^2(\mathcal O_T)}}=\frac{\|R_q\|_{C(\OT)}}
     {\|q\|_{L^2(\mathcal O_T)}}\to0
\qquad\text{as }\|q\|_{L^2(\mathcal O_T)}\to0.
$$
Since $A_h\in\mathcal L\bigl(L^2(\mathcal O_T),C(\OT)\bigr)$,  $\U$ is Fr\'echet differentiable at $h$, with
$
        D\U(h)=A_h={(I-\mathcal V_h)^{-1}}\mathcal B_h.
$
In particular, $D\U(h)$ is compact.}

It remains to prove continuity of $D\U$.  Let $h_i\to h$ in
$L^2(\mathcal O_T)$. For $q\in L^2(\mathcal O_T)$, we denote
$$
\begin{aligned}
        u_i&:=\U(h_i),\qquad u:=\U(h),\\
        v^i_q&:=D\U(h_i)[q],\qquad v_q:=D\U(h)[q],\qquad
        w^i_q:=v^i_q-v_q.
\end{aligned}
$$
Subtracting the two linearized equations for $v^i_q$ and $v_q$ yields
$
        w^i_q=\mathcal{V}_{h_i}w^i_q+\mathcal R_i(q),
$
where
\begin{align}\label{eq:Riq}
\mathcal R_i(q)(t,x)
={}&\int_0^t\int_{\mathcal O}G_{t-s}(x,z)
\bigl[b'(u_i)-b'(u)\bigr](s,z)v_q(s,z)\ud z\ud s\\\notag
&+\int_0^t\int_{\mathcal O}G_{t-s}(x,z)
\bigl[\sigma'(u_i)-\sigma'(u)\bigr](s,z)h_i(s,z)v_q(s,z)\ud z\ud s\\\notag
&+\int_0^t\int_{\mathcal O}G_{t-s}(x,z)
\sigma'(u(s,z))(h_i-h)(s,z)v_q(s,z)\ud z\ud s\\\notag
&+\int_0^t\int_{\mathcal O}G_{t-s}(x,z)
\bigl[\sigma(u_i)-\sigma(u)\bigr](s,z)q(s,z)\ud z\ud s.
\end{align}
Choose $R>0$ such that
$
        \sup_i\|h_i\|_{L^2(\mathcal O_T)}+\|h\|_{L^2(\mathcal O_T)}\le R.
$
Since $h_i\to h$ in
$L^2(\mathcal O_T)$, {Corollary~\ref{lem:continuous-skeleton-local-bounds}(ii)
implies that}
$u_i\to u$ in $C(\OT)$ as $i\to\infty$.  Therefore, by the continuity of $b'$, $\sigma'$ and $\sigma$,
\begin{equation}\label{eq:varbsi}
\begin{aligned}
\varepsilon^i_{b,\sigma}:={}&\|b'(u_i)-b'(u)\|_{C(\OT)}
 +\|\sigma'(u_i)-\sigma'(u)\|_{C(\OT)}\\
&+\|\sigma(u_i)-\sigma(u)\|_{C(\OT)}
\to0,\quad \text{as }i\to\infty.
\end{aligned}
\end{equation}
The global $L^2$-bound \eqref{eq:limiting-kernel-global-l2-bound} on $G$,
{the local bound in
Corollary~\ref{lem:continuous-skeleton-local-bounds}(i)}, and the boundedness
of $\sigma'$ imply
$
        \|\mathcal B_hq\|_{C(\OT)}\le C_R\|q\|_{L^2(\mathcal O_T)},
$ where $\mathcal{B}_hq$ is defined in \eqref{eq:inhomogeneous-control-operator}.
Then applying the estimate \eqref{eq:linear-resolvent-estimate} with
$f=\mathcal B_hq$ leads to
\begin{equation}\label{eq:vq}
       \|v_q\|_{C(\OT)}= \|{(I-\mathcal V_h)^{-1}}(\mathcal B_hq)\|_{C(\OT)}
        \le C_R\|q\|_{L^2(\mathcal O_T)}.
\end{equation}
By \eqref{eq:Riq}, \eqref{eq:vq}, the Cauchy--Schwarz inequality, and 
\eqref{eq:limiting-kernel-global-l2-bound}, we arrive at 
\begin{equation}\label{eq:RiqC}
        \|\mathcal R_i(q)\|_{C(\OT)}
        \le C_R\bigl(\varepsilon^i_{b,\sigma}+\|h_i-h\|_{L^2(\mathcal O_T)}\bigr)
        \|q\|_{L^2(\mathcal O_T)}.
\end{equation}
Since $
        w^i_q=\mathcal V_{h_i}w^i_q+\mathcal R_i(q),
$
applying \eqref{eq:linear-resolvent-estimate}, with $h_i$ in place
of $h$, yields
$$
\|D\U(h_i)[q]-D\U(h)[q]\|_{C(\OT)}=\|w_q^i\|_{C(\OT)}
\le C_R\bigl(\varepsilon^i_{b,\sigma}+\|h_i-h\|_{L^2(\mathcal O_T)}\bigr)
\|q\|_{L^2(\mathcal O_T)},
$$
where the last step is due to \eqref{eq:RiqC}.
Taking the supremum over $\|q\|_{L^2(\mathcal O_T)}\le1$, it follows from $h_i\to h$ in $L^2(\mathcal O_T)$ and \eqref{eq:varbsi} that
\[
        \lim_{i\to\infty}\|D\U(h_i)-D\U(h)\|_{\mathcal L(L^2(\mathcal O_T),C(\OT))}=0.
\]
This proves the continuity of $D\Upsilon$.
\end{proof}
}

The following lemma gives the adjoint representation for the terminal derivative.
\begin{lemma}
\label{lem:adjoint-representation}
Assume that $b$ and $\sigma$ are continuously differentiable with bounded derivatives.
Let the conditions \Ktwo--{\Kfour} hold.
{Define the terminal observation map
$
\mathcal T:L^2(\mathcal O_T)\to\R
$
by}
$
{
\mathcal T(h):=\U(h)(T,x_0).
}
$
{Then $\mathcal T$ is continuously Fr\'echet differentiable.
Moreover, for every $h,q\in L^2(\mathcal O_T)$,}
\begin{equation}\label{eq:terminal-derivative-representation}
        D\mathcal T(h)[q]
        =
        \int_0^T\int_{\mathcal O}
        \sigma(\U(h)(s,z))p_h(s,z)q(s,z) \ud z \ud s,
\end{equation}
where $p_h\in L^2(\mathcal O_T)$ is the unique solution to the backward mild equation
\begin{equation}\label{eq:adjoint-kernel}
\begin{aligned}
p_h(s,z)
=
G_{T-s}(x_0,z)
&+\int_s^T\int_{\mathcal O}
G_{r-s}(\xi,z)a_h(r,\xi)p_h(r,\xi) \ud \xi \ud r.
\end{aligned}
\end{equation}
%where $a_h$ is defined in \eqref{eq:ah}.
\end{lemma}

\begin{proof}
{
Let $u:=\U(h)$ and put
$
        g_{T,x_0}(s,z):=G_{T-s}(x_0,z)
$ in this proof.
Then $g_{T,x_0}\in L^2(\mathcal O_T)$ by
\eqref{eq:limiting-kernel-global-l2-bound}.
For $q\in L^2(\mathcal O_T)$, by Lemma~\ref{lem:skeleton-C1},
$$
        v_q:=D\U(h)[q]={(I-\mathcal V_h)^{-1}}(\mathcal B_hq),
$$
or, equivalently,
$
        v_q=\mathcal V_h v_q+\mathcal B_hq.
$
 Let $\mathscr E_{T,x_0}:C(\OT)\to\R$ denote the terminal evaluation map defined by $\mathscr E_{T,x_0}(f)=f(T,x_0)$.
Since $\mathcal T=\mathscr E_{T,x_0}\circ\U$ and
$\mathscr E_{T,x_0}$ is a continuous linear functional, the chain rule gives
$$
        D\mathcal T(h)[q]
        =\mathscr E_{T,x_0}(D\U(h)[q])
        =v_q(T,x_0).
$$
{Since $D\U$ is continuous by
Lemma~\ref{lem:skeleton-C1} and post-composition with the fixed continuous
linear functional $\mathscr E_{T,x_0}$ is continuous, the identity
$D\mathcal T(h)=\mathscr E_{T,x_0}\circ D\U(h)$ implies that
$D\mathcal T$ is continuous. Hence $\mathcal T\in C^1(L^2(\mathcal O_T),\R)$.}

{
Set $\phi_q:=a_hv_q+\sigma(u)q$.  Since
$u,v_q\in C(\OT)$, $h,q\in L^2(\mathcal O_T)$, and
$b'$ and $\sigma'$ are bounded, we have $\phi_q\in L^2(\mathcal O_T)$.
Moreover, the identity
$v_q=\mathcal V_hv_q+\mathcal B_hq$, together with
\eqref{eq:forward-volterra-operator},
\eqref{eq:inhomogeneous-control-operator}, and \eqref{eq:G}, gives
$v_q=\mathcal G\phi_q$; Lemma~\ref{lem:green-convolution} ensures that
$\mathcal G\phi_q\in C(\OT)$.  Consequently, the
definition of $\phi_q$ yields
\begin{align} \label{eq1}
	\phi_q-a_h\mathcal G\phi_q=\sigma(u)q.
\end{align}        

Moreover, evaluating $v_q=\mathcal G\phi_q$ at $(T,x_0)$ gives
\begin{equation}\label{eq:terminal-evaluation-linearized}
v_q(T,x_0)
 =\langle g_{T,x_0},\phi_q\rangle_{L^2(\mathcal O_T)}.
\end{equation}
To rewrite \(D\mathcal T(h)[q]\) directly in terms of \(q\), without the auxiliary variable \(\phi_q\), we seek a
backward variable $p_h$ such that
\begin{equation}\label{eq:adjoint-cancellation-target}
\langle g_{T,x_0},\phi\rangle_{L^2(\mathcal O_T)}
=
\langle p_h,\phi-a_h\mathcal G\phi\rangle_{L^2(\mathcal O_T)}
\qquad\text{for every }\phi\in L^2(\mathcal O_T).
\end{equation}
Indeed, taking $\phi=\phi_q$ would then give the desired representation in view of \eqref{eq1}.
The following duality identity determines precisely how $p_h$ should be
chosen.}
%For $p,\phi\in L^2(\mathcal O_T)$, the Cauchy--Schwarz inequality and
%\eqref{eq:limiting-kernel-global-l2-bound} yield, for every $(r,\xi)\in\OT$,
%$$
%\int_0^r\int_{\mathcal O}
%|G_{r-s}(\xi,z)\phi(s,z)|\ud z\ud s
%\le
%\left(2\int_0^T\varrho(\theta)\ud\theta\right)^{1/2}
%\|\phi\|_{L^2(\mathcal O_T)}.
%$$
For $p,\phi\in L^2(\mathcal O_T)$, we claim
\begin{equation}\label{eq:forward-backward-duality}
\langle {\mathcal V_h^\dagger}p,\phi\rangle_{L^2(\mathcal O_T)}
=
\langle a_hp,\mathcal G\phi\rangle_{L^2(\mathcal O_T)}.
\end{equation}
Indeed, the triple integral in $\langle {\mathcal V_h^\dagger}p,\phi\rangle_{L^2(\mathcal O_T)}$ is
absolutely integrable. More precisely, by the Cauchy--Schwarz inequality and
\eqref{eq:limiting-kernel-global-l2-bound},
\begin{align*}
&\int_0^T\int_{\mathcal O}|a_h(r,\xi)p(r,\xi)|
\left(\int_0^r\int_{\mathcal O}
|G_{r-s}(\xi,z)\phi(s,z)|\ud z\ud s\right)\ud\xi\ud r\\
&\qquad\le
\left(2\int_0^T\varrho(\theta)\ud\theta\right)^{1/2}
\|a_h\|_{L^2(\mathcal O_T)}\|p\|_{L^2(\mathcal O_T)}\|\phi\|_{L^2(\mathcal O_T)}
<+\infty.
\end{align*}
Thus Fubini's theorem applies and yields
\eqref{eq:forward-backward-duality}.
{In particular, \eqref{eq:forward-backward-duality} implies that
for every $p,\phi\in L^2(\mathcal O_T)$,
$$
\langle p-\mathcal V_h^\dagger p,\phi\rangle_{L^2(\mathcal O_T)}
=
\langle p,\phi-a_h\mathcal G\phi\rangle_{L^2(\mathcal O_T)}.
$$
Consequently, \eqref{eq:adjoint-cancellation-target} holds if $p_h$ is
chosen to satisfy
$$
        p_h-\mathcal V_h^\dagger p_h=g_{T,x_0},
$$
which is precisely
\eqref{eq:adjoint-kernel}. 
Since $g_{T,x_0}\in L^2(\mathcal O_T)$,
Lemma~\ref{lem:backward-volterra-wellposed} gives a unique  $p_h\in L^2(\mathcal O_T)$ with this property.  Accordingly, the $p_h$ defined by \eqref{eq:adjoint-kernel} satisfies \eqref{eq:adjoint-cancellation-target}.  Taking $\phi=\phi_q$ in
\eqref{eq:adjoint-cancellation-target}
and using \eqref{eq1}--\eqref{eq:terminal-evaluation-linearized}, we obtain
$$
\begin{aligned}
D\mathcal T(h)[q]
&=v_q(T,x_0)
=\langle g_{T,x_0},\phi_q\rangle_{L^2(\mathcal O_T)}\\
&=\langle p_h,\phi_q-a_h\mathcal G\phi_q\rangle_{L^2(\mathcal O_T)}
=\langle p_h,\sigma(u)q\rangle_{L^2(\mathcal O_T)}.
\end{aligned}
$$
Since $u=\U(h)$, this is \eqref{eq:terminal-derivative-representation}.}
This completes the proof.}
\end{proof}

\begin{proposition}[Local terminal nondegeneracy]
\label{prop:terminal-nondegeneracy}
Assume that $b$ and $\sigma$ are continuously differentiable with bounded derivatives, and $\sigma(r)\ne0$ for every $r\in\R$. Then under the conditions \Kone--{\Kfive}, for every $h\in L^2(\mathcal O_T)$, there exists
$q_h\in L^2(\mathcal O_T)$
such that
$
D\mathcal T(h)[q_h]\ne0.
$
\end{proposition}

\begin{proof}
Let $p_h\in L^2(\mathcal O_T)$ solve
$
p_h=G_{T-\cdot}(x_0,\cdot)+\mathcal V_h^\dagger p_h.
$
If $\|p_h\|_{L^2(\mathcal O_T)}=0$, then
$
G_{T-\cdot}(x_0,\cdot)=0
~\text{in }L^2(\mathcal O_T),
$
contradicting \Kfive{}.  Hence
{$\|p_h\|_{L^2(\mathcal O_T)}>0$.}  Since $\sigma(u(s,z))\ne0$ for every $(s,z)\in\OT$,
we have {$\|\sigma(u)p_h\|_{L^2(\mathcal O_T)}>0$}.  By
\eqref{eq:terminal-derivative-representation},
$$
D\mathcal T(h)[q]=\langle\sigma(u)p_h,q\rangle_{L^2(\mathcal O_T)},
$$
where $u:=\U(h)$. { Taking
$
q_h:=\sigma(u)p_h
$
gives
$
D\mathcal T(h)[q_h]=\|\sigma(u)p_h\|_{L^2(\mathcal O_T)}^2>0.
$
}
This completes the proof.
\end{proof}

Hereafter, $C_b^1(\mathbb R)$ denotes the class of functions from $\mathbb R$ to $\mathbb R$  which are continuously differentiable with bounded first-order derivatives. 
\begin{theorem}[Diagonal convergence criterion]
\label{thm:closed-diagonal}
{Fix a diagonal refinement sequence $\{(n,m_n)\}_{n\ge1}$ with
$m_n\to+\infty$.}  Assume that $b,\sigma\in C_b^1(\mathbb R)$ and $\sigma(r)\ne0$ for every
$r\in\R$. Under Setting~\ref{set:scalar-spde-realization}, if the conditions
\Kone--{\Kfive} {hold along $\{(n,m_n)\}_{n\ge1}$} and
Assumption~\ref{ass:abstract-control-approx} holds, then for every $y\in\R$,
$$
        I^{n,m_n}(y)\to I(y)
        \qquad\text{as }n\to\infty.
$$
In particular, for every $y\in\R$ with $I(y)<+\infty$, the conclusions of
Corollary~\ref{cor:subsequential-minimizers} hold along the same diagonal
refinement sequence $\{(n,m_n)\}_{n\ge1}$.
\end{theorem}

\begin{proof}
{Under Setting~\ref{set:scalar-spde-realization},
Proposition~\ref{prop:kernel-compactness-stability} verifies
Assumption~\ref{ass:abstract-compact-stability}.
Proposition~\ref{lem:directional-consistency-verification} verifies
Assumption~\ref{ass:abstract-directional-consistency}.
Lemma~\ref{lem:adjoint-representation} shows that the terminal observation map
$\mathcal T$ is continuously Fr\'echet differentiable, while
Proposition~\ref{prop:terminal-nondegeneracy} provides, for every
$h\in L^2(\mathcal O_T)$, a direction $q_h\in L^2(\mathcal O_T)$ such that
$D\mathcal T(h)[q_h]\ne0$.  Hence
Assumption~\ref{ass:abstract-terminal-nondegeneracy} holds.
Together with Assumption~\ref{ass:abstract-control-approx}, all hypotheses of
Theorem~\ref{thm:abstract-framework} hold, and the conclusion follows.}
\end{proof}

\section{Applications to prototype SPDEs and fully discrete schemes}
\label{sec:applications}

We now apply Theorem~\ref{thm:closed-diagonal} in the SPDE setting of
Section~\ref{sec:spde-notation}.  Sections~\ref{sec:K5-prototype-operators}
and~\ref{subsec:grid-controls-LDP} verify the kernel condition \Kfive{} and  Assumption~\ref{ass:abstract-control-approx}, respectively.  Section~\ref{subsec:concrete-fully-discrete-schemes}
states the two fully discrete applications, while the scheme-specific
Green-function estimates and fixed-grid LDP proofs are given in
Appendices~\ref{app:kernel-verifications} and~\ref{app:fixed-grid-LDPs}.

\subsection{Continuous kernels and terminal nondegeneracy}
\label{sec:K5-prototype-operators}
In applications, the continuous Green function needs to satisfy the integrable condition \eqref{eq:limiting-kernel-pointwise-rho}, which imposes the dimensional restrictions 
summarized in Table~\ref{tab:prototype-kernel-dimensions}.

\begin{table}[H]
	\centering
	\small
		\caption{Continuous integrable kernel bounds and dimensional restrictions.}
		\label{tab:prototype-kernel-dimensions}
		\setlength{\tabcolsep}{4pt}
		\begin{tabular}{cccc}
		\toprule
		Operator
		&
		Kernel bound
		&
		Condition ensuring $\varrho\in L^1(0,T)$
		&
		Admissible dimensions
		\\
		\midrule
		$\partial_t-\Delta$
		&
		$\varrho(r)=Cr^{-d/2}$
		&
		$d<2$
		&
		$d=1$
		\\
		$\partial_{tt}-\Delta$
		&
		$\varrho(r)=C$
		&
		$d<2$ 
		&
		$d=1$
		\\
		$\partial_t+\Delta^2$
		&
		$\varrho(r)=Cr^{-d/4}$
		&
		$d<4$
		&
		$d=1,2,3$
		\\
		\bottomrule
	\end{tabular}
\end{table}

The estimates in Table~\ref{tab:prototype-kernel-dimensions} are a common
consequence of the diagonal Dirichlet heat-kernel bound.  
Indeed, set $A=-\Delta$.  By the spectral theory of symmetric elliptic
Dirichlet operators \cite[Section~6.5.1, Theorem~1]{Evans2010PDE}, $A$
admits a sequence of eigenvalues
$
0<\lambda_1\le\lambda_2\le\cdots,
~ \lambda_k\to+\infty,
$
and a corresponding orthonormal basis $\{\varphi_k\}_{k\ge1}$ of
$L^2(\mathcal O)$ satisfying $A\varphi_k=\lambda_k\varphi_k$.  We list
the eigenvalues according to multiplicity.  With this notation, the
spectral representation of $e^{-sA}$ gives its integral kernel as
$$
H_s(x,z)=\sum_{k\ge1}e^{-s\lambda_k}\varphi_k(x)\varphi_k(z),
\qquad s>0.
$$
For a bounded domain with sufficiently regular boundary,
\begin{equation}\label{eq:diagonal-heat-kernel-bound}
	\sup_{x\in\mathcal O}H_s(x,x)\le Cs^{-d/2},
	\qquad 0<s\le T;
\end{equation}
see, for example, \cite[Chapter~3]{Davies1989Heat}.  For the heat equation,
the semigroup property and symmetry give
\[
\sup_{x\in\mathcal O}\int_{\mathcal O}|G_r(x,z)|^2\ud z
=\sup_{x\in\mathcal O}H_{2r}(x,x)
\le Cr^{-d/2}.
\]
For the wave equation, Parseval's identity and
$|\sin(r\sqrt{\lambda_k})|\le1$ yield
\begin{align*}
	\sup_{x\in\mathcal O}\int_{\mathcal O}|G_r(x,z)|^2\ud z
	\le
	\sup_{x\in\mathcal O}
	\sum_{k\ge1}\frac{|\varphi_k(x)|^2}{\lambda_k}	=
	\sup_{x\in\mathcal O}\int_0^\infty H_s(x,x)\ud s
	<\infty
	\qquad\text{if }d<2.
\end{align*}
Here \eqref{eq:diagonal-heat-kernel-bound} controls the integral near
$s=0$.  For $s\ge1$, the spectral expansion and $\lambda_k\ge\lambda_1>0$
give
\begin{align*}
	H_s(x,x)=
	\sum_{k\ge1}e^{-\lambda_k(s-1)}
	e^{-\lambda_k}|\varphi_k(x)|^2
	\le e^{-\lambda_1(s-1)}H_1(x,x),
\end{align*}
and therefore $
\sup_{x\in\mathcal O}\int_1^\infty H_s(x,x)\ud s
\le
\frac{\sup_{x\in\mathcal O}H_1(x,x)}{\lambda_1}
<\infty.$
Thus the spectral gap $\lambda_1>0$ controls only the long-time part,
whereas the restriction $d<2$ comes entirely from the behavior near
$s=0$.

Finally, since $e^{-2y^2}\le e^{1/8}e^{-y}$ for $y\ge0$, the fourth-order 
kernel satisfies (see also \cite[Lemma~1.2]{CC01})
\begin{align*}
	\sup_{x\in\mathcal O}\int_{\mathcal O}|G_r(x,z)|^2\ud z
	=
	\sup_{x\in\mathcal O}
	\sum_{k\ge1}e^{-2r\lambda_k^2}|\varphi_k(x)|^2
	\le
	e^{1/8}\sup_{x\in\mathcal O}H_{\sqrt r}(x,x)
	\le Cr^{-d/4}.
\end{align*}
This proves the bounds and dimensional conditions
recorded in Table~\ref{tab:prototype-kernel-dimensions}.

We next verify \Kfive{}, which concerns only the continuous Green
function, for the same three prototype operators with homogeneous
Dirichlet boundary conditions.

\begin{proposition}
\label{prop:K5-prototype-operators}
Let $\mathcal O\subset\mathbb R^d$ be a bounded connected domain with
sufficiently regular boundary, and let $x_0\in\mathcal O$.
For each of the operators
$$
\partial_t-\Delta,
\qquad
\partial_{tt}-\Delta,
\qquad
\partial_t+\Delta^2,
$$
the corresponding Green
function satisfies \Kfive{}.
\end{proposition}

\begin{proof}
Retain the Dirichlet eigenpairs
$\{(\lambda_k,\varphi_k)\}_{k\ge1}$ introduced above.
Since $\mathcal O$ is connected, the
principal-eigenvalue theorem
\cite[Theorem~8.38]{GilbargTrudinger2001} shows that $\lambda_1$ is simple
and that $\varphi_1$ can be chosen strictly positive in $\mathcal O$.
Consequently,
$$
        \lambda_1<\lambda_2
        \qquad\text{and}\qquad
        \varphi_1(x_0)>0.
$$
Recall from Examples \ref{Example1}--\ref{Example3} that the Green integral operators associated with the three equations are, respectively,
$e^{-rA},
 ~A^{-1/2}\sin(rA^{1/2}),
~e^{-rA^2}.
$
Hence their Green functions have the spectral form
$$
        G_r(x,z)=\sum_{k\ge1}m_k(r)\varphi_k(x)\varphi_k(z),
        \qquad r>0,
$$
with the series understood in the corresponding kernel sense (in
particular, in $L^2(\mathcal O)$ with respect to $z$ whenever that kernel
belongs to $L^2(\mathcal O)$), where
$$
m_k(r)=
\begin{cases}
e^{-\lambda_k r},
& \mathcal L=\partial_t-\Delta,\\[2mm]
\dfrac{\sin(\sqrt{\lambda_k}r)}{\sqrt{\lambda_k}},
& \mathcal L=\partial_{tt}-\Delta,\\[3mm]
e^{-\lambda_k^2r},
& \mathcal L=\partial_t+\Delta^2.
\end{cases}
$$
By orthonormality and Parseval's identity,
\begin{align*}
\int_{\mathcal O}|G_r(x_0,z)|^2\ud z
=\sum_{k\ge1}|m_k(r)|^2|\varphi_k(x_0)|^2
\ge
|m_1(r)|^2|\varphi_1(x_0)|^2.
\end{align*}
The first modal coefficient $m_1$ is nonzero almost everywhere on
$(0,T)$: it is strictly positive for the two parabolic operators, while
$\sin(\sqrt{\lambda_1}r)$ vanishes only on a discrete subset of
$(0,T)$ for the wave operator.  Therefore,
$$
\begin{aligned}
\int_0^T\int_{\mathcal O}|G_{T-s}(x_0,z)|^2\ud z\ud s
&=
\int_0^T\int_{\mathcal O}|G_r(x_0,z)|^2\ud z\ud r\\
&\ge
|\varphi_1(x_0)|^2\int_0^T|m_1(r)|^2\ud r
>0.
\end{aligned}
$$
This is precisely \Kfive{}.
\end{proof}

\subsection{Approximation by space-time grid controls}
\label{subsec:grid-controls-LDP}
Let
$\{C_{n,i}\}_{i=1}^{N_n}$ be the
spatial cells of a mesh of $\mathcal O$ whose maximal diameter tends to zero,
and let $t_j=j\tau_n$, $0\le j\le m_n$.  We use the grid-control space
$$
\mathcal X_n
:=
\left\{
\sum_{j=0}^{m_n-1}\sum_{i=1}^{N_n}c_{j,i}
\mathbf 1_{[t_j,t_{j+1})\times C_{n,i}}:
 c_{j,i}\in\mathbb R
\right\}\subset L^2(\mathcal O_T).
$$
For $h\in L^2(\mathcal O_T)$, let $Q_nh$ be its average on every space-time cell:
$$
(Q_nh)(s,z)
:=
\frac{1}{\tau_n|C_{n,i}|}
\int_{t_j}^{t_{j+1}}\int_{C_{n,i}}h(r,\xi)\ud\xi\ud r,
\qquad
(s,z)\in[t_j,t_{j+1})\times C_{n,i}.
$$

\begin{lemma}[Verification of Assumption~\ref{ass:abstract-control-approx}]
\label{lem:grid-control-approximation}
For every $h\in L^2(\mathcal O_T)$,
$$
        Q_nh\to h
        \qquad\text{in }L^2(\mathcal O_T)
        \quad\text{as }n\to\infty.
$$
\end{lemma}

\begin{proof}
{Set
$D_{j,i}:=[t_j,t_{j+1})\times C_{n,i}$.  For every
$v=\sum_{j,i}c_{j,i}\mathbf 1_{D_{j,i}}\in\mathcal X_n$, the definition
of the cell average gives
$$
\langle h-Q_nh,v\rangle_{L^2(\mathcal O_T)}
=
\sum_{j,i}c_{j,i}\int_{D_{j,i}}(h-Q_nh)\ud s\ud z
=0.
$$
Since $Q_nh\in\mathcal X_n$, this proves that $Q_n$ is the
$L^2(\mathcal O_T)$-orthogonal projection onto $\mathcal X_n$; in particular, it
is a contraction.}  If
$\phi\in C(\OT)$, its uniform continuity and the vanishing maximal
space-time cell diameter imply
$\|Q_n\phi-\phi\|_{L^2(\mathcal O_T)}\to0$.  The conclusion for an arbitrary
$h\in L^2(\mathcal O_T)$ now follows from the density of $C(\OT)$ in $L^2(\mathcal O_T)$
and the contraction property of $Q_n$.
\end{proof}

\subsection{Concrete fully discrete schemes}
\label{subsec:concrete-fully-discrete-schemes}

We now specify the equations, numerical schemes, and continuous and discrete
Green functions for the two applications.  In accordance with
Table~\ref{tab:prototype-kernel-dimensions}, we treat the wave equation in one
spatial dimension and the fourth-order equation in dimensions
$d\in\{1,2,3\}$.

\subsubsection{One-dimensional wave equation: FDM--EEM}
\begin{example}
\label{ex:wave-fdm-eem}
Let $\mathcal O=(0,1)$.  We consider the stochastic wave equation
\begin{equation}\label{sec4.3SWE}
\left\{
\begin{aligned}
&\partial_t^2u_\varepsilon(t,x)-\partial_x^2u_\varepsilon(t,x)
=b(u_\varepsilon(t,x))
+\sqrt\varepsilon\,\sigma(u_\varepsilon(t,x))\dot W(t,x),\\
&u_\varepsilon(t,0)=u_\varepsilon(t,1)=0,\\
&u_\varepsilon(0,x)=u_0(x),\quad
\partial_tu_\varepsilon(0,x)=v_0(x).
\end{aligned}
\right.
\end{equation}
Let $e_\ell(x)=\sqrt2\sin(\ell\pi x)$, $\ell\ge1$.  The continuous
Green function is given by
$$
G_r(x,z)=\sum_{\ell=1}^{\infty}
\frac{\sin(\ell\pi r)}{\ell\pi}e_\ell(x)e_\ell(z),
\qquad r\ge0,
$$
and the homogeneous solution is
$$
g(t,x)=\sum_{\ell=1}^{\infty}
\left[\cos(\ell\pi t)\langle u_0,e_\ell\rangle_{L^2(0,1)}
+\frac{\sin(\ell\pi t)}{\ell\pi}
\langle v_0,e_\ell\rangle_{L^2(0,1)}\right]e_\ell(x).
$$
With these choices of $G$ and $g$, the equation and its controlled
skeleton have the forms \eqref{eq:abstract-spde} and
\eqref{eq:continuous-skeleton}, respectively.

We next describe the FDM--EEM approximation of \eqref{sec4.3SWE}.  Put $x_k=k/n$.
Let $A_n\in\mathbb R^{(n-1)\times(n-1)}$ be the finite-difference
Dirichlet Laplacian, i.e., $(A_n)_{i,i}=-2n^2$, $(A_n)_{i,j}=n^2$ for $|i-j|=1$ and $(A_n)_{i,j}=0$ for $|i-j|>1$.
%$A_n=n^2
%\begin{pmatrix}
%-2&1&&0\\
%1&-2&\ddots&\\
%&\ddots&\ddots&1\\
%0&&1&-2
%\end{pmatrix}.$
Set
$$
\mathcal A_n=\begin{pmatrix}0&I_{n-1}\\ A_n&0\end{pmatrix},
\qquad \Omega_n=\sqrt{-A_n},
$$
and write 
$$
\begin{aligned}
C_n(t)=\cos(t\Omega_n),\qquad
S_n(t)=\Omega_n^{-1}\sin(t\Omega_n),\\
E_n(t)=e^{t\mathcal A_n}
=\begin{pmatrix}C_n(t)&S_n(t)\\A_nS_n(t)&C_n(t)\end{pmatrix}.
\end{aligned}
$$
For $U\in\mathbb R^{n-1}$, define
$$
B_n(U)=\bigl(b([U]_1),\ldots,b([U]_{n-1})\bigr)^\top,
\qquad
D_n(U)=\operatorname{diag}\bigl(\sigma([U]_1),\ldots,
\sigma([U]_{n-1})\bigr),
$$
and
$$
F_n(U)=\begin{pmatrix}0\\B_n(U)\end{pmatrix},
\qquad
\Sigma_n(U)=\begin{pmatrix}0\\D_n(U)\end{pmatrix}.
$$
Let $\tau=\tau_n$, and let
$\boldsymbol\beta=(\beta_1,\ldots,\beta_{n-1})^\top$ be the standard
$(n-1)$-dimensional Brownian motion obtained from the normalized spatial
increments of $W$, namely,
$$
\beta_k(t)=\sqrt n\bigl(W(t,(k+1)/n)-W(t,k/n)\bigr),
\qquad k=1,\ldots,n-1.
$$
For the $j$th time interval, put
$$
\Delta_j\beta_k=\beta_k(t_{j+1})-\beta_k(t_j),
\qquad
\Delta_j\boldsymbol\beta
=\boldsymbol\beta(t_{j+1})-\boldsymbol\beta(t_j)
=(\Delta_j\beta_1,\ldots,\Delta_j\beta_{n-1})^\top.
$$
Thus $j$ indexes the time increment and $k$ indexes its spatial component.
Set
$
X_{\varepsilon,j}^{n,m_n}
=\bigl(U_{\varepsilon,j}^{n,m_n},V_{\varepsilon,j}^{n,m_n}\bigr)^\top.
$
The fully discrete scheme of \eqref{sec4.3SWE} is
\begin{equation}\label{eq:wave-fdm-eem-random-recursion}
\begin{aligned}
X_{\varepsilon,j+1}^{n,m_n}
={}&E_n(\tau)X_{\varepsilon,j}^{n,m_n}
+\tau E_n(\tau)F_n(U_{\varepsilon,j}^{n,m_n})\\
&+\sqrt{\varepsilon n}\,E_n(\tau)
\Sigma_n(U_{\varepsilon,j}^{n,m_n})\Delta_j\boldsymbol\beta,
\qquad j=0,\ldots,m_n-1.
\end{aligned}
\end{equation}
The deterministic initial vector
$X_{\varepsilon,0}^{n,m_n}=X_0^n=(U_0^n,V_0^n)^\top$ is obtained from
$u_0$ and $v_0$ on the spatial grid.
The discrete Green function used for the continuous interpolation is
$G_n(t,\lfloor s\rfloor_n;x,z)
=\mathbf{G}^n(t-\lfloor s\rfloor_n,x,z)$, where
(see p.~404 in \cite{Cohen16})
\begin{equation}\label{eq:wave-discrete-green}
{\mathbf{G}^n(t,x,z)=\sum_{\ell=1}^{n-1}
\frac{\sin(\ell\pi t\sqrt{c_\ell^n})}
{\ell\pi\sqrt{c_\ell^n}}e_\ell^n(x)e_\ell(\kappa_n(z))}
\end{equation}
with $c_\ell^n=\sin^2(\ell\pi/2n)/(\ell\pi/2n)^2$,
{$\kappa_n(z)=\lfloor nz\rfloor/n$}, $e_\ell(x)={\sqrt2\sin(\ell\pi x)}$,
and
$$
e_\ell^n(x)=e_\ell(\kappa_n(x))+n(x-\kappa_n(x))
\bigl(e_\ell(\kappa_n(x)+1/n)-e_\ell(\kappa_n(x))\bigr).
$$
Let $g_n(t,\cdot)$ be the piecewise-linear spatial interpolation of
$C_n(t)U_0^n+S_n(t)V_0^n$.  We define the continuous space-time
interpolation of \eqref{eq:wave-fdm-eem-random-recursion} by
\begin{equation}\label{eq:wave-mild-extension}
\begin{aligned}
&u_\varepsilon^{n,m_n}(t,x)
={}g_n(t,x)+\int_0^t\int_0^1
\mathbf G^n(t-\lfloor s\rfloor_n,x,z)
b\!\left(u_\varepsilon^{n,m_n}(\lfloor s\rfloor_n,\kappa_n(z))\right)
\ud z\ud s\\
&+\sqrt\varepsilon\int_0^t\int_0^1
\mathbf G^n(t-\lfloor s\rfloor_n,x,z)
\sigma\!\left(u_\varepsilon^{n,m_n}(\lfloor s\rfloor_n,\kappa_n(z))\right)
W(\ud s\ud z),
\quad (t,x)\in\OT.
\end{aligned}
\end{equation}
Then it holds that $u_\varepsilon^{n,m_n}(t_j,x_k)=[U_{\varepsilon,j}^{n,m_n}]_k$, as verified  in the proof of
Proposition~\ref{prop:grid-control-LDP}. Hereafter, for a vector $a$, $[a]_k$ denotes its $k$th component.
For $x_0\in[x_k,x_{k+1})$, define the terminal interpolation vector by
\begin{align}\label{eq:elln}
\ell_{n,x_0}^{\top}U
=(k+1-nx_0)[U]_k+(nx_0-k)[U]_{k+1},
\qquad [U]_0=[U]_n=0.
\end{align}
\end{example}

\begin{corollary}[Diagonal convergence for the FDM--EEM scheme]
\label{cor:wave-diagonal-convergence}
Assume that $u_0\in H^\alpha([0,1])$ and
$v_0\in H^{\alpha-1}([0,1])$ for some $\alpha>3/2$, that
$b,\sigma\in C_b^1(\mathbb R)$, and that $\sigma$ is nowhere zero.  Then
$\{u_\varepsilon^{n,m_n}(T,x_0)\}_{\varepsilon>0}$ satisfies, for every fixed $n$, an LDP with good rate function
$$
I_{\mathrm w}^{n,m_n}(y)
=\inf_{\{h_n\in\mathcal X_n:
\Un(h_n)(T,x_0)=y\}}
\frac12\|h_n\|_{L^2(\mathcal O_T)}^2,
$$
and
$$
I_{\mathrm w}^{n,m_n}(y)\rightarrow I_{\mathrm w}(y)
\qquad\text{as }n\to\infty,
\quad y\in\mathbb R.
$$
In particular, for every $y\in\mathbb R$ with
$I_{\mathrm w}(y)<+\infty$, the conclusions of
Corollary~\ref{cor:subsequential-minimizers} hold.
\end{corollary}

\begin{proof}
The kernel conditions \Kone--\Kfive{} are verified in
Appendix~\ref{app:wave-kernel-verification}, and
Proposition~\ref{prop:grid-control-LDP} proves the fixed-grid LDP and
identifies its rate function.  The conclusion follows from
Theorem~\ref{thm:closed-diagonal} and
Lemma~\ref{lem:grid-control-approximation}.
\end{proof}

\subsubsection{Multidimensional fourth-order parabolic equation:  FDM--implicit Euler method}

\begin{example}
\label{ex:fourth-order-fdm}

Let $d\in\{1,2,3\}$, $\mathcal O=(0,\pi)^d$, and $A=-\Delta$.  Consider
the fourth-order parabolic equation
\begin{align}
\partial_tu_\varepsilon+A^2u_\varepsilon
=b(u_\varepsilon)+\sqrt{\varepsilon}\,
\sigma(u_\varepsilon)\dot W,
\qquad
u_\varepsilon(0)=0, \label{sec4.5parabolic}
\end{align}
with the boundary conditions  $u_\varepsilon=\Delta u_\varepsilon=0$ on
$\partial\mathcal O$.   Since the
initial value is zero, $g\equiv0$. Then the exact mild solution has the form
stated in Section~\ref{sec:spde-notation}, with
$$
G_r(x,z)
=\sum_{\mathbf k\in\mathbb N^d}e^{-r|\mathbf k|^4}
\varphi_{\mathbf k}(x)\varphi_{\mathbf k}(z),
\qquad r>0,
$$
where $\varphi_{\mathbf k}(x)=(2/\pi)^{d/2}
\prod_{\ell=1}^d\sin(k_\ell x_\ell)$.

We use the uniform mesh $\Delta x_n=\pi/n$ and write
$\tau_n=T/m_n$, where $m_n\to\infty$.  Put
$$
\mathbb I_n=\{0,\ldots,n-1\}^d,
\qquad \mathbb J_n=\{1,\ldots,n-1\}^d,
$$
and, for $\mathbf k=(k_1,\ldots,k_d)\in\mathbb I_n$, define
$$
C_{n,\mathbf k}
=\prod_{\ell=1}^d[k_\ell \Delta x_n,(k_\ell+1)\Delta x_n),
\qquad
\kappa_n(z)=\Delta x_n\mathbf k\quad\text{for }z\in C_{n,\mathbf k}.
$$
Let $L_n$ be the positive one-dimensional finite-difference Dirichlet Laplacian, i.e.,
$$
[L_nv]_i=\frac{2v_i-v_{i-1}-v_{i+1}}{(\Delta x_n)^2},
\qquad i=1,\ldots,n-1,\qquad v_0=v_n=0,
$$
and let
$$
L_{n,d}
=\sum_{\ell=1}^{d}
I_{n-1}^{\otimes(\ell-1)}\otimes L_n
\otimes I_{n-1}^{\otimes(d-\ell)}.
$$
Here \(I_{n-1}\) denotes the \((n-1)\times(n-1)\) identity matrix and \(\otimes\) denotes the Kronecker product.
For every $\mathbf k\in\mathbb J_n$, the nodal values of
$\varphi_{\mathbf k}$ form an eigenvector of $L_{n,d}$:
$$
L_{n,d}
\bigl(\varphi_{\mathbf k}(\Delta x_n\mathbf p)\bigr)_{\mathbf p\in\mathbb J_n}
=\mu_{\mathbf k,n}
\bigl(\varphi_{\mathbf k}(\Delta x_n\mathbf p)\bigr)_{\mathbf p\in\mathbb J_n},
$$
where the corresponding eigenvalue is
$
\mu_{\mathbf k,n}
=\frac{4}{(\Delta x_n)^2}\sum_{\ell=1}^d
\sin^2\!\left(\frac{k_\ell\pi}{2n}\right).$ 
For the finite-difference approximation of the bi-Laplacian, we also refer to \cite{HongJinSheng2024} for the one-dimensional case.

For $x\in C_{n,\mathbf p}$, $\mathbf p\in\mathbb I_n$, set
$\vartheta_\ell(x)=(x_\ell-p_\ell \Delta x_n)/\Delta x_n$.  We define the multilinear
grid interpolant $\varphi_{\mathbf k}^n$ by
$$
\varphi_{\mathbf k}^n(x)
=\sum_{\alpha\in\{0,1\}^d}
\left[
\prod_{\ell=1}^d
\bigl((1-\alpha_\ell)(1-\vartheta_\ell(x))
+\alpha_\ell\vartheta_\ell(x)\bigr)
\right]
\varphi_{\mathbf k}\bigl(\Delta x_n(\mathbf p+\alpha)\bigr),
$$
and extend it continuously to $\overline{\mathcal O}$.
For $U\in\mathbb R^{(n-1)^d}$, set
$$
[B_n(U)]_{\mathbf k}=b([U]_{\mathbf k}),
\qquad
D_n(U)=\operatorname{diag}
\bigl(\sigma([U]_{\mathbf k})\bigr)_{\mathbf k\in\mathbb J_n}.
$$
Writing $\mathcal W(A)=\int_AW(ds,dz)$ for the Gaussian random measure
induced by the Brownian sheet, define
$$
[\Delta_j\boldsymbol\beta]_{\mathbf k}
=(\Delta x_n)^{-d/2}\mathcal W\bigl([t_j,t_{j+1})\times C_{n,\mathbf k}\bigr),
\qquad \mathbf k\in\mathbb J_n.
$$
Denoting
$
R_n=(I+\tau_nL_{n,d}^2)^{-1}
$, the nodal FDM--implicit Euler scheme for \eqref{sec4.5parabolic} reads
\begin{equation}\label{eq:fourth-order-implicit-scheme}
U_{\varepsilon,j+1}^n
=R_n\left[
U_{\varepsilon,j}^n+\tau_nB_n(U_{\varepsilon,j}^n)
+\sqrt{\varepsilon}\,(\Delta x_n)^{-d/2}
D_n(U_{\varepsilon,j}^n)\Delta_j\boldsymbol\beta
\right],
\end{equation}
with $U_{\varepsilon,0}^n=0\in\mathbb R^{(n-1)^d}$.
For $\mu>0$, set
$\chi_n(\mu)=(1+\tau_n\mu^2)^{-1}$ and define
\begin{equation}\label{eq:fourth-q-definition}
q_n(r,\mu)=
\begin{cases}
\chi_n(\mu),&0\le r\le\tau_n,\\[1mm]
(1-\theta)\chi_n(\mu)^\ell
+\theta\chi_n(\mu)^{\ell+1},
&\begin{gathered}
r\in[\ell\tau_n,(\ell+1)\tau_n],\ell\ge1,
\theta=r/\tau_n-\ell.
\end{gathered}
\end{cases}
\end{equation}
Introduce the fully discrete Green function 
\begin{equation}\label{eq:fourth-discrete-green}
G_n(t,r;x,z)
=\sum_{\mathbf k\in\mathbb J_n}
q_n(t-r,\mu_{\mathbf k,n})
\varphi_{\mathbf k}^n(x)
\varphi_{\mathbf k}(\kappa_n(z)),
\qquad 0\le r<t\le T.
\end{equation}
We then define a continuous space-time interpolation of the nodal scheme  \eqref{eq:fourth-order-implicit-scheme} by
\begin{equation}\label{eq:fourth-mild-extension}
\begin{aligned}
u_{\varepsilon,n}(t,x)
={}&\int_0^t\int_{\mathcal O}
G_n(t,\lfloor s\rfloor_n;x,z)
b\!\left(u_{\varepsilon,n}(\lfloor s\rfloor_n,\kappa_n(z))\right)
\ud z\ud s\\
&+\sqrt{\varepsilon}\int_0^t\int_{\mathcal O}
G_n(t,\lfloor s\rfloor_n;x,z)
\sigma\!\left(u_{\varepsilon,n}(\lfloor s\rfloor_n,\kappa_n(z))\right)
\mathcal W(ds,dz).
\end{aligned}
\end{equation}
Then \(u_{\varepsilon,n}(t_j,\Delta x_n \mathbf k)
=
[U_{\varepsilon,j}^n]_{\mathbf k},
~
j=0,\ldots,m_n,~
\mathbf k\in\mathbb J_n\), as verified  in the proof of
Proposition~\ref{prop:fourth-fixed-grid-LDP}.
For $x_0\in\mathcal O$, define
$\mathbf i=(i_1,\ldots,i_d)$ with $i_{\ell}=\lfloor x_{0,\ell}/\Delta x_n\rfloor$ and
$\theta_\ell=(x_{0,\ell}-i_\ell \Delta x_n)/\Delta x_n$ for $\ell=1,\ldots,d$.  The terminal value is obtained
from the multilinear interpolation operator
$$
\ell_{n,x_0}(U)
=
\sum_{\alpha\in\{0,1\}^d}
\left[
\prod_{\ell=1}^d
\bigl((1-\alpha_\ell)(1-\theta_\ell)+\alpha_\ell\theta_\ell\bigr)
\right][U]_{\mathbf i+\alpha},
\qquad U\in\mathbb R^{(n-1)^d},
$$
where $[U]_{\mathbf k}=0$ if $\mathbf k\notin\mathbb J_n$.
\end{example}

\begin{corollary}[Diagonal convergence for the fourth-order FDM--implicit Euler scheme]
\label{cor:fourth-diagonal-convergence}
Let $d\in\{1,2,3\}$.  Assume that
$b,\sigma\in C_b^1(\mathbb R)$ and that $\sigma$ is nowhere zero.  Then  $\{u_{\varepsilon,n}(T,x_0)\}_{\varepsilon>0}$ satisfies, for every fixed $n$, an LDP with good rate function
$$
I_d^{n,m_n}(y)
=\inf_{\{h_n\in\mathcal X_n:
\Un(h_n)(T,x_0)=y\}}
\frac12\|h_n\|_{L^2(\mathcal O_T)}^2,
$$
and, for every $x_0\in\mathcal O$ and $y\in\mathbb R$,
$$
I_d^{n,m_n}(y)\to I_d(y)
\qquad\text{as }n\to\infty.
$$
In particular, for every $y\in\mathbb R$ with $I_d(y)<+\infty$, the
conclusions of Corollary~\ref{cor:subsequential-minimizers} hold.
\end{corollary}

\begin{proof}
The kernel conditions \Kone--\Kfive{} are verified in
Appendix~\ref{app:fourth-kernel-verification}, while
Proposition~\ref{prop:fourth-fixed-grid-LDP} proves the fixed-grid LDP and
identifies its rate function.  Applying
Theorem~\ref{thm:closed-diagonal} and
Lemma~\ref{lem:grid-control-approximation} completes the proof.
\end{proof}

The two corollaries establish diagonal convergence without prescribing a
relative rate of spatial and temporal refinement.  

\section{Numerical experiments for the stochastic wave equation}
\label{sec:numerical-experiments}
In this section, we present two numerical experiments for the wave FDM--EEM
scheme in Example~\ref{ex:wave-fdm-eem} to illustrate the diagonal convergence
established in Corollary~\ref{cor:wave-diagonal-convergence}.  
%In the linear
%additive-noise case, both the continuous and fully discrete rate functions are
%available explicitly, so the pointwise discrepancy
%$I_{\mathrm w}^{n,n}(y)-I_{\mathrm w}(y)$ can be evaluated without using a
%reference-grid approximation.  In the nonlinear case with state-dependent
%diffusion, the continuous rate function is not available in closed form.  We
%therefore compute the discrete rate functions by solving equality-constrained
%minimum-energy problems and assess their convergence relative to the finest
%grid used in the experiment.

\subsection{Experimental setting and discrete control problem}
\label{subsec:numerical-setting}
 We consider the stochastic wave equation on  $[0,T]\times(0,1)$:
\begin{equation} \label{SWE}
\left\{
\begin{aligned}
&\partial_t^2u_\varepsilon(t,x)-\partial_x^2u_\varepsilon(t,x)
=b(u_\varepsilon(t,x))
+\sqrt{\varepsilon}\,\sigma(u_\varepsilon(t,x))\dot W(t,x),\\
&u_\varepsilon(t,0)=u_\varepsilon(t,1)=0,\\
&u_\varepsilon(0,x)=0,\quad
\partial_tu_\varepsilon(0,x)=0.
\end{aligned}
\right.
\end{equation}
The final time and observation point are  set to
\[
T=1,\qquad x_0=\frac12,
\]
and we use the diagonal choice $m_n=n$, so that
$\tau=T/m_n=1/n$.  All spatial resolution parameters \(n\) used below are even; hence $x_0$
is a grid point and the observation vector $\ell_{n,x_0}$, defined by \eqref{eq:elln}, selects the
$n/2$-th interior component without interpolation.

At the grid level, let
$H=(H_0,\ldots,H_{m_n-1})\in(\R^{n-1})^{m_n}$.  Define 
$(U_j^{\tau H},V_j^{\tau H})$, $j=0,1,\ldots, m_n$, recursively by
\begin{equation}
\label{eq:num-fdm-eem-recursion}
\begin{aligned}
U_{j+1}^{\tau H}={}&C_n(\tau)U_j^{\tau H}
+S_n(\tau)V_j^{\tau H}
+\tau S_n(\tau)
\bigl(B_n(U_j^{\tau H})+D_n(U_j^{\tau H})H_j\bigr),\\
V_{j+1}^{\tau H}={}&A_nS_n(\tau)U_j^{\tau H}
+C_n(\tau)V_j^{\tau H}
+\tau C_n(\tau)
\bigl(B_n(U_j^{\tau H})+D_n(U_j^{\tau H})H_j\bigr),
\end{aligned}
\end{equation}
with $U_0^{\tau H}=V_0^{\tau H}=0$.  Here $A_n$,
$B_n$, $D_n$, $C_n$, and $S_n$ are defined in
Example~\ref{ex:wave-fdm-eem}.   For a prescribed terminal value
$y$, the rate function $I^{n,m_n}_{\mathrm w}$ in Corollary \ref{cor:wave-diagonal-convergence} has the following  equivalent expression
\begin{equation}
\label{eq:num-control-optimization}
I_{\mathrm w}^{n,m_n}(y)
=
\min_{H\in(\R^{n-1})^{m_n}}
\left\{
\frac{\tau}{2n}\sum_{j=0}^{m_n-1}\lVert H_j\rVert_2^2:
\ell_{n,x_0}^{\top}U_{m_n}^{\tau H}=y
\right\};
\end{equation}
see Appendix~\ref{app:wave-fixed-grid-LDP} for more details.   Here, $\|\cdot\|_2$ denotes the $2$-norm of a vector.

\subsection{An exactly verifiable linear additive-noise benchmark}
\label{subsec:numerical-linear}

We first set
\[
b(u)=0,\qquad \sigma(u)=1.
\]
For the continuous equation, the rate function is given by the
variational problem
\[
I_{\mathrm w}(y)
=\inf_{\{h\in L^2((0,1)\times(0,1)):
\int_0^1\int_0^1G_{1-s}(1/2,z)h(s,z)\ud z\ud s=y\}}
\frac12\int_0^1\int_0^1|h(s,z)|^2\ud z\ud s.
\]
Set
$V
:=\int_0^1\int_0^1|G_{1-s}(1/2,z)|^2\ud z\ud s
=\sum_{k\ge1}
\frac{e_k(1/2)^2}{(k\pi)^2}
\int_0^1\sin^2(k\pi r)\ud r
=\frac18$. 
By the equality case of the Cauchy--Schwarz inequality, for every
$y\in\mathbb R$ the unique minimizer of this variational problem is
\[
h_y^*(s,z)=\frac{y}{V}G_{1-s}(1/2,z),
\]
and the minimum value is
\[
I_{\mathrm w}(y)=\frac{y^2}{2V}=4y^2.
\]

For the fully discrete equation, set
$a_j:=S_n(T-t_j)^{\top}\ell_{n,x_0},
~ j=0,\ldots,m_n-1.$
It follows from \eqref{eq:wave-iterated-controlled-displacement},
$U_0^n=V_0^n=0$, $b=0$, and $\sigma=1$ that
$\ell_{n,x_0}^{\top}U_{m_n}^{\tau H}
=\tau\sum_{j=0}^{m_n-1}a_j^{\top}H_j$. Thus, 
the discrete rate function is given by
\[
I_{\mathrm w}^{n,m_n}(y)
=\min_{\{H\in(\mathbb R^{n-1})^{m_n}:
\tau\sum_{j=0}^{m_n-1}a_j^{\top}H_j=y\}}
\frac{\tau}{2n}\sum_{j=0}^{m_n-1}\lVert H_j\rVert_2^2.
\]
Define
$
V_{n,m_n}=n\tau\sum_{j=0}^{m_n-1}\lVert a_j\rVert_2^2.
$
Again by the equality case of the Cauchy--Schwarz inequality, the unique
minimizer of this finite-dimensional variational problem is
\[
H_j^*=\frac{ny}{V_{n,m_n}}a_j,
\qquad j=0,\ldots,m_n-1,
\]
and the minimum value is
\[
I_{\mathrm w}^{n,m_n}(y)=\frac{y^2}{2V_{n,m_n}}.
\]

Table~\ref{tab:num-linear-errors} reports the diagonal errors.  The maximum
rate error is evaluated over 401 equally spaced values in $[-1,1]$.
\begin{table}[H]
\centering
\caption{Linear additive-noise benchmark with $m_n=n$.}
\label{tab:num-linear-errors}
\begin{tabular}{rrrr}
\toprule
$n$ & $V_{n,m_n}$ & $|V_{n,m_n}-1/8|$
& $\max_{|y|\le1}|I_{\mathrm w}^{n,m_n}(y)-I_{\mathrm w}(y)|$\\
\midrule
8   & 0.126737845222 & $1.737845\times10^{-3}$ & $5.484850\times10^{-2}$\\
16  & 0.125713379194 & $7.133792\times10^{-4}$ & $2.269859\times10^{-2}$\\
32  & 0.125297425681 & $2.974257\times10^{-4}$ & $9.495029\times10^{-3}$\\
64  & 0.125107275791 & $1.072758\times10^{-4}$ & $3.429882\times10^{-3}$\\
128 & 0.125039410622 & $3.941062\times10^{-5}$ & $1.260742\times10^{-3}$\\
256 & 0.125014578947 & $1.457895\times10^{-5}$ & $4.664719\times10^{-4}$\\
\bottomrule
\end{tabular}
\end{table}

A log--log least-squares fit over $n=16,\ldots,256$ gives empirical
orders $1.414$ for the error in $V_{n,m_n}$ and $1.412$ for the
maximum rate error. Thus, the linear benchmark directly illustrates the diagonal convergence
$I_{\mathrm w}^{n,m_n}(y)\to I_{\mathrm w}(y)$ asserted by
Theorem~\ref{thm:closed-diagonal} for the FDM--EEM approximation in
Example~\ref{ex:wave-fdm-eem}.

\begin{figure}[!htbp]
\centering
\begin{minipage}[t]{0.49\textwidth}
\centering
\includegraphics[width=\linewidth]
{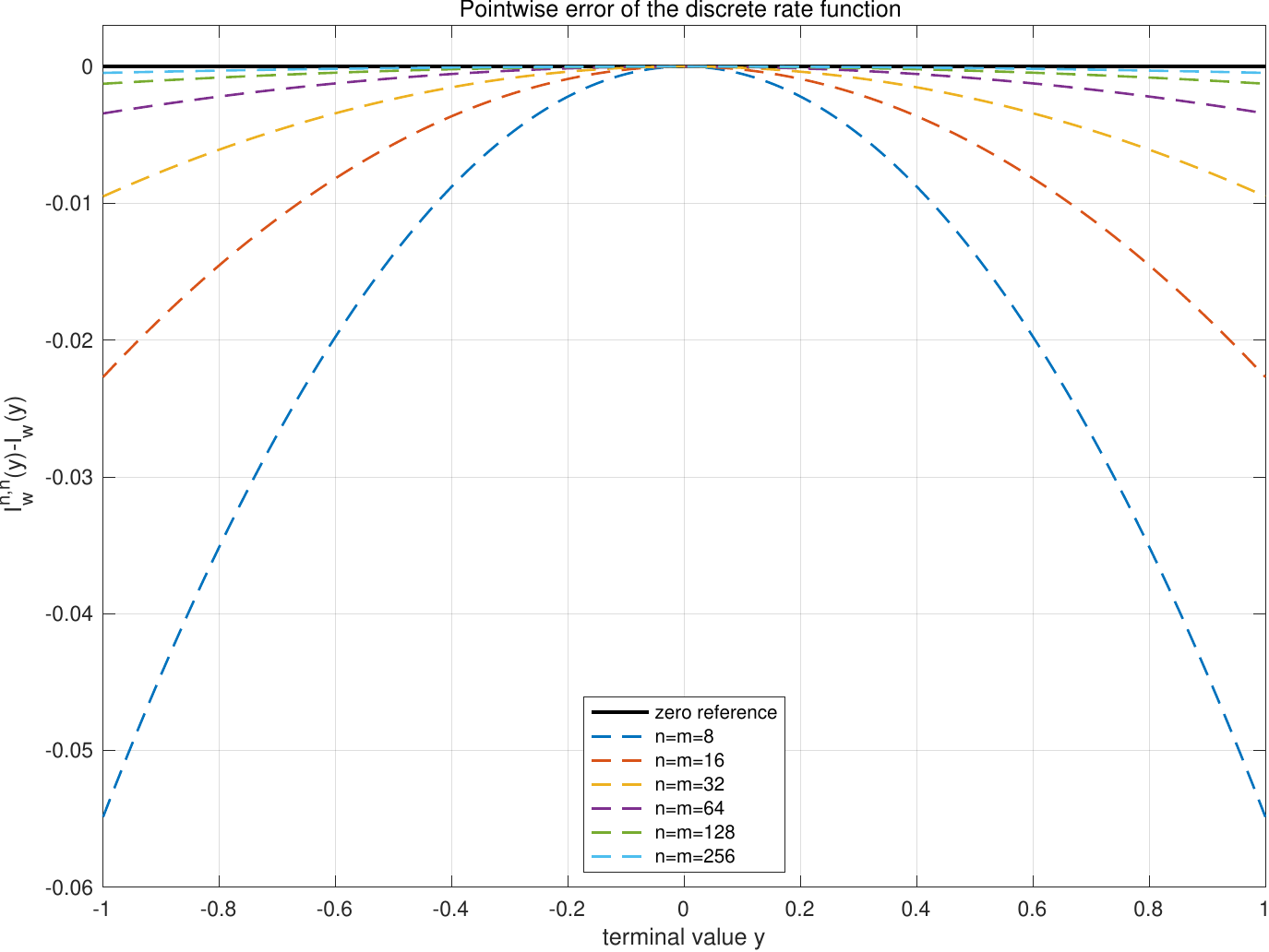}\\[-1mm]
{\small (a) Pointwise errors
$I_{\mathrm w}^{n,m_n}(y)-I_{\mathrm w}(y)$.}
\end{minipage}\hfill
\begin{minipage}[t]{0.49\textwidth}
\centering
\includegraphics[width=\linewidth]
{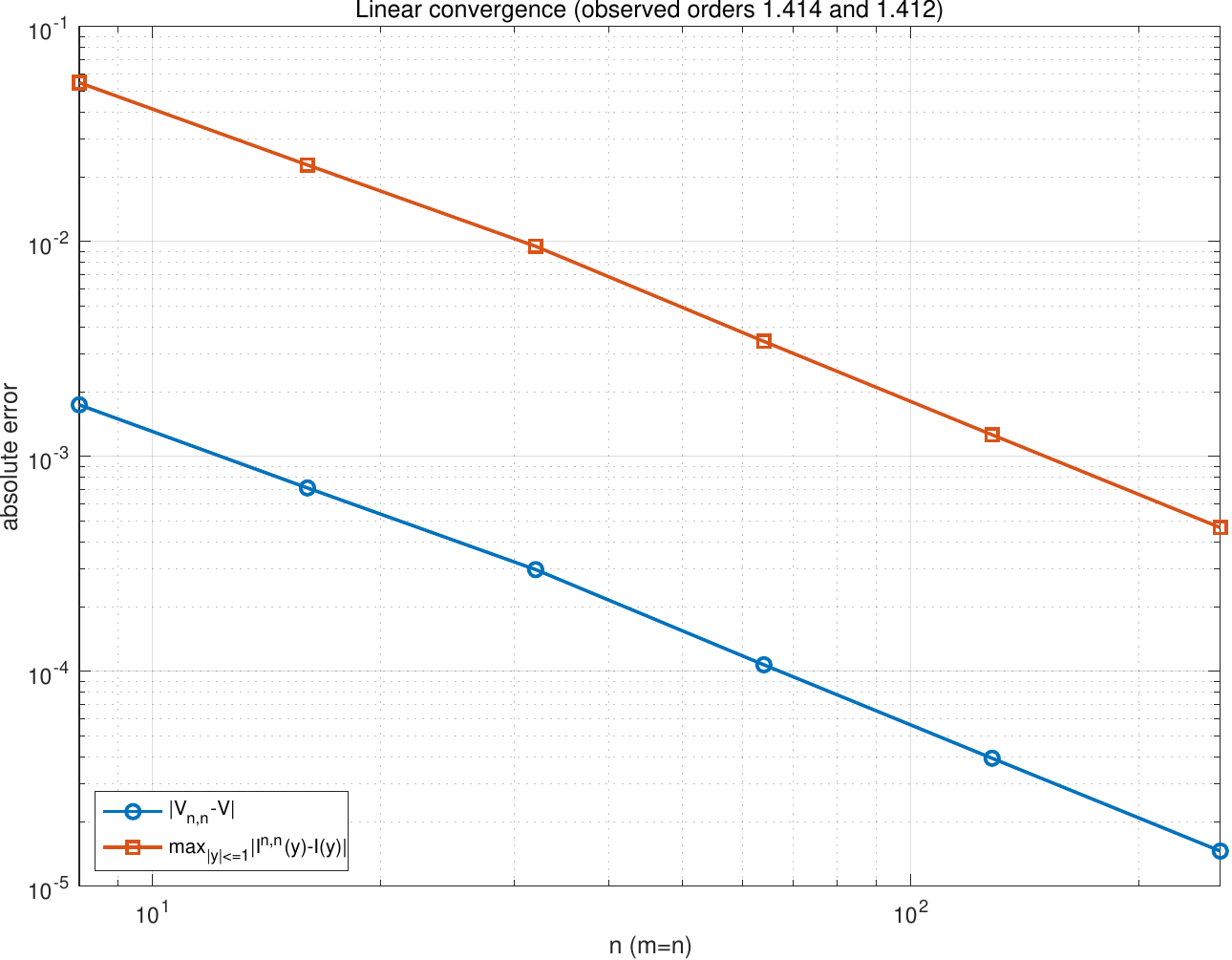}\\[-1mm]
{\small (b) Diagonal errors under simultaneous grid refinement.}
\end{minipage}
\caption{Linear additive-noise experiment.    All values
are obtained from the explicit  formulas of $I_{\mathrm w}=4y^2$ and $I_{\mathrm w}^{n,m_n}=\frac{y^2}{2V_{n,m_n}}$.}
\label{fig:num-linear-results}
\end{figure}

\subsection{Nonlinear drift and state-dependent diffusion}
\label{subsec:numerical-nonlinear}

For the nonlinear experiment, we use the same wave equation, initial data,
boundary conditions, terminal time, and observation point, but take
\[
b(u)=-\frac12\sin u,
\qquad
\sigma(u)=1+\frac14\tanh u.
\]
The functions $b$ and $\sigma$ belong to $C_b^1(\R)$, with
$b'(u)=-\frac12\cos u$ and
$\sigma'(u)=\frac14\mathrm{sech}^2u$, while
$3/4<\sigma(u)<5/4$.  Thus the coefficients satisfy the hypotheses used for
Corollary~\ref{cor:wave-diagonal-convergence}.

We evaluate \eqref{eq:num-control-optimization} for
$$
m_n=n\in\{8,12,16,24,32,48\},
\qquad
y\in\{-0.6,-0.5,\ldots,0.5,0.6\}.
$$
Each equality-constrained problem is solved by MATLAB's \texttt{fmincon},
which carries out sequential quadratic programming (SQP) iterations.  The
objective and equality-constraint gradients passed to \texttt{fmincon} are
computed analytically in our implementation.  Set
\[
J(H):=\frac{\tau}{2n}\sum_{j=0}^{m_n-1}\lVert H_j\rVert_2^2,
\qquad
\Phi(H):=\ell_{n,x_0}^{\top}U_{m_n}^{\tau H}.
\]
Thus \eqref{eq:num-control-optimization} is the problem of minimizing $J(H)$
subject to $\Phi(H)=y$.  The objective gradient is
$\nabla J(H)=(\tau/n)H$, while the gradient of the terminal map is evaluated
by the discrete adjoint recursion below, which we implement separately from
the optimization routine.  This recursion is obtained from
\eqref{eq:num-fdm-eem-recursion}.  We write $U_j=U_j^{\tau H}$ for brevity.
For fixed $H_j,\ldots,H_{m_n-1}$, let $\Phi_j(U,V)$ denote the
terminal observation obtained by starting \eqref{eq:num-fdm-eem-recursion}
from the state $(U,V)$ at time $t_j$, and define
\[
p_j^U:=\nabla_U\Phi_j(U_j,V_j),
\qquad
p_j^V:=\nabla_V\Phi_j(U_j,V_j).
\]
Since $\Phi_{m_n}(U,V)=\ell_{n,x_0}^{\top}U$, we have
$
p_{m_n}^U=\ell_{n,x_0},
~
p_{m_n}^V=0.$
We then set, backward for
$j=m_n-1,\ldots,0$,
\[
\begin{aligned}
w_j={}&S_n(\tau)^{\top}p_{j+1}^U
+C_n(\tau)^{\top}p_{j+1}^V,\\
p_j^U={}&C_n(\tau)^{\top}p_{j+1}^U
+\bigl(A_nS_n(\tau)\bigr)^{\top}p_{j+1}^V\\
&+\tau\bigl(b'(U_j)+\sigma'(U_j)\odot H_j\bigr)\odot w_j,\\
p_j^V={}&w_j,
\qquad
\nabla_{H_j}\bigl(\ell_{n,x_0}^{\top}U_{m_n}^{\tau H}\bigr)
=\tau\sigma(U_j)\odot w_j,
\end{aligned}
\]
where $\odot$ denotes componentwise multiplication.  

The optimizer for each
nonzero target was initialized using a previously computed solution at a
neighboring target value.  We used an SQP constraint
tolerance of $10^{-10}$, a first-order optimality tolerance of $2\times10^{-8}$,
a step tolerance of $10^{-12}$, and at most 300 iterations in the primary solve.
All computed solutions satisfied the prescribed terminal-constraint and
first-order optimality tolerances.  Since the
continuous nonlinear rate function is not available in closed form, the
$n=48$ curve is used only as the finest-grid numerical reference.  The
maximum discrepancies from that curve are listed in
Table~\ref{tab:num-nonlinear-errors}.
\begin{table}[H]
\centering
\caption{Nonlinear rate curves relative to the $m_n=n=48$ reference over the
13 target values in $[-0.6,0.6]$, where
$Err_n:=\max_y|I_{\mathrm w}^{n,m_n}(y)-I_{\mathrm w}^{48,48}(y)|$.}
\label{tab:num-nonlinear-errors}
\scriptsize
\setlength{\tabcolsep}{3pt}
\begin{tabular}{@{}c*{5}{c}@{}}
\toprule
$m_n=n$ & 8 & 12 & 16 & 24 & 32\\
\midrule
$Err_n$
& $2.908178\times10^{-2}$
& $1.677266\times10^{-2}$
& $9.995724\times10^{-3}$
& $5.179326\times10^{-3}$
& $2.468802\times10^{-3}$\\
\bottomrule
\end{tabular}
\end{table}

\begin{figure}[H]
\centering
\begin{minipage}[t]{0.49\textwidth}
\centering
\includegraphics[width=\linewidth]
{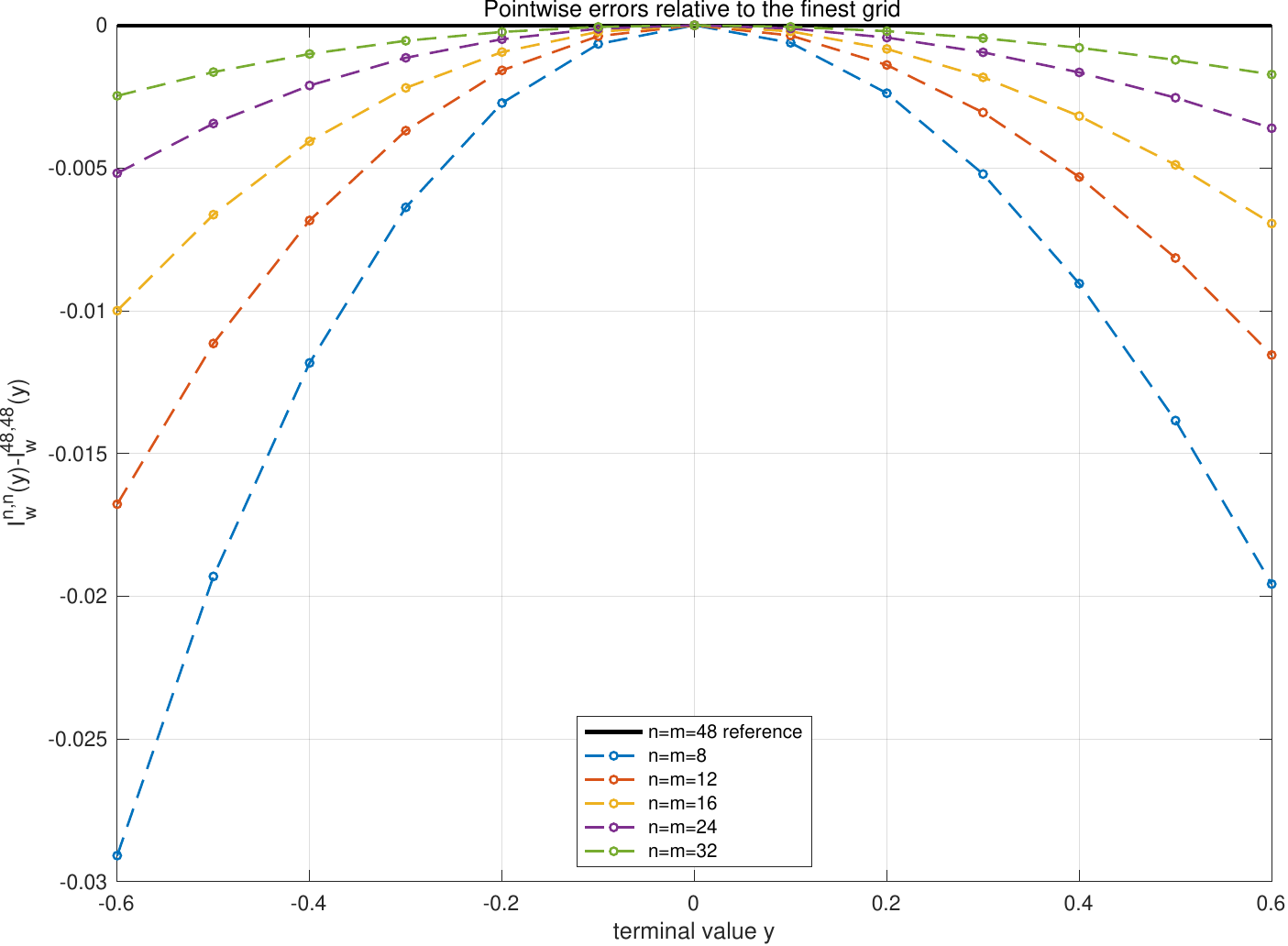}\\[-1mm]
{\small (a) Pointwise discrepancies
$I_{\mathrm w}^{n,m_n}(y)-I_{\mathrm w}^{48,48}(y)$.}
\end{minipage}\hfill
\begin{minipage}[t]{0.49\textwidth}
\centering
\includegraphics[width=\linewidth]
{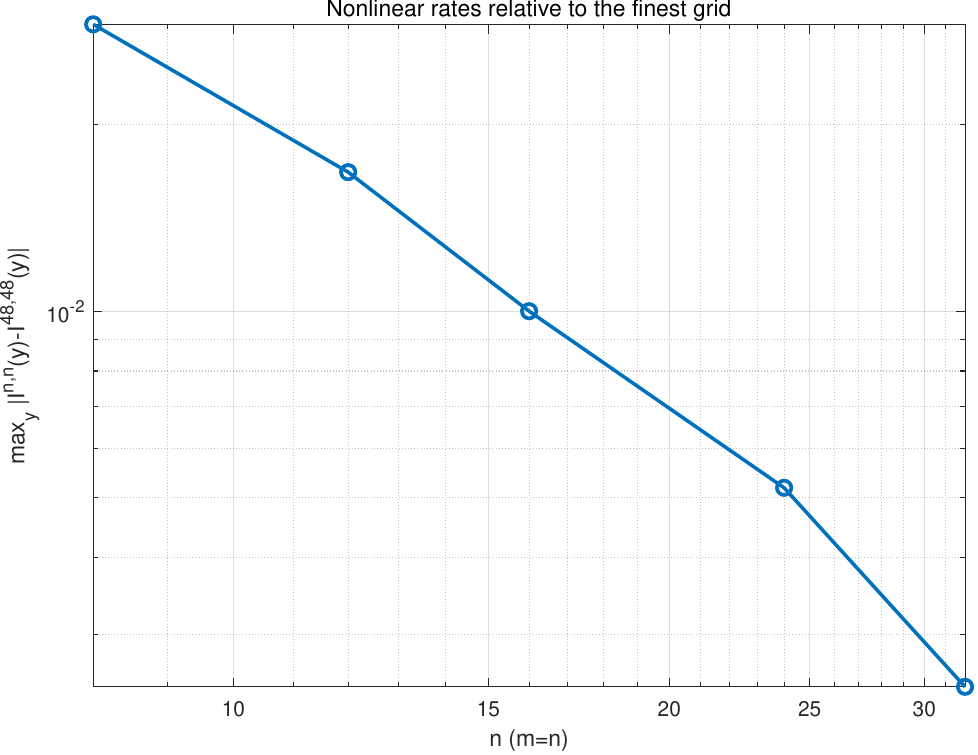}\\[-1mm]
{\small (b) Maximum discrepancy from the $m_n=n=48$ reference curve.}
\end{minipage}
\caption{Nonlinear wave experiment with
$b(u)=-\frac12\sin u$ and $\sigma(u)=1+\frac14\tanh u$.  }
\label{fig:num-nonlinear-results}
\end{figure}

Figure~\ref{fig:num-nonlinear-results} summarizes the pointwise and maximum
discrepancies relative to the $m_n=n=48$ reference curve.

\begin{figure}[!htbp]
\centering
\includegraphics[width=0.88\textwidth]
{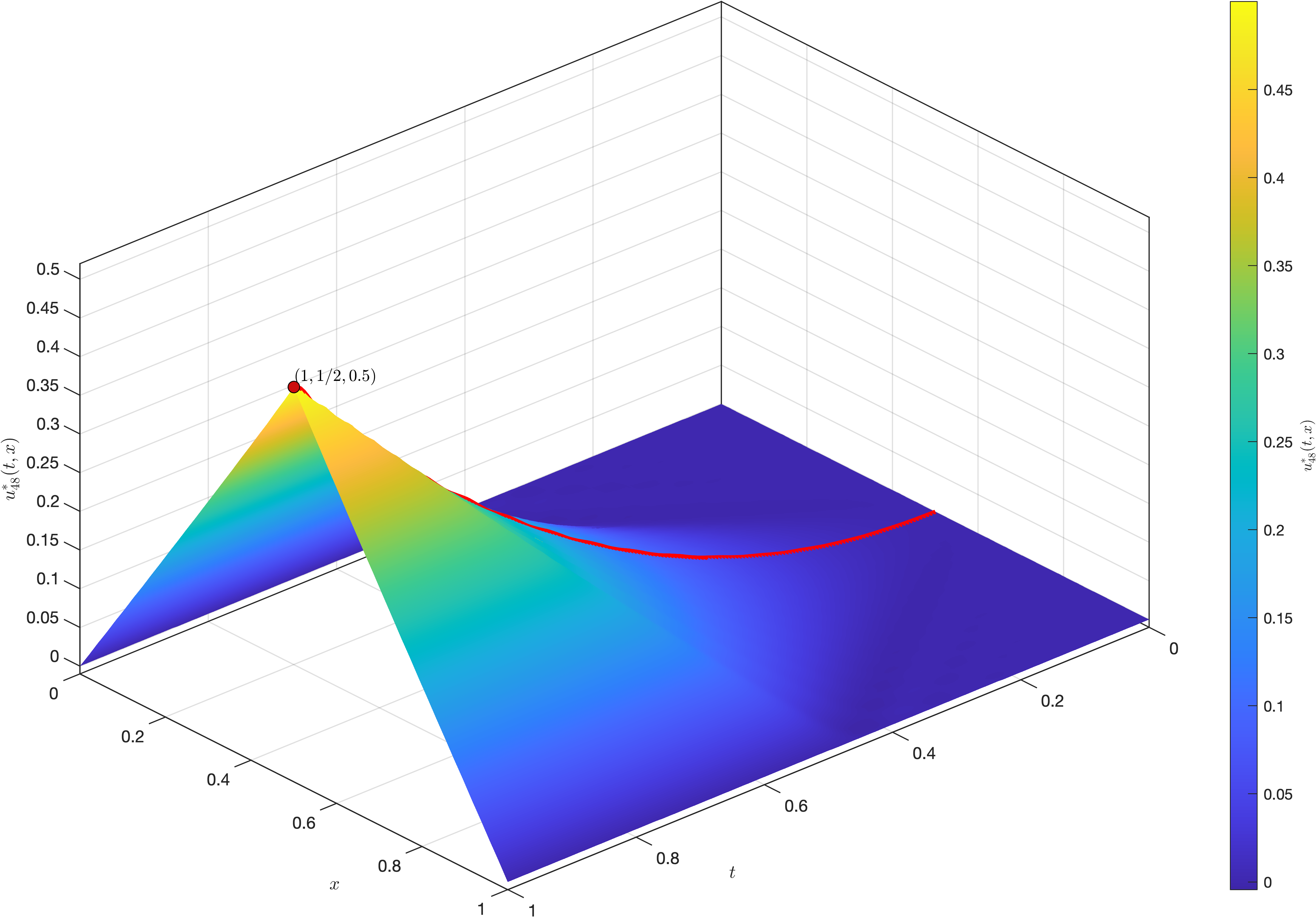}
\caption{The computed reference minimum action path associated with
$I_{\mathrm w}^{48,48}(0.5)$.  The surface is the continuous space-time
interpolation generated by the computed minimizing control.  The red curve
shows $t\mapsto u_{48}^*(t,1/2)$, and the marked terminal point satisfies
$u_{48}^*(1,1/2)=0.5$.}
\label{fig:num-reference-path}
\end{figure}

The finest-grid optimization also provides a computed minimizing control for
the target value $y=0.5$.  Let $H_{48}^*$ denote the control returned by the
SQP solve of \eqref{eq:num-control-optimization} with $m_n=n=48$ and $y=0.5$,
and let $h_{H_{48}^*}$ be the associated piecewise-constant space-time grid
control.  We define the corresponding continuous interpolation by
$u_{48}^*(t,x):=\Upsilon_{48}(h_{H_{48}^*})(t,x)$ for
$(t,x)\in[0,1]\times[0,1]$.
The computed action is $I_{\mathrm w}^{48,48}(0.5)=0.998411078$, and the
terminal residual is $7.82\times10^{-12}$.  Figure~\ref{fig:num-reference-path}
shows the resulting reference minimum action path.  Since $b(0)=0$ and the
initial displacement and velocity vanish, the zero solution is an equilibrium state
of the deterministic equation, i.e., the equation \eqref{SWE} with $\varepsilon=0$.  The computed path starts from this equilibrium state
and reaches
$u_{48}^*(1,1/2)=0.5$.  
\par

At the large-deviation scale, this surface gives the computed most probable transition of $u_\varepsilon$ for $\varepsilon\ll 1$,  which starts from the equilibrium state (zero solution)  and reaches $0.5$ at $(t,x)=(1,1/2)$.  Since $m_n=n=48$ serves as the finest-grid
numerical reference, Figure~\ref{fig:num-reference-path} is an approximation of 
the exact minimum action path of the continuous equation. 
\par

Thus, the decreasing discrepancies under simultaneous space-time refinement
are consistent with the diagonal convergence asserted by
Theorem~\ref{thm:closed-diagonal} for the nonlinear wave equation in
Example~\ref{ex:wave-fdm-eem}.

\FloatBarrier
\section{Concluding remarks}
\label{sec:concluding-remarks}

In this paper, we  propose a variational framework for proving the convergence of fully discrete, scheme-induced MAMs for SPDEs under simultaneous space-time refinement. We establish an abstract variational convergence theorem and a Green-function criterion. We then verify the criterion for the fully discrete schemes considered for the one-dimensional stochastic wave equation and the multidimensional fourth-order parabolic equation. These results yield convergence of the corresponding discrete rate functions, together with subsequential convergence of asymptotically minimizing controls and the associated minimum action paths.

The convergence results established in this paper are qualitative and do not provide convergence rates. In the present \(\Gamma\)-convergence argument, the control-approximation and Green-function consistency errors are only required to tend to zero, while the terminal correction inherits this qualitative convergence and the lower-bound argument relies on compactness and weak convergence. Deriving quantitative rates would require explicit estimates for these errors and additional analytical tools, since the standard theory of \(\Gamma\)-convergence itself does not provide  convergence rates for variational problems. Establishing such rates for specific nonlinear SPDEs and their fully discrete schemes would be an important direction, which we leave for future work.
%For the wave and fourth-order parabolic examples considered here, the
%criterion gives pointwise convergence of the fully discrete constrained
%minimum values to their continuous counterparts along arbitrary space-time
%diagonals.  The analysis is deliberately confined to a
%scalar one-point observation and a nonvanishing scalar diffusion coefficient.
%Profile constraints, degenerate noises, and systems of SPDEs would require a
%corresponding finite-codimensional controllability condition.  Rates and
%convergence of optimal controls would require stability estimates beyond the
%compactness argument used here.

\section*{Declarations}
\subsection*{Conflict of Interest}
The authors declared that they have no conflict of interest.

\subsection*{Author Contributions}
All authors contributed equally to this work.

\subsection*{Data Availability}
The numerical data generated during the current study are available from the corresponding author upon reasonable request.

\appendix
\section{Verification of the Green-function conditions}
\label{app:kernel-verifications}

This appendix contains the estimates used in the two corollaries of
Section~\ref{subsec:concrete-fully-discrete-schemes}.  The equations,
schemes, and Green functions are those defined in
Examples~\ref{ex:wave-fdm-eem} and~\ref{ex:fourth-order-fdm}.

\subsection{One-dimensional wave equation}
\label{app:wave-kernel-verification}
We now verify \Kone--\Kfive{}.  If
$u_0\in H^\alpha([0,1])$ and $v_0\in H^{\alpha-1}([0,1])$ for some
$\alpha>\frac32$, it follows from
\cite[Propositions 4 and 5]{FDMofSWE} that, for every
$\delta_1\in(0,(\alpha-\frac32)\wedge\frac13)$,
\begin{align*}
\sup_{(t,x)\in[0,T]\times[0,1]}|g_n(t,x)-g(t,x)|
\le Cn^{-\delta_1},
\end{align*}
which proves \Kone{}.

For \Ktwo{}, by \cite[Lemma 1]{FDMofSWE}, for every
$\delta\in(0,\frac23)$, there exists $C:=C(\delta)$ such that for every
$n\ge1$,
\begin{equation}\label{eq:G1Wspace}
\sup_{(t,x)\in[0,T]\times[0,1]}\int_0^1|G_t(x,y)-\mathbf{G}^n(t,x,y)|^2\ud y\le Cn^{-\delta}.
\end{equation}
In view of \cite[Lemma 2.1]{Cohen16}, there exists $C$ independent of $n$ such that for any $0<s<t\le T$,
\begin{align}\label{eq:G1Wtime}
\int_0^1|\mathbf{G}^n(t,x,y)-\mathbf{G}^n(s,x,y)|^2\ud y\le C(t-s).
\end{align}
Combining \eqref{eq:G1Wspace}--\eqref{eq:G1Wtime} and recalling $G_n(t,\lfloor s\rfloor_n;x,z)
=\mathbf{G}^n(t-\lfloor s\rfloor_n,x,z)$ in Example \ref{ex:wave-fdm-eem}, we have
\begin{align*}
\varepsilon_n:=\sup_{(t,x)\in\OT}
{\int_0^t\int_{\mathcal O}
\left|
G_n(t,\lfloor s\rfloor_n;x,z)-G_{t-s}(x,z)
\right|^2 \ud z \ud s}\le C(n^{-\delta}+m_n^{-1}),
\end{align*}
which, together with $m_n\to+\infty$ as $n\to\infty$, yields \Ktwo{}.
For \Kthree{}, it is shown in \cite[p. 320 (20)]{FDMofSWE} that 
\begin{align*}
\sup_{n\ge1}\sup_{(t,x)\in[0,T]\times[0,1]}\int_0^1|\mathbf{G}^n(t,x,y)|^2\ud y<\infty,
\end{align*}
and thus \Kthree{} holds with $\varrho\equiv C$ for some $C>0$.

To verify \Kfour{}, we give the details of the corresponding estimate for the
continuous Dirichlet wave Green function stated above.
Since
$
|e_j(x)-e_j(x')|
\le C\bigl(j|x-x'|\wedge1\bigr),
$
Parseval's identity gives, with
$N=\lfloor |x-x'|^{-1}\rfloor$ when $x\ne x'$,
\begin{align*}
\int_0^1|G_r(x,z)-G_r(x',z)|^2\ud z
&\le
C\sum_{j\ge1}
\left(|x-x'|^2\wedge j^{-2}\right)\\
&\le
C\sum_{j\le N}|x-x'|^2
+C\sum_{j>N}j^{-2}
\le C|x-x'|.
\end{align*}
Similarly,
$|\sin(j\pi r)-\sin(j\pi r')|
\le C(j|r-r'|\wedge1)$ and
$\sup_{j,x}|e_j(x)|<\infty$ imply
\begin{align*}
\int_0^1|G_r(x,z)-G_{r'}(x,z)|^2\ud z
&\le
C\sum_{j\ge1}
\left(|r-r'|^2\wedge j^{-2}\right)
\le C|r-r'|.
\end{align*}
For $|r-r'|>1$, the last estimate follows after enlarging $C$ from the
uniform $L^2$ bound for $G_r(x,\cdot)$.

Now assume without loss of generality that $t\ge t'$.  Using the preceding two estimates yields
\begin{align*}
&\;\int_0^{t'}\int_0^1
\bigl|G_{t-s}(x,z)-G_{t'-s}(x',z)\bigr|^2\ud z\ud s\\
\le &\;
2\int_0^{t'}\int_0^1
|G_{t-s}(x,z)-G_{t-s}(x',z)|^2\ud z\ud s
+2\int_0^{t'}\int_0^1
|G_{t-s}(x',z)-G_{t'-s}(x',z)|^2\ud z\ud s\\
\le &\;
C_T\bigl(|x-x'|+|t-t'|\bigr).
\end{align*}
Thus the desired space-time modulus holds for all
$t,t'\in[0,T]$ and $x,x'\in[0,1]$, and hence \Kfour{} holds.
Proposition~\ref{prop:K5-prototype-operators} gives \Kfive{}, and
Lemma~\ref{lem:grid-control-approximation} verifies the required control
approximation.

\subsection{Fourth-order parabolic equation}
\label{app:fourth-kernel-verification}
We use the mesh, discrete operator, eigenvalues, interpolated
eigenfunctions, and Green function introduced in
Example~\ref{ex:fourth-order-fdm}.  The elementary sine bounds give,
uniformly in $n$ and $\mathbf k\in\mathbb J_n$,
\begin{equation}\label{eq:fourth-discrete-spectrum}
c|\mathbf k|^2\le\mu_{\mathbf k,n}\le C|\mathbf k|^2.
\end{equation}
The discrete sine orthogonality at the grid points and the definition of
the multilinear interpolant give
\begin{equation}\label{eq:fourth-discrete-orthogonality}
\int_{\mathcal O}
\varphi_{\mathbf k}(\kappa_n(z))
\varphi_{\mathbf l}(\kappa_n(z))\ud z
=\mathbf 1_{\{\mathbf k=\mathbf l\}},
\qquad
\sup_{n,\mathbf k,x}|\varphi_{\mathbf k}^n(x)|\le C.
\end{equation}

We now verify \Kone--\Kfive{}.  Since the initial value is zero,
$g_n=g=0$, and hence \Kone{} holds.

The definition of $\chi_n(\mu)$ and \((1+x)^\ell\ge1+\ell x\) $(x\ge-1,\ \ell\ge1)$ give
\begin{equation}\label{eq:fourth-rational-bounds}
|q_n(r,\mu)|
\le \frac{C}{1+r\mu^2},
\qquad
\sup_{0\le t\le T}
\int_0^t
|q_n(t-\lfloor s\rfloor_n,\mu)|^2\ud s
\le C\mu^{-2}.
\end{equation}
We next verify \Ktwo{}.  Fix $K\ge1$ and take $n>K$.  For
$0\le s\le t\le T$, define the low-frequency parts
\begin{align*}
G_n^{\le K}(t,s;x,z)
&:=\sum_{\substack{\mathbf k\in\mathbb J_n\\|\mathbf k|\le K}}
q_n(t-\lfloor s\rfloor_n,\mu_{\mathbf k,n})
\varphi_{\mathbf k}^n(x)\varphi_{\mathbf k}(\kappa_n(z)),\\
G^{\le K}(t,s;x,z)
&:=\sum_{|\mathbf k|\le K}
e^{-(t-s)|\mathbf k|^4}
\varphi_{\mathbf k}(x)\varphi_{\mathbf k}(z),
\end{align*}
and let $G_n^{>K}$ and $G^{>K}$ denote the corresponding sums over
$|\mathbf k|>K$.  Thus
\[
G_n(t,\lfloor s\rfloor_n;x,z)=G_n^{\le K}(t,s;x,z)+G_n^{>K}(t,s;x,z)
\]
and $G_{t-s}(x,z)=G^{\le K}(t,s;x,z)+G^{>K}(t,s;x,z)$.
Set
\begin{align*}
E_{n,K}&:=\sup_{(t,x)\in\OT}\int_0^t\int_{\mathcal O}
|G_n^{\le K}(t,s;x,z)-G^{\le K}(t,s;x,z)|^2\ud z\ud s,\\
R_{n,K}^{\mathrm d}&:=\sup_{(t,x)\in\OT}\int_0^t\int_{\mathcal O}
|G_n^{>K}(t,s;x,z)|^2\ud z\ud s,\\
R_K^{\mathrm c}&:=\sup_{(t,x)\in\OT}\int_0^t\int_{\mathcal O}
|G^{>K}(t,s;x,z)|^2\ud z\ud s.
\end{align*}
The definition of $\varepsilon_n$ in \eqref{eq:green-consistency} and
$|a+b+c|^2\le3(|a|^2+|b|^2+|c|^2)$ give
\begin{equation}\label{eq:fourth-K2-decomposition}
\varepsilon_n\le3E_{n,K}+3R_{n,K}^{\mathrm d}+3R_K^{\mathrm c}.
\end{equation}

For each fixed $K$,  the
eigenvalue convergence $\mu_{\mathbf k,n}\to|\mathbf k|^2$, the uniform convergence of the interpolated eigenfunctions,
and $0\le s-\lfloor s\rfloor_n<\tau_n$ yield
\begin{equation}\label{eq:fourth-low-mode-convergence}
\begin{aligned}
&\max_{\substack{\mathbf k\in\mathbb N^d\\|\mathbf k|\le K}}
\sup_{0\le s\le t\le T}
\left|q_n(t-\lfloor s\rfloor_n,\mu_{\mathbf k,n})
-e^{-(t-s)|\mathbf k|^4}\right|\rightarrow0,\\
&\max_{\substack{\mathbf k\in\mathbb N^d\\|\mathbf k|\le K}}
\left(
\|\varphi_{\mathbf k}^n-\varphi_{\mathbf k}\|_{C(\overline{\mathcal O})}
+\|\varphi_{\mathbf k}\circ\kappa_n-\varphi_{\mathbf k}\|_{L^2(\mathcal O)}
\right)\rightarrow0.
\end{aligned}
\end{equation}
Since the sums contain only finitely many modes,
\eqref{eq:fourth-low-mode-convergence} implies
$E_{n,K}\to0$ as $n\to\infty$.

For the two remainders, from
\eqref{eq:fourth-discrete-orthogonality},
\eqref{eq:fourth-discrete-spectrum}, and
\eqref{eq:fourth-rational-bounds} we obtain
\begin{equation}\label{eq:fourth-K2-tail-bounds}
R_{n,K}^{\mathrm d}
\le C\sum_{|\mathbf k|>K}|\mathbf k|^{-4},
\qquad
R_K^{\mathrm c}
\le C\sum_{|\mathbf k|>K}|\mathbf k|^{-4},
\end{equation}
where the second inequality also uses
$\int_0^t e^{-2(t-s)|\mathbf k|^4}\ud s
\le(2|\mathbf k|^4)^{-1}$.
Combining \eqref{eq:fourth-K2-decomposition}--
\eqref{eq:fourth-K2-tail-bounds} gives
\[
\limsup_{n\to\infty}\varepsilon_n
\le C\sum_{|\mathbf k|>K}|\mathbf k|^{-4}.
\]
Since $\sum_{\mathbf k\in\mathbb N^d}|\mathbf k|^{-4}<\infty$ if and
only if $d<4$, letting $K\to\infty$ proves \Ktwo{} for $d=1,2,3$.

The first bound in \eqref{eq:fourth-rational-bounds}, together with
\eqref{eq:fourth-discrete-orthogonality} and
\eqref{eq:fourth-discrete-spectrum}, also gives
$$
\sup_{n}\sup_{\substack{0\le r<t\le T\\x\in\overline{\mathcal O}}}
\int_{\mathcal O}|G_n(t,r;x,z)|^2\ud z
\le C(t-r)^{-d/4}.
$$
Thus \Kthree{} holds with
$\varrho(r)=Cr^{-d/4}\in L^1(0,T)$.  To verify \Kfour{}, for $K\ge1$ set
$$
G_t^{K}(x,z)
=\sum_{|\mathbf k|\le K}e^{-t|\mathbf k|^4}
\varphi_{\mathbf k}(x)\varphi_{\mathbf k}(z),
\qquad R_t^{K}=G_t-G_t^{K}.
$$
By the orthonormality and uniform boundedness of
$\{\varphi_{\mathbf k}\}_{\mathbf k\in\mathbb N^d}$,
$$
\sup_{0\le t\le T}\sup_{x\in\overline{\mathcal O}}
\int_0^t\int_{\mathcal O}|R_{t-s}^{K}(x,z)|^2\ud z\ud s
\le C\sum_{|\mathbf k|>K}|\mathbf k|^{-4}\rightarrow0.
$$
By symmetry, assume $t\ge t'$.  For fixed $K$, the finite sum $G^K$ is
uniformly continuous, and hence
$$
\lim_{\delta\downarrow0}
\sup_{|t-t'|+|x-x'|\le\delta}
\int_0^{t'}\int_{\mathcal O}
|G_{t-s}^{K}(x,z)-G_{t'-s}^{K}(x',z)|^2\ud z\ud s=0.
$$
Since $G=G^K+R^K$, it follows that
$$
\limsup_{\delta\downarrow0}
\sup_{|t-t'|+|x-x'|\le\delta}\mathcal K_1(t,t',x,x')
\le C\sum_{|\mathbf k|>K}|\mathbf k|^{-4}.
$$
Letting $K\to\infty$ proves \Kfour{} for $d=1,2,3$.

Finally, Proposition~\ref{prop:K5-prototype-operators} gives \Kfive{}.
Together with Lemma~\ref{lem:grid-control-approximation}, this verifies all
kernel and control-approximation hypotheses used in
Corollary~\ref{cor:fourth-diagonal-convergence}.

\section{Fixed-grid LDPs}
\label{app:fixed-grid-LDPs}

The schemes and their continuous interpolations are defined in
Section~\ref{subsec:concrete-fully-discrete-schemes} and are not restated
here.  For each fixed grid, the numerical terminal value is a continuous
function of finitely many Gaussian increments.  We prove the resulting LDPs
and identify the finite-dimensional rates with the scheme-induced control
energies.

\subsection{Wave FDM--EEM scheme}
\label{app:wave-fixed-grid-LDP}
\begin{proposition}[Fixed-grid LDP: one-dimensional wave FDM--EEM]
\label{prop:grid-control-LDP}
Consider the solution $u_\varepsilon^{n,m_n}$ of \eqref{eq:wave-mild-extension} for approximating \eqref{sec4.3SWE}.  Let the settings of Section \ref{subsec:grid-controls-LDP} and Example \ref{ex:wave-fdm-eem} hold.  Then, for every $x_0\in(0,1)$,
$\{u_\varepsilon^{n,m_n}(T,x_0)\}_{\varepsilon>0}$ satisfies an LDP on
$\mathbb R$ with speed $1/\varepsilon$ and good rate function
\begin{align}\label{B.1eq1}
I_{\mathrm w}^{n,m_n}(y)
=
\inf_{\{h_n\in\mathcal X_n:
		\Un(h_n)(T,x_0)=y\}}
\frac12\|h_n\|_{L^2(\mathcal O_T)}^2,
\qquad y\in\mathbb R.
\end{align}
\end{proposition}

\begin{proof}
{
We use the notation introduced in
Example~\ref{ex:wave-fdm-eem}; in particular, $j$ indexes the time
increment and $k$ indexes its spatial component.
Set
$$
Z^\varepsilon
:=
\sqrt{\varepsilon n}\,
(\Delta_0\boldsymbol\beta,\ldots,\Delta_{m_n-1}\boldsymbol\beta)
\in(\mathbb R^{n-1})^{m_n}.
$$
Since the vectors $\Delta_j\boldsymbol\beta$ are independent and have distribution
$N(0,\tau I_{n-1})$, for every
$\theta=(\theta_0,\ldots,\theta_{m_n-1})\in(\mathbb R^{n-1})^{m_n}$,
$$
\lim_{\varepsilon\downarrow0}
\varepsilon\log
\mathbb E\exp\left(\frac{\langle\theta,Z^\varepsilon\rangle}
{\varepsilon}\right)
=\frac{n\tau}{2}\sum_{j=0}^{m_n-1}|\theta_j|^2.
$$
The finite-dimensional G\"artner--Ellis theorem
\cite[Theorem~2.3.6]{Dembo} therefore shows that
$Z^\varepsilon$ satisfies an LDP with speed $1/\varepsilon$ and good
rate function
$$
S^{n,m_n}(z)
=
\frac{1}{2n\tau}\sum_{j=0}^{m_n-1}|z_j|^2.
$$

For $z\in(\mathbb R^{n-1})^{m_n}$, replace
$\sqrt{\varepsilon n}\,\Delta_j\boldsymbol\beta$ in
\eqref{eq:wave-fdm-eem-random-recursion} by $z_j$ and denote the displacement
component of the resulting finite recursion by $U_j^z$.  If
$x_0\in[k/n,(k+1)/n)$, use the terminal interpolation vector
$\ell_{n,x_0}$ defined in Example~\ref{ex:wave-fdm-eem}.  The terminal map
$$
\mathscr F_{n,x_0}(z):=\ell_{n,x_0}^{\top}U_{m_n}^z
$$
is continuous because it is obtained from finitely many continuous
operations.  Hence the contraction principle
\cite[Theorem~4.2.1]{Dembo} yields an LDP for
$u_\varepsilon^{n,m_n}(T,x_0)=\mathscr F_{n,x_0}(Z^\varepsilon)$ with
good rate function
$$
\widehat I_{\mathrm w}^{n,m_n}(y)
=
\inf_{\{z\in(\mathbb R^{n-1})^{m_n}:
\mathscr F_{n,x_0}(z)=y\}}
\frac{1}{2n\tau}\sum_{j=0}^{m_n-1}|z_j|^2.
$$

It remains to identify the above rate function with \eqref{B.1eq1}.  For
$H=(H_0,\ldots,H_{m_n-1})\in(\mathbb R^{n-1})^{m_n}$, write
$$
H_j=(H_{j,1},\ldots,H_{j,n-1})^\top,
\qquad H_{j,k}=[H_j]_k,
$$
so that $H_{j,k}$ denotes the $k$th spatial component of the vector $H_j$
associated with the $j$th time interval.  Substituting
$z_j=\tau H_j$ gives the deterministic recursion
\begin{equation}
\label{eq:wave-terminal-map-control-step}
\begin{aligned}
X_{j+1}^{\tau H}
={}&E_n(\tau)X_j^{\tau H}
+\tau E_n(\tau)F_n(U_j^{\tau H})
+\tau E_n(\tau)\Sigma_n(U_j^{\tau H})H_j,~j=0,\ldots,m_n-1.
\end{aligned}
\end{equation}
Iterating the full-state recursion
\eqref{eq:wave-terminal-map-control-step} from time $0$ to $t_j$ gives
$$
X_j^{\tau H}
=E_n(t_j)X_0^n
+\tau\sum_{r=0}^{j-1}E_n(t_j-t_r)
\bigl(F_n(U_r^{\tau H})+\Sigma_n(U_r^{\tau H})H_r\bigr).
$$
Taking its displacement component yields
\begin{equation}
\label{eq:wave-iterated-controlled-displacement}
\begin{aligned}
U_j^{\tau H}
={}&C_n(t_j)U_0^n+S_n(t_j)V_0^n\\
&+\tau\sum_{r=0}^{j-1}S_n(t_j-t_r)
\bigl(B_n(U_r^{\tau H})+D_n(U_r^{\tau H})H_r\bigr),
\qquad j=0,\ldots,m_n,
\end{aligned}
\end{equation}
where the sum is empty when $j=0$.
This is the deterministic-control counterpart of the iteration from the
full-state recursion to the displacement formula in
\cite[(3.1)--(3.4)]{Cohen16}.

Put $\mathcal C_0=[0,1/n)$ and introduce the effective grid-control subspace
$$
\mathcal X_n^\circ
:=
\bigl\{h\in\mathcal X_n:
h=0\ \text{a.e.~on }(0,T)\times \mathcal C_0\bigr\}.
$$
The correspondence $H\mapsto h_H$ defined by
$$
h_H(s,z)=[H_j]_k=H_{j,k}
\quad\text{on }[t_j,t_{j+1})\times[k/n,(k+1)/n),
\quad k=1,\ldots,n-1,
$$
and $h_H=0$ on $(0,T)\times \mathcal C_0$ is a bijection from
$(\mathbb R^{n-1})^{m_n}$ onto $\mathcal X_n^\circ$.  We next verify directly,
including the
normalization factors, that it produces the controlled FDM--EEM recursion.
We compare \eqref{eq:wave-iterated-controlled-displacement} with the
controlled version of the interpolation \eqref{eq:wave-mild-extension}.
Write $x_k=k/n$ and
$\mathcal C_k=[x_k,x_{k+1})$, $k=0,\ldots,n-1$.  For the grid control $h_H$, put
$$
u_n^H:=\Un(h_H).
$$
For $j=0,\ldots,m_n$, define only the nodal displacement vector
$$
U_j^{n,h_H}
:=
\bigl(u_n^H(t_j,x_1),\ldots,u_n^H(t_j,x_{n-1})\bigr)^\top.
$$
Thus
$$
[U_j^{n,h_H}]_i
=u_n^H(t_j,x_i)
=\Un(h_H)(t_j,x_i),
\qquad i=1,\ldots,n-1.
$$

The normalized eigenpairs of $A_n$ in Example \ref{ex:wave-fdm-eem} satisfy
\begin{align*}
A_nq_\ell&=-(\omega_\ell^n)^2q_\ell,
&\omega_\ell^n&=2n\sin\!\left(\frac{\ell\pi}{2n}\right)
=\ell\pi\sqrt{c_\ell^n},\\
[q_\ell]_k&=n^{-1/2}e_\ell(x_k),
&e_\ell(x)&=\sqrt2\sin(\ell\pi x).
\end{align*}
Hence, spectral calculus gives
\begin{align*}
S_n(t)
&=\sum_{\ell=1}^{n-1}
\frac{\sin(t\omega_\ell^n)}{\omega_\ell^n}q_\ell q_\ell^\top,
&
C_n(t)
&=\sum_{\ell=1}^{n-1}
\cos(t\omega_\ell^n)q_\ell q_\ell^\top.
\end{align*}
For $z\in \mathcal C_k$, we have $\kappa_n(z)=x_k$,
$e_\ell^n(x_i)=e_\ell(x_i)$, and $|\mathcal C_k|=1/n$.
Comparing the preceding expansions with
\eqref{eq:wave-discrete-green}, and differentiating that formula in its
time-lag variable for the second identity, yields
$$
\int_{\mathcal C_k}\mathbf G^n(t,x_i,z)\ud z
=[S_n(t)]_{ik},~
\int_{\mathcal C_k}\partial_t\mathbf G^n(t,x_i,z)\ud z
=[C_n(t)]_{ik},
~ k=1,\ldots,n-1.
$$
Moreover, the definition of $g_n$ in
Example~\ref{ex:wave-fdm-eem} gives
$$
g_n(t_j,x_i)
=\bigl[C_n(t_j)U_0^n+S_n(t_j)V_0^n\bigr]_i.
$$
Evaluate the controlled equation \eqref{eq:fully-discrete-mild} at
$(t_j,x_i)$ with the wave kernel \eqref{eq:wave-discrete-green}.
On $[t_r,t_{r+1})\times \mathcal C_k$ one has
$\lfloor s\rfloor_n=t_r$, $\kappa_n(z)=x_k$, and
$h_H(s,z)=H_{r,k}$.  The contribution of $\mathcal C_0$ vanishes because
$e_\ell(\kappa_n(z))=e_\ell(0)=0$ there.  Splitting the integral into
the space-time cells and using the first kernel identity above gives
\begin{equation}
\label{eq:wave-grid-skeleton-iterated-displacement}
\begin{aligned}
&U_j^{n,h_H}
={}C_n(t_j)U_0^n+S_n(t_j)V_0^n\\
&+\tau\sum_{r=0}^{j-1}S_n(t_j-t_r)
\bigl(B_n(U_r^{n,h_H})+D_n(U_r^{n,h_H})H_r\bigr),
\qquad j=0,\ldots,m_n.
\end{aligned}
\end{equation}

Equations \eqref{eq:wave-iterated-controlled-displacement} and
\eqref{eq:wave-grid-skeleton-iterated-displacement} have the same initial
value and the same Volterra recursion.    Induction therefore gives
$
U_j^{n,h_H}=U_j^{\tau H},
~ j=0,\ldots,m_n,
$
and hence
$$
\Un(h_H)(t_j,x_i)=[U_j^{n,h_H}]_i=[U_j^{\tau H}]_i,
\qquad i=1,\ldots,n-1.
$$
This nodal identification is the controlled analogue of the passage from
the global displacement formula to the continuous interpolation in
\cite[(3.4)--(3.5)]{Cohen16}.
If $x_0\in[x_k,x_{k+1})$, the polygonal spatial interpolation at
$T=t_{m_n}$ consequently gives
\begin{align*}
\Un(h_H)(T,x_0)
&=(k+1-nx_0)[U_{m_n}^{n,h_H}]_k
+(nx_0-k)[U_{m_n}^{n,h_H}]_{k+1}\\
&=(k+1-nx_0)[U_{m_n}^{\tau H}]_k
+(nx_0-k)[U_{m_n}^{\tau H}]_{k+1}\\
&=\ell_{n,x_0}^{\top}U_{m_n}^{\tau H}
=\mathscr F_{n,x_0}(\tau H).
\end{align*}
Moreover,
$
S^{n,m_n}(\tau H)
=\frac{\tau}{2n}\sum_{j=0}^{m_n-1}|H_j|^2
=\frac12\|h_H\|_{L^2(\mathcal O_T)}^2.
$
For any $h\in\mathcal X_n$, define
$$
h^\circ:=h\,\mathbf 1_{(0,T)\times((0,1)\setminus \mathcal C_0)}
\in\mathcal X_n^\circ.
$$
If $z\in \mathcal C_0$, then $\kappa_n(z)=0$ and $e_\ell(0)=0$; hence
\eqref{eq:wave-discrete-green} gives
$\mathbf G^n(r,x,z)=0$.  It follows from the discrete skeleton equation
\eqref{eq:fully-discrete-mild} that
$$
\Un(h)=\Un(h^\circ),
\qquad
\|h^\circ\|_{L^2(\mathcal O_T)}\le\|h\|_{L^2(\mathcal O_T)}.
$$
Since $\mathcal X_n^\circ\subset\mathcal X_n$, we therefore have
$$
\inf_{\{h\in\mathcal X_n:\Un(h)(T,x_0)=y\}}
\frac12\|h\|_{L^2(\mathcal O_T)}^2
=
\inf_{\{h\in\mathcal X_n^\circ:\Un(h)(T,x_0)=y\}}
\frac12\|h\|_{L^2(\mathcal O_T)}^2.
$$
Using the bijection $H\mapsto h_H$, we conclude that
\begin{align*}
\widehat I_{\mathrm w}^{n,m_n}(y)
&=\inf_{\{H\in(\mathbb R^{n-1})^{m_n}:
\mathscr F_{n,x_0}(\tau H)=y\}}
\frac{\tau}{2n}\sum_{j=0}^{m_n-1}|H_j|^2\\
&=\inf_{\{h\in\mathcal X_n^\circ:
\Un(h)(T,x_0)=y\}}
\frac12\|h\|_{L^2(\mathcal O_T)}^2
=I_{\mathrm w}^{n,m_n}(y),
\end{align*}
which proves the claim.}
\end{proof}

\subsection{Fourth-order \texorpdfstring{FDM--implicit Euler scheme}{FDM--implicit Euler scheme}}
\label{app:fourth-fixed-grid-LDP}

\begin{proposition}[Fixed-grid LDP: fourth-order FDM--implicit Euler]
\label{prop:fourth-fixed-grid-LDP}
 Let the settings of Section \ref{subsec:grid-controls-LDP} and Example \ref{ex:fourth-order-fdm} hold.
Fix $n\ge2$, write $m=m_n$ and $\tau=\tau_n$, and let
$u_{\varepsilon,n}$ be the continuous interpolation defined by
\eqref{eq:fourth-mild-extension}.  Then, for every $x_0\in\mathcal O$,
$\{u_{\varepsilon,n}(T,x_0)\}_{\varepsilon>0}$ satisfies an LDP on
$\mathbb R$ with speed $1/\varepsilon$ and good rate function
$$
I_d^{n,m}(y)
=\inf_{\{h_n\in\mathcal X_n:
\Un(h_n)(T,x_0)=y\}}
\frac12\|h_n\|_{L^2(\mathcal O_T)}^2.
$$
\end{proposition}

\begin{proof}
We first record the identity that connects the interpolation
\eqref{eq:fourth-mild-extension} with the nodal scheme.  If
$s\in[t_j,t_{j+1})$ and $i\ge j+1$, then
$q_n(t_i-\lfloor s\rfloor_n,\mu)=\chi_n(\mu)^{i-j}$.
For $x_{\mathbf p}=\Delta x_n\mathbf p$, discrete sine orthogonality and the
spectral representation of $R_n=(I+\tau_nL_{n,d}^2)^{-1}$ therefore give
\begin{equation}\label{eq:fourth-grid-green}
\begin{aligned}
G_n(t_i,t_j;x_{\mathbf p},z)
&=(\Delta x_n)^{-d}[R_n^{\,i-j}]_{\mathbf p\mathbf q},\\
\int_{C_{n,\mathbf q}}G_n(t_i,t_j;x_{\mathbf p},z)\ud z
&=[R_n^{\,i-j}]_{\mathbf p\mathbf q},
\quad z\in C_{n,\mathbf q}.
\end{aligned}
\end{equation}
Substituting \eqref{eq:fourth-grid-green} into
\eqref{eq:fourth-mild-extension} and summing over the space-time cells
shows that its nodal values satisfy
\eqref{eq:fourth-order-implicit-scheme}.  The two constructions have the
same zero initial value and hence coincide at every time node.

We finally identify the fixed-grid rate function in the
present multidimensional setting.  For fixed $n$, we define
\[
Z^\varepsilon
:=
\sqrt{\varepsilon}\,(\Delta x_n)^{-d/2}
(\Delta_0\boldsymbol\beta,\ldots,\Delta_{m-1}\boldsymbol\beta).
\]
By the same argument as in the proof of Proposition  \ref{prop:grid-control-LDP}, $\{Z^\varepsilon\}$ satisfies an LDP on
$(\mathbb R^{(n-1)^d})^m$ with speed $1/\varepsilon$ and good rate
function
$
S_d^{n,m}(z)
=
\frac{(\Delta x_n)^d}{2\tau}
\sum_{j=0}^{m-1}|z_j|^2.
$

For $z=(z_0,\ldots,z_{m-1})$, let
$$
U_{j+1}^z
=
R_n\left[
U_j^z+\tau B_n(U_j^z)+D_n(U_j^z)z_j
\right],
\qquad U_0^z=0.
$$
Using the terminal interpolation operator $\ell_{n,x_0}$ defined in
Example~\ref{ex:fourth-order-fdm}, define
$$
\mathscr F_{d,n,m,x_0}(z):=\ell_{n,x_0}(U_m^z).
$$
This map is continuous.  Since
$u_{\varepsilon,n}(T,x_0)
=\mathscr F_{d,n,m,x_0}(Z^\varepsilon)$, the contraction principle
gives the good rate function
$$
\widehat I_d^{n,m}(y)
=
\inf_{\{z\in(\mathbb R^{(n-1)^d})^m:
\mathscr F_{d,n,m,x_0}(z)=y\}}
S_d^{n,m}(z).
$$

It remains to identify this contraction-principle rate with the
scheme-induced control rate.  Put
$$
\mathcal O_n^\circ
:=\bigcup_{\mathbf k\in\mathbb J_n}C_{n,\mathbf k},
\qquad
\mathcal X_{n,d}^\circ
:=\left\{h\in\mathcal X_n:
h=0\ \text{a.e.~on }(0,T)\times
(\mathcal O\setminus\mathcal O_n^\circ)\right\}.
$$
For
$H=(H_0,\ldots,H_{m-1})\in(\mathbb R^{(n-1)^d})^m$, define
$$
h_H(s,z)=[H_j]_{\mathbf k}
\quad\text{on }[t_j,t_{j+1})\times C_{n,\mathbf k},
\quad \mathbf k\in\mathbb J_n,
$$
and set $h_H=0$ on
$(0,T)\times(\mathcal O\setminus\mathcal O_n^\circ)$.  Then
$H\mapsto h_H$ is a bijection from
$(\mathbb R^{(n-1)^d})^m$ onto $\mathcal X_{n,d}^\circ$.

Set $u_n^H:=\Un(h_H)$ and define its nodal vector by
$$
[U_i^{n,h_H}]_{\mathbf k}
:=u_n^H(t_i,\Delta x_n\mathbf k)
=\Un(h_H)(t_i,\Delta x_n\mathbf k),
\qquad \mathbf k\in\mathbb J_n.
$$
Evaluating the discrete skeleton equation
\eqref{eq:fully-discrete-mild} at the grid points and using
\eqref{eq:fourth-grid-green} gives
\begin{equation}\label{eq:fourth-grid-controlled-skeleton}
U_i^{n,h_H}
=
\tau\sum_{j=0}^{i-1}R_n^{\,i-j}
\left[
B_n(U_j^{n,h_H})
+D_n(U_j^{n,h_H})H_j
\right],
\qquad i=1,\ldots,m.
\end{equation}
On the other hand, substituting $z_j=\tau H_j$ into the deterministic
recursion and iterating it gives the same Volterra recursion
\eqref{eq:fourth-grid-controlled-skeleton}.  Induction therefore yields
$$
U_i^{n,h_H}=U_i^{\tau H},
\qquad i=0,\ldots,m.
$$
{Since $g_n=0$, equations
\eqref{eq:fully-discrete-mild} and \eqref{eq:fourth-discrete-green} show
that $u_n^H(T,\cdot)$ is a linear combination of the functions
$\varphi_{\mathbf k}^n$ and is therefore multilinear on every spatial
cell.  Moreover, since $T=t_m$,
$$
u_n^H(T,\Delta x_n\mathbf k)=[U_m^{n,h_H}]_{\mathbf k},
\qquad \mathbf k\in\mathbb J_n.
$$
Hence the definition of $\ell_{n,x_0}$ gives}
\begin{equation}\label{eq:fourth-terminal-identification}
\Un(h_H)(T,x_0)
=\ell_{n,x_0}(U_m^{n,h_H})
=\ell_{n,x_0}(U_m^{\tau H})
=\mathscr F_{d,n,m,x_0}(\tau H).
\end{equation}
Moreover,
$$
S_d^{n,m}(\tau H)
=\frac{\tau\,(\Delta x_n)^d}{2}\sum_{j=0}^{m-1}|H_j|^2
=\frac12\|h_H\|_{L^2(\mathcal O_T)}^2.
$$

For an arbitrary $h\in\mathcal X_n$, set
$
h^\circ:=h\,\mathbf 1_{(0,T)\times\mathcal O_n^\circ}
\in\mathcal X_{n,d}^\circ.
$
If $z\in\mathcal O\setminus\mathcal O_n^\circ$, the lower-left index of
the cell containing $z$ has at least one zero component.  Hence
$\varphi_{\mathbf k}(\kappa_n(z))=0$ for every
$\mathbf k\in\mathbb J_n$, and \eqref{eq:fourth-discrete-green} gives
$G_n(t,r;x,z)=0$.  It follows that
$$
\Un(h)=\Un(h^\circ),
\qquad
\|h^\circ\|_{L^2(\mathcal O_T)}\le\|h\|_{L^2(\mathcal O_T)}.
$$
Therefore, using the bijection $H\mapsto h_H$,
\begin{align*}
\widehat I_d^{n,m}(y)
&=\inf_{\{H\in(\mathbb R^{(n-1)^d})^m:
\mathscr F_{d,n,m,x_0}(\tau H)=y\}}
\frac{\tau\,(\Delta x_n)^d}{2}\sum_{j=0}^{m-1}|H_j|^2\\
&=\inf_{\{h\in\mathcal X_{n,d}^\circ:
\Un(h)(T,x_0)=y\}}
\frac12\|h\|_{L^2(\mathcal O_T)}^2\\
&=\inf_{\{h\in\mathcal X_n:
\Un(h)(T,x_0)=y\}}
\frac12\|h\|_{L^2(\mathcal O_T)}^2
=I_d^{n,m}(y).
\end{align*}
This proves the fixed-grid LDP and identifies its good rate function as
$I_d^{n,m}$.
\end{proof}

\bibliographystyle{abbrv}
\bibliography{references}

\end{document}